\documentclass[11pt]{article}

\usepackage{color}
\usepackage{amsmath}
\usepackage{amssymb}
\usepackage{latexsym}
\usepackage{comment}
\usepackage{xr}
\def\qed{ \ \vrule width.2cm height.2cm depth0cm\smallskip}

\newcommand{\eps}{\varepsilon}

\newcommand{\ba}{\begin{array}}
\newcommand{\ea}{\end{array}}
\newcommand{\be}{\begin{equation}}
\newcommand{\ee}{\end{equation}}
\newcommand{\bea}{\begin{eqnarray}}
\newcommand{\eea}{\end{eqnarray}}
\newcommand{\beaa}{\begin{eqnarray*}}
\newcommand{\eeaa}{\end{eqnarray*}}

\newcommand{\beal}{\begin{align}}
\newcommand{\eeal}{\end{align}}
\newcommand{\beaal}{\begin{align*}}
\newcommand{\eeaal}{\end{align*}}

\def\dbE{\mathbb{E}}
\def\dbF{\mathbb{F}}

\def\dbL{\mathbb{L}}

\def\dbN{\mathbb{N}}
\def\dbP{\mathbb{P}}
\def\dbR{\mathbb{R}}

\def\a{\alpha}
\def\b{\beta}
\def\g{\gamma}
\def\d{\delta}
\def\e{\varepsilon}

\def\l{\lambda}

\def\si{\sigma}

\def\f{\varphi}
\def\th{\theta}
\def\o{\omega}

\def\G{\Gamma}
\def\D{\Delta}
\def\Th{\Theta}

\def\O{\Omega}

\def\uo{\underline\o}

\def\cA{{\cal A}}

\def\cC{{\cal C}}
\def\cD{{\cal D}}
\def\cE{{\cal E}}
\def\cF{{\cal F}}

\def\cI{{\cal I}}

\def\cL{{\cal L}}

\def\cS{{\cal S}}

\def\cU{{\cal U}}
\def\cV{{\cal V}}

\def\cX{{\cal X}}
\def\cY{{\cal Y}}
\def\cZ{{\cal Z}}

\def\ol{\overline}
\def\ul{\underline}

\def\no{\noindent}

\def\ms{\medskip}

\def\q{\quad}
\def\qq{\qquad}

\def\pa{\partial}
\def\cd{\cdot}
\def\cds{\cdots}

\def\td{\nabla}

\def\qed{ \hfill \vrule width.25cm height.25cm depth0cm\smallskip}

\newcommand{\basa}{\begin{assum}}
\newcommand{\easa}{\end{assum}}

\newcommand{\bas}{\begin{assum}}
\newcommand{\eas}{\end{assum}}

\def\limsup{\mathop{\overline{\rm lim}}}

\def\pa{\partial}

 \def\cd{\cdot}
\def\cds{\cdots}

\def\dis{\displaystyle}

\def\cag{c\`{a}gl\`{a}d~}

\def\1{{\bf 1}}

\newcommand{\lo}{\stackrel{\leftarrow}{\o}}
\newcommand{\luo}{\stackrel{\leftarrow}{\uo}}
\newcommand{\lu}{\stackrel{\leftarrow}{u}}

\def\:{\!:\!}
\def\eqref#1{{\rm(\ref{#1})}}
\def \proof{{\noindent \bf Proof\quad}}

at 9pt

\newtheorem{thm}{Theorem}[section]
\newtheorem{lem}[thm]{Lemma}
\newtheorem{cor}[thm]{Corollary}
\newtheorem{prop}[thm]{Proposition}
\newtheorem{rem}[thm]{Remark}
\newtheorem{eg}[thm]{Example}
\newtheorem{defn}[thm]{Definition}
\newtheorem{assum}[thm]{Assumption}

\numberwithin{equation}{section}

\def\&{and}

\date{}
\begin{document}

\title{\bf Fully nonlinear stochastic and rough PDEs:
Classical and viscosity solutions\footnote{An earlier version of this paper was entitled "Pathwise Viscosity Solutions of Stochastic PDEs and Forward Path-Dependent PDEs".
The authors would like to thank Peter Baxendale and Remigijus Mikulevicius for very helpful discussions. 
}}
\author{Rainer {\sc Buckdahn} 
\footnote{
Universit\'e de Bretagne-Occidentale
and
Shandong University, D\'epartement de Math\'ematiques, Rainer.Buckdahn@univ- brest.fr. Research supported in part by PGMO-Gaspard Monge Program for Optimisation and operational research-Fondation Math\'ematiques Jacques Hadamard.}
\and Christian  {\sc Keller}\footnote{University of Central Florida, Department of Mathematics, christian.keller@ucf.edu.}   
\and Jin {\sc Ma} \footnote{University of Southern California and Shanghai Lixin University of Accounting and Finance, Department of Mathematics, jinma@usc.edu. Research supported in part by NSF grant \#DMS 1106853.}
       \and Jianfeng {\sc Zhang}\footnote{University of Southern California, Department of Mathematics, jianfenz@usc.edu. Research supported in part by NSF grant DMS 1413717.}
}\maketitle

\begin{abstract} 
We study fully nonlinear second-order (forward) stochastic partial differential equations (SPDEs). They  
 can also be viewed as forward path-dependent PDEs (PPDEs) and will be treated  as rough 
PDEs (RPDEs) under a unified framework. 
 We develop first a local theory of classical solutions   and  define then viscosity solutions  through smooth test functions.
Our notion of viscosity solutions
is equivalent to the alternative one using semi-jets. Next, we prove  basic properties such as consistency, stability, and a partial comparison principle in the general setting.  
When the diffusion coefficient is 
 semi-linear (but the drift can be fully nonlinear), 
 we  establish a  complete theory, including global existence and comparison principle.
 Our methodology
 relies heavily on the method of characteristics.
\end{abstract}



\vfill\eject

\tableofcontents

\section{Introduction}
\label{sect-introduction}

 We study the   fully nonlinear second order 
 SPDE
\begin{equation}
\label{SPDE0} 
d u(t,x,\o)   = f(t,x,\o, u,\pa_xu,\pa^2_{xx}u)\,dt+  g(t,x, \o, u, \pa_xu) \circ dB_t, 
\end{equation}
with initial condition $u(0,x,\o) = u_0(x)$, where $(t,x)\in[0,\infty)\times\dbR$, $B$ is a standard Brownian motion defined on a probability space $(\O, \cF, \dbP)$, $f$ and $g$ are $\dbF^B$-progressively measurable random fields, and $\circ$ denotes the Stratonovic integration. Using the pathwise analysis based on
Dupire's  path derivatives (see \cite{Dupire}),  we showed in  Buckdahn, Ma, \& Zhang \cite{BMZ} that  SPDE (\ref{SPDE0}) can  be
rewritten as the (forward) PPDE
\begin{equation}
\label{PPDE0}
\begin{split}
\pa^\o_t u(t,x,\o)  = f(t,x,\o, u,\pa_xu,\pa^2_{xx}u),~ \pa_\o u(t,x,\o) =  g(t,x, \o, u, \pa_xu).
\end{split}
\end{equation}
Here, $\pa^\o_t$ and $\pa_\o$ are  temporal and spatial path derivatives in the spirit of \cite{Dupire}. On the other hand, in the
recent work Keller \& Zhang \cite{KZ}, we also investigated SPDE (\ref{SPDE0}) in the framework of rough path theory initiated by Lyons \cite{Lyons98}.
We showed that Gubinelli's derivative (see \cite{Gubinelli})  for 
``controlled rough paths" is essentially equivalent to Dupire's path derivatives. 
Hence, SPDE \eqref{SPDE0} and  also PPDE (\ref{PPDE0}) can be viewed as the rough PDE 
\begin{equation}
\label{RPDE0}
d u(t,x,\o)   = f(t,x,\o, u,\pa_xu,\pa^2_{xx}u)\,dt +g(t,x, \o, u, \pa_xu) \,d\o_t.
\end{equation}
Here, $\o$ is a geometric rough path corresponding to Stratonovic integration. 
 Our investigation  will be built upon these observations.

SPDE (\ref{SPDE0}), especially when both 
$f$ and $g$ are linear or semilinear, has been studied extensively in the literature. We refer to the well-known reference  Rozovsky \cite{Rozovski} for a fairly complete theory on linear SPDEs, and to Krylov \cite{Krylov} for an $L^p$-theory of linear and some  semi-linear SPDEs.  In the case when SPDE (\ref{SPDE0}) is fully nonlinear, as often encountered in applications such as stochastic control theory and
many other fields (cf.~the lecture notes of Souganidis \cite{Souganidis}), the situation is quite different. In fact, in such a case one can hardly expect (global) classical solutions,  even in the Sobolev sense.
Some weaker forms of solutions will have to come
into play. 

In a series of works, Lions \& Souganidis \cite{LS1, LS2, LS3, LS4} initiated the 
investigation 
 of fully nonlinear
SPDEs 
and proposed  
general schemes for finding  
viscosity solutions, especially when $g = g(\pa_x u)$ or $g= g(u)$. Roughly speaking, they suggested two approaches. The first one is to use the method of {\it stochastic characteristics} (cf.~Kunita \cite{Kunita}) to remove the stochastic integrals of SPDE \eqref{SPDE0}, either by defining test functions along the characteristics, or by  defining  
viscosity solutions via  transformed $\o$-wise (i.e., deterministic)  PDEs. 
The second approach is to approximate the Brownian sample paths 
by smooth functions and define a 
solution as the limit, if  it exists, of the solutions to the approximating equations, which are standard PDEs. 

Since the seminal works \cite{LS1} - \cite{LS4}, there have been many efforts to develop the theory of stochastic viscosity solutions using both approaches. 
Along the lines of stochastic characteristics,  Buckdahn \& Ma \cite{BM1, BM2, BM3}, Buckdahn, Bulla, \& Ma \cite{BBM}, and  Matoussi,  Possamai, \& Sabbagh \cite{MPS} proposed different notions of stochastic viscosity solutions in various 
situations when $g$ is independent of $\pa_x u$; see also 
 Seeger \cite{Seeger16Perron}, who investigated the case
$g=g(\pa_x u)$ in the rough path framework. 
The main challenges in this approach seem to be
the following: 

(i) The stochastic characteristics 
exist only locally for general SPDEs;

(ii) the transformed $\o$-wise PDE 
becomes rather complex and generally falls out of the standard PDE literature  (unless $g$ is linear in $u$ and $\pa_x u$). As a consequence, many fundamental issues, such as a comparison principle for general nonlinear SPDEs, remain open. 

On the other hand, as  rough path theory gradually took shape in the early 2000s, many works emerged using the second approach. 
Publications along this line include Caruana, Friz, \& Oberhauser \cite{CFO}, Friz \& Oberhauser \cite{FrizOberhauser11JDE, FO},  Diehl \& Friz \cite{DF},   Diehl, Friz, \& Oberhauser \cite{DFO14proceedings},  Diehl, Oberhauser, \& Riedel \cite{DOR15SPA}, Diehl, Friz, \& Gassiat \cite{DFG17AMO},  Friz, Gassiat, Lions, \& Souganidis \cite{FGLS},  and  Seeger \cite{Seeger17homog, Seeger18schemes}, to mention a few.
In many of these works 
 $g$ is linear in $u$ and $\pa_x u$, or  $g= g(x, u)$ and $f$  is linear in $\pa^2_{xx} u$ (see \cite{DF}), or  $g=g(\pa_xu)$ (see \cite{FGLS,  Seeger17homog,
Seeger18schemes}).
 With this approach, the uniqueness often comes free, provided the uniqueness for the approximating standard PDEs is established. The main challenge here
 is the existence of the  solution (i.e., the existence of the limit).   Note that the two approaches 
 are often  combined; e.g.,  path-approximation was also used in the first approach  and   stochastic characteristics 
 are often used to prove the existence of the limit.

In this paper, we  take the first approach but in a rough path framework, by taking the advantage that  SPDE (\ref{SPDE0}),
PPDE (\ref{PPDE0}), and RPDE (\ref{RPDE0}) can now be considered as equal. More precisely, we shall investigate the notion  of viscosity solution for the general fully nonlinear RPDE \eqref{RPDE0}
 as an $\o$-wise viscosity solution for  SPDE \eqref{SPDE0}. 
 We utilize  PPDE \eqref{PPDE0} by requiring that  smooth test functions $\f$ satisfy
\bea
\label{pao=g}
\pa_\o\f(t,x) = g(t,x,\f,\pa_x\f).
\eea
Involving the coefficient $g$ in the definition of test functions has already been noted in many previous works; 
e.g., in \cite{BM2, BM3,BMZ} the notion of ``$g$-jets" was introduced to reflect the intrinsic  $g$-dependence of the terms involving  ``path derivatives" in the stochastic Taylor expansion. Another reason that the function $g$ will be treated differently than $f$, 
as we will do throughout the paper, is that,  in general, there is no comparison principle in terms of
the  diffusion term $g$  as far as stochastic
differential equations are concerned. 

Using the rough path language, we  define viscosity solutions directly for RPDE \eqref{RPDE0} as well as PPDE \eqref{PPDE0}  in
a completely local manner in all three variables $(t,x, \o)$. We   prove that our definition is equivalent to an alternative definition through semi-jets, 
and show, by using pathwise characteristics, that RPDE \eqref{RPDE0} can be transformed into a standard PDE (with parameter $\o$) without the $d\o_t$ term, as expected. When $g$ is semilinear (i.e., linear in $\pa_x u$), it is not hard to see that our notion of viscosity solution to RPDE \eqref{RPDE0} is also equivalent to the viscosity solution to the 
 transformed PDE in the standard sense of Crandall, Ishii, \& Lions \cite{CIL}.  In the general case, however, the issue becomes quite subtle due to some intrinsic singularity of the transformed PDE, and thus we are not able to obtain the desired equivalence for viscosity solutions. In fact, 
%
 at this point it is even not  clear to us  how to define viscosity solution for the transformed
 PDE when the initial condition $u_0$ is not smooth. 

Modulo some technical conditions as well as  differences of  language, our definition is very
similar or essentially equivalent to the ones in the literature; e.g., in the case $g=g(\pa_x u)$ and $g = g(u)$,
resp., our definition is essentially equivalent to that introduced in \cite{LS1} and \cite{LS3}, resp. When $f$ does not depend on $\pa^2_{xx} u$ (i.e., in the case of  first order RPDEs), our definition is essentially the same as the one in Gubinelli,  Tindel, \& Torrecilla \cite{GTT},  which also uses the rough path language.
They also introduced the alternative definition through semi-jets but left 
the equivalence of these two definitions open. 

We continue with the next goal of this paper. We  establish most, if not all, of the important properties of viscosity solutions. These include consistency (with classical solutions), stability, and a (partial) comparison principle  between a viscosity semi-solution and a classical semi-solution. We follow the arguments of our previous works on backward PPDEs (e.g., Ekren, Keller, Touzi, \& Zhang \cite{EKTZ} and Ekren, Touzi, \& Zhang \cite{ETZ1, ETZ2}).  
The main difference between the forward case and the backward case is the 
additional requirement \eqref{pao=g}.
This requirement is the main source of subtlety when small perturbations on the test function $\f$ are needed. 
 These subtleties  show up in the proofs of stability and the partial comparison principle  especially in the case when $g$ depends on $u$. With the help of 
local classical solutions for first order RPDEs and higher order pathwise Taylor expansions along the lines of \cite{BMZ}, we 
establish crucial auxiliary results (Lemmas \ref{L:Stability} and \ref{lem-PartCompGeneral}) to obtain the desired perturbations. 
The involved techniques 
 have a certain generic efficiency for dealing with the constraint \eqref{pao=g}. In fact, the equivalence of our definitions (test functions vs.~semi-jets) also follows similar lines.

As in all  studies involving viscosity solutions, the most challenging part is the comparison principle. Our case is no exception. 
An inevitable outcome of the  stochastic characteristics approach is to study a fairly complicated $\o$-wise PDE, which will
(naturally, due to the nature of the transformation) have quadratic growth in $\pa_x u$ and will never satisfy a Lipschitz condition in the variable $u$, except for some trivial linear cases. The latter deficiency in the coefficient clearly violates the technical requirements in the standard viscosity solution literature as it is no longer {\it proper} in the sense of  \cite{CIL}.  
Thus special considerations are 
in order. Our plan of attack is  the following. We  establish first  a comparison principle on small time intervals.
Then we extend our comparison principle to arbitrary duration  by using a combination of uniform {\it a priori} estimates for PDEs and BMO estimates for backward stochastic differential equations
 with quadratic growth. Such a ``cocktail" approach enables us to prove  the comparison principle in the general fully nonlinear case under an extra condition, see \eqref{olu=ulu}. In the case when $g$ is semilinear however, even when $f$ is fully nonlinear (e.g., Hamilton-Jacobi-Bellman type), 
 we  verify 
the extra condition (\ref{olu=ulu}) and  establish a complete theory
 including existence and comparison principle.
Thereby, we extend
 the result of Diehl \& Friz \cite{DF}, which follows the second approach and studies the case when both $g$ and $f$ are semilinear.  
 However,  the verification of (\ref{olu=ulu}) in general cases is a 
 challenging issue and requires further investigation. 

Another main result of the paper 
is  the equivalence between local classical solutions to RPDE~\eqref{RPDE0} and those of the corresponding 
 transformed PDE  in the general fully nonlinear case. We
provide sufficient conditions for the existence of  local classical solutions to this PDE.  
Similar results have been proven in the stochastic setting by Da Prato \& Tubaro \cite{DaPratoTubaro} when $g$ is linear (linear in $u$ and $\pa_x u$).  But to our best knowledge, the results for the general fully nonlinear case 
are  new. 
We emphasize again that our PDE involves some serious singularity issues so that the local existence interval depends on the regularity of the classical solution 
(which in turn depends on the regularity  of $u_0$). Consequently, our results are only valid for classical solutions. 
Whether 
the comparison principle  for viscosity solutions of RPDE \eqref{RPDE0} with such  generality can be shown directly  remains an open problem.

The  paper is organized as follows. In Section \ref{sect-rough}, we review the basic theory of rough paths,
  rough differential equations (RDEs), and the crucial rough Taylor expansions. Some results go beyond the standard literature 
and are proven in Appendix.  In Section \ref{sect-RPDE}, we set up the framework for SPDEs, RPDEs, and PPDEs. In Section \ref{sect-classical}, we introduce the crucial characteristic equations and transform the RPDE into a PDE. We  establish the equivalence of their local classical solutions and provide sufficient conditions for their existence. Sections \ref{sect-viscosity} and \ref{sect-comparison} are devoted to  viscosity solutions in the general case.
In Section \ref{sect-semilinear}, we establish the complete viscosity  theory in the case that $g$ is semilinear.  
Finally, in Section \ref{sect-summary}, we summarize our results.


{\color{black}
To conclude this section, we remark that since the paper already involves rather heavy notations, to focus on the main idea and 
to simplify the presentation, we  make three simplifications throughout the paper:

1) We restrict   our work to a finite time horizon $[0, T]$.  
Our results can be easily extended to the infinite horizon case except that the generic constant $C$ in various estimates will then depend on $T$ and  may explode when $T\to\infty$. 

2) We restrict   our work to  a one dimensional setting. The extension to multidimensional rough paths is non-trivial but is standard in the literature. Most results in the paper hold true in multidimensional settings. We  provide further remarks 
when 
the extension to the multidimensional case is crucial or quite tricky.  In particular, 
 Proposition~\ref{P:CharExistence}  relies on results for multidimensional RDEs. Throughout the paper, we use the notation:
 \begin{equation}
 \label{RT}
 \dbR_T := [0, T]\times \dbR\q\mbox{and}\q \dbR^k_T := [0, T]\times \dbR^k.
 \end{equation}

3) The paper  involves higher order derivatives and  related norms. Although we may track the proofs and figure out what regularity exactly is needed, 
we   use  norms involving all  partial derivatives up to the same order. 
Our estimates  suffice for our purpose but they are not necessarily sharp. 
Also, many results involve higher order regularities,
whose orders are denoted with generic constants $k$ 
without specifying the precise value of $k$. 
}

\section{ Preliminary results from rough path theory}
\label{sect-rough}
We  introduce the framework for rough path theory that is used in this paper.
Only those results that are directly needed  are considered. 
 We mainly follow Keller \& Zhang \cite{KZ} (see Friz \& Hairer \cite{FH} and the references therein
 for the general theory).
 Some less standard results will be proved  in Appendix. 
 
 First, we introduce  general notation. Given  normed spaces $E$ and $V$, put
\begin{align*}
 \dbL^\infty(E; V):=\big\{ u: E\to V:\, \|u\|_\infty 
:= \mbox{$\sup_{x\in E}$} \|u(x)\|_V<\infty\big\}.
\end{align*}
When $V=\dbR$, we  omit  $V$ and  just write  $\dbL^\infty(E)$.  For a constant $\a>0$, set
\begin{align*}
C^\a(E; V):= \big\{u\in \dbL^\infty(E; V):  [u]_\a :=\sup_{x, y\in E, x\neq y}  {\|u(x)-u(y)\|_V \over \|x-y\|_E^\a }<\infty\big\}.
\end{align*}
Given functions $u: [0, T]\to \dbR$ and   $\ul u: [0, T]^2\to \dbR$,  we write the time variable as subscript, i.e., $u_t = u(t)$ and $\ul u_{s,t} = \ul u(s, t)$, and  we define
\begin{equation}
\label{ust} 
u_{s,t} := u_t - u_s,~ s, t\in [0, T],\qq [\ul u]_\a := \sup\limits_{s, t\in [0, T], s\neq t}  {|\ul u(s,t)| \slash |s-t|^\a}.
\end{equation}
Moreover, we shall use $C$ to denote a generic constant in various estimates, which will typically depend on $T$ and possibly on other parameters as well.



Furthermore, we define the standard H\"older spaces  and parabolic H\"older spaces (cf.~Lunardi \cite[Chapter~5]{Lunardi}):
Given $k\in \dbN_0$ and $\beta\in (0,1]$, set
\begin{align*}
 C^{k+\beta}_b(\dbR) &:=\big\{ u:\dbR\to\dbR:\, 
 \| u\|_{C^{k+\beta}_b(\dbR)}<\infty\big\},\\
 C^\beta_b(\dbR_T)&:=\big\{u\in C^0(\dbR_T):\,\| u\|_{C^{\beta}_b(\dbR_T)} <\infty\big\},\\
 C^{2+\beta}_b(\dbR_T)&:=\big\{u\in C^{1,2}(\dbR_T):
\| u\|_{C^{2+\beta}_b(\dbR_T)}  <\infty \big\},
\end{align*}
\begin{align*}
\mbox{where}\q \| u\|_{C^{k+\beta}_b(\dbR)}&:=\mbox{$\sum_{j=0}^k$} \| \partial^j_x u\|_\infty+ [\partial^k u]_\b,\\
 \| u\|_{C^{\beta}_b(\dbR_T)}&:= \| u\|_\infty+ \sup_{t\in [0,T]} [u(t,\cdot)]_\b+
\sup_{x\in\dbR} [u(\cdot,x)]_{\b/2},\\
\| u\|_{C^{2+\beta}_b(\dbR_T)}&:=
\mbox{$ \sum_{j=0}^{\color{black} 1}$} \| \pa^j_x  u\|_\infty+ 
\| \pa^2_{xx} u \|_{C^\b_b(\dbR_T)}+
\| \pa_t u\|_{C^{{\color{black}\b}}_b(\dbR_T)}.
 \end{align*}

\subsection{Rough path differentiation and  integration}
\label{ss-roughintegration}
From now on, unless otherwise stated,  fix two parameters 
$ \a \in ({1/ 3}, {1/ 2}]$ and $\b\in (0, 1]$ satisfying
\bea
\label{ab}
\a(2+\b) > 1.
\eea
Let $\hat \o := (\o, \ul \o)$ be a rough path in the following sense: (i) $\o \in C^\a([0, T])$ and  $\ul \o_{s,t} := {(1/ 2)} |\o_{s, t}|^2$;
 (ii) $\hat \o$ is \textit{truly rough}, i.e.,  there is a set $A$ such that
\begin{align}
\label{trulyrough}
\limsup_{t \downarrow  s}
\frac{\lvert  \omega_{s,t} \rvert}
{\lvert t-s \rvert^{\alpha(1+\beta)}}
=\infty\text{ for all $s\in A$ and $A$ is dense in $[0,T]$.}
\end{align}

\begin{rem}
\label{rem-roughpath}
{\rm (i) The second component $\ul \o$ maps  $[0, T]^2$ to $\dbR$ with $[\ul \o]_{2\a} <\infty$.
Note that $\ul\o_{s,t}$ should not be understood as $\ul\o_t - \ul\o_s$ as in \eqref{ust}. 

(ii) For a general $d$-dimensional rough path, $\ul\o: [0, T]^2\to \dbR^{d\times d}$ has $2\a$-regularity in the sense $[\ul\o]_{2\a}<\infty$, and it satisfies  Chen's relation, i.e., 
$\underline{\omega}_{s,t}-
\underline{\omega}_{s,r}-
\underline{\omega}_{r,t}=
\omega_{s,r}\,\omega^\top_{r,t}$,
$s$, $t$, $r\in [0,T]$,
where $^\top$ denotes the transpose. 
Moreover, $\hat\o$ is called a geometric rough path if  
$\ul\o_{s,t} +\ul \o_{s,t}^\top = \o_{s,t}\,\o_{s,t}^\top$, $s$, $t\in [0, T]$.
In our setting, $(\o,\ul\o)$ is a geometric rough path
and the related integration theory  corresponds to  Stratonovic integration.  

(iii) In standard rough path theory, it is typically not required that $\hat\o$ is truly rough 
as defined in \eqref{trulyrough}.
  But it is convenient for us because under \eqref{trulyrough} the rough path derivatives we  define next will be unique.} 
\end{rem}

Next, we  introduce  path derivatives with respect to our rough path. 
To this end, we introduce spaces of multi-indices
\begin{equation*}
\cV_n := \{ 0, 1\}^n,\q  \|\nu\|:=  \sum_{i=1}^n [ 2\1_{\{\nu_i=0\}}+ \1_{\{\nu_i=1\}}]~\mbox{for}~ \nu = (\nu_1,\cds, \nu_n)\in \cV_n.
\end{equation*}

\begin{defn}
\label{defn-C1}
 (a) Let $u\in C^\a([0, T])$ and $C^0_{\a, \b}([0, T]) := C^{\a\b}([0, T])$.

(i) A first order spatial derivative of $u$ is a  $\pa_\o u \in C^{0}_{\a,\b}([0, T])$ that satisfies
\begin{align}
\label{C1}
u_{s,t}=\partial_\omega u_s\,\omega_{s,t} + R^{1,u}_{s,t},\q s,t\in [0, T], \q\mbox{with}\q [R^{1, u}]_{\a(1+\b)} <\infty.
\end{align}

(ii)  Assume that $\pa_\o u \in C^\a([0, T])$ exists and has a derivative $\pa_\o\pa_\o u$,
 then 
a temporal derivative of $u$  is a   $\pa^\o_t u \in C^{0}_{\a,\b}([0, T])$ that
satisfies 
\begin{equation}
\label{C2}
\begin{split}
u_{s,t}&=\pa^\o_t u_s~ [t-s] +  \pa_\o u_s~ \o_{s,t} + \pa_{\o}\pa_\o u_s~ \ul\o_{s,t} + R^{2,u}_{s,t},\,
s,t\in [0,T],\\
&\qquad\mbox{with}~ [R^{2, u}]_{\a(2+\b)} <\infty.
\end{split}
\end{equation}

(iii) For $\nu\in \cV_n$,  $\cD_\nu u := \pa_{\nu_1} \cds \pa_{\nu_n} u$, where $\pa_{0} := \pa^\o_t$ and $\pa_1 := \pa_\o$.

 (b) For $k\ge 1$, let
 $C^k_{\a,\b}([0,T]):=\{ u\in C^\alpha([0,T]):\, \text{$\cD_\nu u$ exists  $\forall\|\nu\|\le k$\}}$.
\end{defn}

\begin{rem}
\label{rem-C1}
{\rm (i) In the rough path literature, a first order spatial derivative $\pa_\o u$ is typically called  
Gubinelli derivative and the corresponding function $u$ is called a \emph{controlled rough path}. 
 In our case, the path derivatives  defined above are unique due to $\hat\o$  being truly rough 
(see 
\cite[Proposition~6.4]{FH}).

(ii) The derivative $\pa_\o u$ depends on $\o$, but not on $\ul \o$.
The derivative $\pa^\o_t u$ depends on $\ul \o$ as well and should  be denoted 
by $\pa^{\hat\o}_t u$. However, in our 
setting $\ul \o$ is a function of $\o$ and 
thus we can write $\pa^\o_t u$ instead.

(iii) When $\pa_\o u=0$,  it follows from \eqref{C2} and \eqref{ab}
that $u$ is differentiable in $t$ and $\pa^\o_t u = \pa_t u$, 
the standard derivative with respect to $t$.

(iv) In the multidimensional case, $\pa_{\o\o} u \in \dbR^{d\times d}$ could be symmetric if  $u$ is smooth enough (see 
\cite[Remark 3.3]{BMZ}); i.e.,  $\pa_{\o^i}$ and $\pa_{\o^j}$ commute for $1\le i$, $j\le d$. However, typically $\pa^\o_t$ and $\pa_\o$ do not commute, even when $d=1$.}
\end{rem}

\begin{rem}
\label{rem-pat}
{\rm Note that in \eqref{C2} the term $t-s$ is the difference of the identity function $t\mapsto t$, which is Lipschitz continuous.  For all  estimates below, it suffices to assume 
$\pa^\o_t u \in C^{\a(2+\b)-1}([0, T])$. However, to make the estimates more homogeneous, we  only use the H\"{o}lder-$2\a$ regularity of $t$  and thus require $\pa^\o_t u \in C^{\a\b}([0, T])$.  For this same reason,  all our estimates  will  actually hold true if we replace $t$ with a H\"older-$2\a$ continuous path $\zeta\in C^{2\alpha}([0,T])$. To be more precise,  we define 
a path derivative of $u$ with respect to $\zeta$ as a function $\pa^\o_\zeta u\in C^0_{\a,\b}([0,T])$ that satisfies $ [R^{2, u}]_{\a(2+\b)} <\infty$, where 
\begin{equation}
\begin{split}
\label{pazeta}
u_{s,t}&=\pa^\o_\zeta u_s ~\zeta_{s,t} +  \pa_\o u_s ~\o_{s,t} + \pa_{\o\o} u_s ~\ul\o_{s,t} + R^{2,u}_{s,t},
\end{split}
\end{equation}
then  Lebesgue integration  $dt$ should be replaced  with  Young integration $d \zeta_t$.} 
\end{rem}

Next, we equip  $C^k_{\a,\b}([0, T])$ with a norm $\|\cd\|_k$. Given $u\in C^k_{\a,\b}([0, T])$, put
\bea
\label{norm}
\left.\ba{lll}
~ \|u\|_0 := |u_0| +[u]_{\a\b}; \q \|u\|_1 := \|u\|_0 + \|\pa_\o u\|_0 + [R^{1, u}]_{\a(1+\b)};\\
 ~  \|u\|_2 :=\|u\|_1+ \|\pa_\o u\|_1 + \|\pa^\o_t u\|_0 + [R^{2, u}]_{\a(2+\b)};\\
~  \|u\|_k :=\|u\|_{k-1}+ \|\pa_\o u\|_{k-1} + \|\pa^\o_t u\|_{k-2},\q k\ge 3.
\ea\right.
\eea
We emphasize that, besides $k$, the norms depend on $T$, $\o$, $\a$, and $\b$ as well. 
To simplify the notation, we do not indicate these dependencies explicitly. 
In some places we  restrict $u$ to some subinterval $[t_1, t_2]\subset [0, T]$.  
Corresponding spaces $C^k_{\a,\b}([t_1, t_2])$ are defined 
in an obvious way. To not further complicate the notation, 
the corresponding norm  is still denoted by $\|\cd\|_k$.
Note  that, for $u\in C^1_{\a,\b}([0, T])$ and for a constant $C$ depending on $\o$,
\begin{align}
\label{norminfinity}
|u_t|\le |u_0| + |\pa_\o u_0|[\o]_\a t^{\a} + [R^{1,u}]_{\a(1+\b)} t^{\a(1+\b)} \le |u_0| + C\|u\|_1 t^\a. 
\end{align}

Finally, we  define  the rough integral of $u\in C^1_{\a,\b}([0, T])$.  Let $\pi: 0=t_0<\cds < t_n =T$ be a time partition and  $|\pi| := \max_{0\le i\le n-1} |t_{i+1}-t_i|$. 
By  \cite{Gubinelli}, 
\bea
\label{roughintegration}
\int_0^t u_s \,d\o_s := \lim_{|\pi|\to 0} \sum_{i=0}^{n-1} \big[u_{t_i} \o_{t_i\wedge t, ~t_{i+1}\wedge t} + \pa_\o u_{t_i} ~\ul\o_{t_i\wedge t, ~t_{i+1}\wedge t}\big].
\eea
exists and defines the rough integral.
The integration path $U_t := \int_0^t u_s\, d \o_s$ 
belongs to $C^1_{\a, \b}([0, T])$ with $\pa_\o U_t = u_t$ 
and we define $\int_s^t u_r \,d\o_r := U_{s,t}$.

 In this context, we  define  iterated integrals as follows.  For $\nu \in \cV_n$, set
\begin{align}
\label{Ith}
\cI^\nu_{s,t} := \int_s^t \int_s^{t_n} \cds \int_s^{t_2} d_{\nu_1} t_1 \cds \,d_{\nu_n} t_n,\,\mbox{where $d_0 t := dt$, $d_1 t = d\o_t$.}
\end{align}
One can 
check that
$ \cI^\mu_{s,t} = \int_s^t \cI^{(\mu_1,\cds, \mu_n)}_{s,r}\, d_{\mu_{n+1}} r$
for  $\mu = (\mu_1,\cds, \mu_{n+1})\in \cV_{n+1}$.
 In the multidimensional case,
 defining  iterated integrals is  
not trivial.
 Nevertheless, by   \cite[Theorem~2.2.1]{Lyons98}, this can be accomplished
via  uniquely determined  (higher-order) extensions of the
 geometric rough path $\hat\o=(\o,\ul\o)$.

By \eqref{C2} and \eqref{ab}, the following result is obvious and we omit the proof.
\begin{lem}
\label{lem-Ito}
(i) If $u \in C^2_{\a,\b}([0, T])$, then
\bea
\label{Ito}
d u_t = \pa^\o_t u_t \,dt + \pa_\o u_t\, d\o_t.
\eea

(ii) Suppose that $u_t = u_0 + \int_0^t a_s\, ds + \int_0^t \eta_s\, d\o_s$ with $a\in C^0_{\a, \b}([0, T])$ and $\eta \in C^1_{\a,\b}([0, T])$. Then $u \in  C^2_{\a,\b}([0, T])$ with $\pa^\o_t u = a$ and $\pa_\o u = \eta$.  Moreover, $[R^{2,u}]_{\alpha(2+\beta)}\le C(\|a\|_0+\|\eta\|_1)$.
 \end{lem}

Finally, we introduce  backward rough paths. Fix $t_0 \in (0, T]$. Set
\bea
\label{back}
\lo^{t_0}_t := \o_{t_0}-\o_{t_0-t},\q \luo^{t_0}_{s,t}~ :=  ~{1\over 2}  |\lo^{t_0}_{s,t}|^2.
\eea
Then $(\lo^{t_0}, \luo^{t_0})$ is 
a rough path on $[0, t_0]$. Moreover, for $u\in C^1_{\a,\b}([0, t_0])$,
the function $\lu^{t_0}$ defined by
 $\lu^{t_0}_t := u_{t_0-t}$ belongs  to $C^1_{\a,\b}([0, t_0])$ with $\hat \o$ replaced by $(\lo^{t_0}, \luo^{t_0})$
 in Definition~\ref{defn-C1}. In this case, 
$\int_0^{t_0} \lu^{t_0}_s \,d\lo^{t_0}_s = \int_0^{t_0} u_s \,d\o_s$.

\subsection{Rough differential equations}
\label{sect-RDE}
We start with controlled rough paths with parameter $x\in \dbR^d$. They  serve as solutions to RPDEs and coefficients for RDEs and RPDEs. For this purpose, we  have to allow $d>1$ here. Consider 
a function $u: \dbR_T^d \to \dbR$. If, for fixed $x\in \dbR^d$, the mapping
 $t\mapsto u(t,x)$ is a  controlled rough path, 
we use the notations $\pa_\o u$, $\pa^\o_t u$, $\cD_\nu u$ to denote the path derivatives as in 
the previous subsection. For fixed $t$, 
 we  use $\pa_x u$, $\pa^2_{xx} u$, etc., to denote the derivatives 
 of 
 $x \mapsto u(t,x)$  
 with respect to $x$. 
Now, we  introduce the appropriate spaces, extending Definition \ref{defn-C1}.

\begin{defn}
\label{defn-Cxloc}
Let  $[t_1, t_2]\subset [0, T]$, $O \subset \dbR^d$ be  convex,  $u\in C^0([t_1, t_2]\times O)$.

(i) We say  $u\in C^{0,loc}_{\a,\b}([t_1, t_2]\times O)$ if the following holds:
\begin{itemize}
 \item  $x \mapsto u(\cd, x)$ maps $O$ into  $C^0_{\a,\b}([t_1, t_2])$ and
 is continuous under $\|\cd\|_0$. 
 \item $x\mapsto u(t,x)$ are locally H\"{o}lder-$\b$ continuous,  uniformly in
$t\in [t_1, t_2]$.
\end{itemize}

 (ii) We say  $u\in C^{1,loc}_{\a,\b}([t_1, t_2]\times O)$ if the following holds:
 \begin{itemize}
 \item 
 $x\mapsto u(\cd, x)$ maps $O$ into $C^1_{\a,\b}([t_1, t_2])$ and
  is continuous under $\|\cd\|_1$.
  \item 
 $\pa_\o u\in C^{0, loc}_{\a,\b}([t_1, t_2]\times O)$ and 
 $\pa_{x} u\in C^{0, loc}_{\a,\b}([t_1, t_2]\times O; \dbR^d)$, in the sense that each component $\pa_{x_i} u\in C^{0, loc}_{\a,\b}([t_1, t_2]\times O)$, $i=1$, $\ldots$, $d$.
\end{itemize}
 
  (iii) We say  $u\in C^{2,loc}_{\a,\b}([t_1, t_2]\times O)$ if the following holds:
  \begin{itemize}
  \item 
   $x\mapsto u(\cd, x)$ maps  $O$ into $C^2_{\a,\b}([t_1, t_2])$ and is continuous under $\|\cd\|_2$.
   \item 
    $\pa_\o u\in C^{1,loc}_{\a,\b}([t_1, t_2]\times O)$, 
    $\pa_x u \in  C^{1,loc}_{\a,\b}([t_1, t_2]\times O; \dbR^d)$, and
  $\pa_t^\o u\in C^{0,loc}_{\a,\b}([t_1, t_2]\times O)$;
    for all $x\in O$, $[\pa_x R^{1,u}(x)]_{\a(1+\b)} < \infty$. 
\end{itemize}

  (iv) For $k\ge 3$, we say  $u\in C^{k, loc}_{\a,\b}([t_1, t_2]\times O)$ if  $u$,  $\pa_\o u\in C^{k-1, loc}_{\a,\b}([t_1, t_2]\times O)$, $\pa_x u\in C^{k-1, loc}_{\a,\b}([t_1, t_2]\times O;\dbR^d)$, and $\pa^\o_t u  \in C^{k-2, loc}_{\a,\b}([t_1, t_2]\times O)$. 
 \end{defn}
 
\begin{lem}
\label{lem-commute}
(i) Let $u \in C^{2, loc}_{\a,\b}(\dbR_T^d)$. Then $\pa_\o$ and $\pa_x$ commute, i.e., 
\bea
\label{commute1}
\pa_\o \pa_x u = \pa_x \pa_\o u. 
\eea
Assume further that $u\in C^{3, loc}_{\a,\b}(\dbR_T^d)$. Then $\pa^\o_t$ and $\pa_x$ 
commute, i.e.,
\bea
\label{commute2}
\pa^\o_t \pa_x \f = \pa_x \pa^\o_t \f.
\eea

(ii) If $u\in C^{1, loc}_{\a,\b}(\dbR_T^d)$, then, 
for any bounded domain $O\subset \dbR^d$,
\bea
\label{Fubini}
\int_s^t \int_O u(r,x) \,dx \,d\o_r = \int_O  \int_s^t u(r,x) \,d\o_r \,dx.
\eea
\end{lem}

The next result is the crucial chain rule 
(cf.~\cite[Theorem~3.4]{KZ}). 

\begin{lem}
\label{lem-chain}
Assume that $\f\in C^{1, loc}_{\a,\b}(\dbR_T^d)$ and $X\in C^1_{\a,\b}([0, T]; \dbR^d)$. Let  $Y_t := \f(t, X_t)$. Then $Y\in  C^{1}_{\a,\b}([0, T])$ and it holds
\bea
\label{chain1}
\pa_\o Y_t = \pa_\o \f (t, X_t) + \pa_x \f(t, X_t) \cd \pa_\o X_t.
\eea
If  $\f\in C^{2, loc}_{\a,\b}(\dbR_T^d)$, $X \in C^2_{\a,\b}([0, T]; \dbR^d)$, then $Y\in  C^{2}_{\a,\b}([0, T])$ and 
\bea
\label{chain2}
\pa^\o_t Y_t = \pa^\o_t \f (t, X_t) + \pa_x \f(t, X_t) \cd \pa^\o_t X_t.
\eea
\end{lem}

Our study  relies heavily on the following  rough Taylor expansion. The result holds true for multidimensional cases as well and we emphasize that the numbers $\d$ below can be negative. 

\begin{lem}
\label{lem-Taylor}
 Let $u \in C^{k, loc}_{\a,\b}(\dbR_T)$ and  $K\subset\dbR$ be compact.
Then, for every $(t,x)\in \dbR_T$ and $(\delta,h)\in\dbR^2$ with
$t+\delta\in [0,T]$ and $x+h\in K$, we have $|R^{k,u}_{t, x; \d, h}| \le  C(K,x)\, (|\d|^\a+|h|)^{k+\b}$, where
\begin{equation}
\label{Taylor}
 u(t+\d, x+h) = \sum_{m=0}^k 
\sum_{\|\nu\|\le k-m}  
\frac{1}{m!}
\cD_\nu \pa_x^m  u(t,x) ~ h^m ~\cI^\nu_{t, t+\d} + R^{k,u}_{t, x; \d, h}.
\end{equation}
\end{lem}

To study RDEs, uniform properties 
for the functions in $C^{k,loc}_{\a,\b}([t_1, t_2]\times O)$ are needed.  In the next definition, we 
abuse the notation $\|\cd\|_k$ from  \eqref{norm}.

\begin{defn}
\label{defn-Cx}
(i)  We say  $u\in C^{k}_{\a,\b}([t_1, t_2]\times O)
\subset C^{k, loc}_{\a,\b}([t_1, t_2]\times O)$ if
\begin{equation}
\label{normx}
\|u\|_k := \sum_{i=0}^k\sup_{x\in O} \|\pa^i_x u(\cd, x)\|_{k-i} < \infty.
\end{equation}

(ii) For 
solutions to standard PDEs (recall Remark~\ref{rem-C1} (iii)), we use
\bea
\label{Ckb}
 C^{k,0}_{\a,\b}([t_1, t_2]\times O) := \Big\{u\in C^k_{\a,\b}([t_1, t_2]\times O): \pa_\o u = 0\Big\}.
\eea
\end{defn} 
We remark that in (i)   we do not require $\sup_{t\in [t_1,t_2]} [\pa^k_x u(t,\cd)]_\b<\infty$, but restrict ourselves to local
H\"older continuity with respect to $x$ (uniformly in $t$), which suffices
for our rough Taylor expansion above.

Although functions in $C^{k,0}_{\a,\b}(\dbR_T)$ are, in general,
 only at most once differentiable in time, they behave in our rough path
framework as if they were $k$ times
differentiable in time
(cf.~also the 
discussion in 
\cite[Section~13.1]{FH}).

\begin{rem}
\label{rem-Cx}
{\rm (i)  If  $u:[t_1,t_2]\times O\to\dbR$ 
satisfies $\| u\|_{k+1}<\infty$ (as in \eqref{normx}), then
$u(t, x+h) - u(t,x) =  h \int_0^1 \pa_x u (t, x+ \l h)\,d\l$.
Thus the mapping $x \mapsto u(\cd, x)$, $O\to C^k_{\alpha,\beta}([t_1,t_2])$, is continuous under $\|\cd\|_k$ (as defined in \eqref{norm})
 and,   for $\|\nu\|=k$, $\cD_{\nu} u(t,\cdot)$ is  H\"older-$\b$ continuous, uniformly in $t$.
Hence, the continuity required in  the definition of $C^{k, loc}_{\a,\b}(\dbR_T^d)$ is automatic.

(ii) Similarly, if $u\in C^{k+1}_{\a,\b}(\dbR_T^d \times \dbR^{d'})$, then 
$y   \mapsto u(\cd,y)$, $\dbR^{d'}\to C^{k+1}_{\a,\b}(\dbR_T^d)$,
 is continuous under $\|\cd\|_k$ (as defined in \eqref{normx}).} 
\end{rem}

Now, we  study rough differential equations  of the form
\bea
\label{RDE}
u_t = u_0 + \int_0^t f(s, u_s)\, ds + \int_0^t g(s, u_s)\, d\o_s.
\eea

\begin{lem}
\label{lem-RDE}
If $f \in C^{k}_{\a,\b}(\dbR_T)$   and $g \in C^{k+1}_{\a,\b}(\dbR_T)$ for some $k\ge 2$, then RDE \eqref{RDE} has a unique solution $u \in C^{k+2}_{\a,\b}([0, T])$ and 
\bea
\label{RDEest}
\|u-u_0\|_{k+2} \le C(T, \|f\|_k, \|g\|_{k+1}).
\eea
\end{lem}

In the following linear case, we have a representation formula for $u$:
\begin{align}
\label{RDElinear}
u_t = u_0 + \int_0^t [f_0(s) + f_1(s) u_s]\, ds + \int_0^t [g_0(s) + g_1(s)u_s]\, d\o_s.
\end{align}

\begin{lem}
\label{lem-RDElinear}
If $f_0$, $f_1 \in C^k_{\a,\b}([0, T])$  and $g_0$, $g_1\in C^{k+1}_{\a,\b}([0, T])$ for some $k\ge 2$, 
then RDE \eqref{RDElinear} has a unique solution $u\in C^{k+2}_{\a, \b}([0, T])$ given by
\begin{align}
\label{linear-rep}
\dis u_t =  \G_t\big[u_0 + \int_0^t {f_0(s)\over \G_s}\,ds + \int_0^t {g_0(s) \over\G_s}\,d\o_s\big], 
\end{align}
where $\G_t := \exp\big\{\int_0^t f_1(s)\,ds + \int_0^t g_1(s)\,d\o_s\big\}$.
\end{lem}

{\color{black}
\begin{rem}
\label{rem-RDElinear}
{\rm This representation holds  only true in the one dimensional case. For multidimensional  linear RDEs, Keller \& Zhang \cite{KZ} derived a semi-explicit representation formula. Moreover,  note that \eqref{RDElinear} actually does  not satisfy the technical conditions in Lemma \ref{lem-RDE} ($f$ and $g$ are not bounded). But nevertheless, due to its special structure,  RDE \eqref{RDElinear} is wellposed as shown in this lemma.
}
\end{rem}
}

Finally, we extend Lemma \ref{lem-RDE} to RDEs with parameters   of the form
\begin{align}
\label{RDEx}
u(t,x) = u_0(x) + \int_0^t f(s, x, u(s, x)) \,ds +  \int_0^t g(s, x, u(s, x)) \,d\o_s.
\end{align}
\begin{lem} 
\label{lem-RDEx}
Assume that $u_0 \in  C^{k+\b}(\dbR)$, 
$f\in C^{k+1}_{\a,\b}(\dbR_T^2) $, and $g\in   C^{k+1}_{\a,\b}(\dbR_T^2)$ for some $k\ge 3$. 
Then 
$u \in C^{k, {loc}}_{\a,\b}(\dbR_T)$ and $\pa_x u$ solves  
\begin{align}
\label{paxuv}
\pa_x u(t,x) = \pa_x u_0(x) +  \int_0^t \big[\pa_x f(s, x, u(s, x)) + \pa_y f(s, x, u(s, x)) \cd 
\pa_x u(s,x)\big] \,ds \nonumber\\
+  \int_0^t \big[\pa_x g(s, x, u(s, x)) + \pa_y g(s, x, u(s, x)) \cd \pa_x u(s,x)\big]\,d\o_s.
\end{align}
 If all the related derivatives of $u_0$ are bounded (but not necessary $u_0$ itself), then $u - u_0 \in C^{k}_{\a,\b}(\dbR_T)$. 
  If $u_0$ 
  is 
  bounded, then $u\in  C^{k}_{\a,\b}(\dbR_T)$.
\end{lem}

\section{Stochastic PDEs, rough PDEs, and  path-dependent PDEs}
\label{sect-RPDE}
\setcounter{equation}{0}
Our initial goal is to study the 
 fully nonlinear stochastic PDE
\begin{equation}
\label{SPDE}
\begin{split}
du(t,x, B_\cd)&=f(t,x, B_\cd, u, \pa_x u, \pa^2_{xx} u)\, +g(t,x,B_\cd,u,\pa_x u) \circ dB_t.
\end{split}
\end{equation}
Here, $B$ is a standard Brownian motion, $\circ$ denotes the Stratonovich integral, and $f$ and $g$ are $\dbF^B$-progressively measurable. 

To transform \eqref{SPDE} into a rough PDE, we need to introduce some notation. 
Let  $\O_0 := \{\o \in C^0([0, T]): \o_0 = 0\}$, 
 $B$ the canonical process on $\O_0$, i.e., 
$B(\o) := \o$, $\dbP_0$ the Wiener measure, and
\begin{align}
\label{Omega}
 \O := \mbox{$\bigcup_{{1/3 } < \a < {1/ 2}}$} \O_\a,\quad 
  \O_\a := \big\{\o\in\O_0:  
[\o]_\a <\infty ~\mbox{and \eqref{trulyrough} holds}\big\}.
\end{align}
Then $\dbP_0(\O) = 1$ (cf.~\cite[Theorem~6.6]{FH}). Moreover, consider the space
\begin{align}
\label{CO}
\dis \cC(\O) &:=  \bigcup\big\{\cC_{\a,\b}(\O): \text{$\a\in ({1/ 3}, {1/2 })$, $\b\in (0, 1]$, and \eqref{ab} holds}\big\},
\end{align}
where $\cC_{\a,\b}(\O)$ is the set of all $\dbF$-progressively measurable real processes $(u_t)_{t\in [0,T]}$
on $\Omega$
 with  $\dbE^{\dbP_0}[|u\|_{B;1}^2] <\infty$ and with
 $u(\o) \in C^1_{\o;\a,\b}([0, T])$ for all $\o\in\O$. Here $\|\cdot\|_{\o;1}$ and $C^1_{\o;\a,\b}([0, T])$ are defined by  \eqref{norm} and Definition~\ref{defn-C1}, respectively, 
with indication of the dependence on $\o$. 
 
Then, for $u \in \cC(\O)$, we have
\begin{align}
\label{Stratonovic}
\Big(\int_0^t u_s \circ dB_s\Big)(\o) = \int_0^t u_s(\o)\,d\o_s,
\q 0\le t\le T,~\mbox{for $\dbP_0$-a.e.~$\o\in \O$}.
\end{align}
Here, the left hand side is a Stratonovic integral while the right hand side is a rough path integral. In this sense, we may write SPDE \eqref{SPDE} as the RPDE
\begin{align}
\label{RPDE1}
du(t,x, \o)=f(t,x, \o, u, \pa_x u, \pa^2_{xx} u)\,dt+g(t,x,\o,u,\pa_x u)\, d\o_t,\, \o\in \O.
\end{align}

\begin{rem}
\label{rem-SPDE}
{\rm (i) If $u$ is a classical solution of 
\eqref{SPDE} with $g(\cdot,x,B_\cd,u,\pa_x u) \in \cC(\O)$ for all $x\in \dbR$, 
then, by \eqref{Stratonovic}, 
RPDE~\eqref{RPDE1} holds true for $\dbP_0$-a.e.~$\o\in \O$. 

(ii) In the earlier version of this paper (see arXiv:1501.06978v1), we studied  pathwise viscosity solutions of SPDE~\eqref{SPDE} in the a.s.~sense. In this version, we 
 study instead the wellposedness of RPDE \eqref{RPDE1} for fixed $\o$. This is easier and 
 more convenient.
Moreover, the rough path framework
allows us to prove crucial perturbation results such as Lemma~\ref{L:Stability}.

(iii)  If we have obtained a solution (in classical or viscosity sense) $u(\cd, \o)$ of RPDE \eqref{RPDE1} for each $\o$, to go back  to SPDE \eqref{SPDE}, one needs to verify the measurability and integrability of the mapping $\o \mapsto u(\cd,\o)$. 
{\color{black}
To do so, one can, in principle, apply the strategy in \cite[Section~3]{DaPratoTubaro},
which relies on construction of solutions to SDEs via iteration so that adaptedness is
preserved. This strategy can  be applied in our setting and does not require $f$
and $g$ to be continuous in $\omega$. Another possible approach is to follow the argument
in \cite[Section~9.1]{FH}, which}
is in the direction of stability and norm estimates but requires at least
 $g$ to be continuous in $\omega$. Since the paper is already too lengthy, we do not pursue these approaches here in detail.} 
\end{rem}

From now on, we shall fix $(\a, \b)$ and $\o$ as in Subsection \ref{ss-roughintegration} and
omit $\o$ in $f$, $g$, and $u$. To be precise, the goal of this paper is to study the  
RPDE 
\bea
\label{RPDE}
du(t,x)=f(t,x,u, \pa_x u, \pa^2_{xx} u)\,dt+g(t,x,u,\pa_x u) d\o_t
\eea
 with initial condition $u(0,x) = u_0(x)$.  Note that $u(t,x)$ implicitly depends on $\o$. In particular, $\pa^\o_t u$  is different from $\pa_t u$  in the standard PDE literature.  Moreover, by Lemma \ref{lem-Ito}, we may  write \eqref{RPDE} as the path dependent PDE 
\begin{align}
\label{PPDE}
\pa^\o_t u(t,x) = f(t,x,u, \pa_x u, \pa^2_{xx} u),\q \pa_\o u(t,x) = g(t,x,u,\pa_x u).
\end{align}
with initial condition.
The arguments of $f$ and $g$  are implicitly denoted as $f(t,x,y,z,\g)$ and $g(t,x,y,z)$.
Throughout the paper, the following assumptions shall be employed.

\begin{assum}\label{assum-g}
Let $g\in C^{k_0, loc}_{\a,\b}(\dbR_T^3)$ for some 
sufficiently large regularity index $k_0\in\dbN$.

(i) $\pa_y g \in C^{k_0-1}_{\a,\b}(\dbR_T^3)$, $\pa_z g \in C^{k_0-1}_{\a,\b}(\dbR_T^3)$.

(ii) For $i=0$, $\ldots$, 
$k_0$ and $(y,z) \in \dbR^2$,  $\pa^i_x g(\cd,y,z) \in C^{k_0-i}_{\a,\b}(\dbR_T)$ with 
\begin{align*}
\|\pa^i_x g(\cd,y,z) \|_{k_0-i} \le C [1+|y|+|z|].
\end{align*}
\end{assum}
Note  that, for any bounded set $Q\subset \dbR^2$,  
$g \in C^{k_0}_{\a,\b}(\dbR_T\times Q)$.

\begin{assum}
\label{assum-f}
(i) $f$ is nondecreasing in $\g$. 

(ii)  $f\in C^0(\dbR_T^4)$ and $|f(t,x, 0,0,0)|\le K_0$ for all $(t,x)\in \dbR_T$. 

(iii) $f$ is uniformly Lipschitz 
in $(y,z,\g)$ with Lipschitz constant $L_0$.  
\end{assum}

\begin{assum}
\label{assum-u0}
Let $u_0$ be uniformly continuous and
 $\|u_0\|_\infty \le K_0$. 
\end{assum}

We remark that there is no comparison principle  in terms of $g$.  Hence, a smooth approximation of $g$ does not help for our purpose and thus we require $g$ to be smooth. By more careful arguments, we may figure out the precise value of $k_0$, but that would make the paper less readable. In the rest of the paper, we  use $k$ to denote a generic index for regularity, which may vary from line to line. We assume 
always that $k$ is large enough so that we can apply all the results in Section~\ref{sect-rough} freely, and we assume that
the regularity index $k_0$ in Assumption \ref{assum-g}  is large enough so that we  have the desired $k$-regularity in the related results.

\begin{defn}
\label{defn-classical} Let $u\in C^{2,loc}_{\a,\b}(\dbR_T)$. We say $u$ is a classical solution (resp. subsolution, supersolution) of RPDE \eqref{RPDE} if 
\begin{equation}
\label{classical}
\begin{split}
&\qq\qq\qq \pa_\o u(t,x) = g(t,x,u,\pa_x u),\\
&\cL u(t,x) := \pa^\o_t u (t,x) - f(t,x,u, \pa_x u, \pa^2_{xx} u) = ~(\mbox{resp.} \le, \ge ) ~ 0.
\end{split}
\end{equation}
\end{defn}
Note again that there is no comparison principle in terms of $g$. So  the first line in \eqref{classical} is an equality even for sub/super-solutions.

\section{Classical solutions of rough PDEs}
\label{sect-classical}
\setcounter{equation}{0}

We establish wellposedness of classical solutions for RPDE \eqref{RPDE}.
To this end, we have to require that  the coefficients $f$, $g$ and the initial value $u_0$ are 
sufficiently smooth.
For general RPDEs, most results are  valid
 only locally in time. 
 However, this will turn out to be sufficient for our study of viscosity solutions in the next sections. 

\subsection{The characteristic equations}
Our main tool is the   {\em method of characteristics}
 (see Kunita \cite{Kunita} for the stochastic setting). It will be used to get rid of the diffusion term $g$ and to transform the RPDE to a standard PDE. Given $\th:= (x, y,z)\in \dbR^3$, consider the 
 coupled system of RDEs
\bea
\label{E:Char}
\begin{split}
X_t&=x- \int_0^t \pa_z g(s,\Theta_s)\,d\omega_s,\\
Y_t&=y + \int_0^t \big[g(s,\Theta_s)-Z_s\pa_z g(s,\Theta_s)\big]\,d\omega_s,\\
Z_t&=z +\int_0^t \big[\pa_x g(s,\Theta_s)+ Z_s \pa_yg(s,\Theta_s) \big]\,d\omega_s.
\end{split}
\eea
Its solution is denoted by $\Theta_t(\th):=(X_t(\th),Y_t(\th),Z_t(\th))$. Fix $K_0>0$ and put
\begin{align}
\label{Q}
 Q := \dbR \times Q_2,\q Q_2:= \{(y,z)\in\dbR^2: 
{\color{black} \max\{|y|, |z|\}}\le K_0 +  1\} \subset \dbR^2.
\end{align}

\begin{prop}\label{P:CharExistence}
Let Assumption~\ref{assum-g} hold and let $K_0\ge 0$ be a constant. Then there exist constants $\d_0>0$ and $C_0$, depending only on $K_0$ and the  $k_0$-th norm of $g$ 
 (in the sense of Definition~\ref{defn-Cx}~(i)) on $[0, T]\times Q$, such that for all $\th\in Q$, the system  \eqref{E:Char}  has a unique solution $\Theta(\th)$ such that 
\begin{align}
\label{Charest}
\tilde \Th \in C^{k_0}_{\a,\b}([0,\delta_0]\times Q; \dbR^3),\, \|\tilde \Th\|_{k_0}\le C_0,\,
\text{where }\tilde \Th_t(\th) := \Th_t(\th) - \th.  
\end{align}
\end{prop}

\proof Uniqueness follows directly from   
an appropriate multi-dimensional extension of Lemma \ref{lem-RDE} for each $\th\in Q$. 
To prove existence, we note that the main difficulty here is 
that some coefficients in \eqref{E:Char} are not bounded.  To deal with this difficulty,
we introduce,   
for each $N>0$,  
a smooth truncation function $\iota^N:\dbR\to \dbR$ with $\iota^N(z)=z$ for $\lvert z\rvert\le N$ and $\iota^N(z)=0$ for $\lvert z\rvert> N+1$, and  consider $g^N (t, \th) := g(t, x, \iota^N(y), \iota^N(z))$.  Then, by Assumption~\ref{assum-g}, $g^N\in C^{k_0}_{\a,\b}(\dbR_T^3)$.  Next, for  each $\th\in\dbR^3$,  consider the system
\begin{equation*}
\begin{split}
&\mbox{$X_t^N=x- \int_0^t \pa_z g^N (s, \Theta_s^N)\,d\o_s,$}\\
&\mbox{$Y_t^N=y+ \int_0^t \big[g^N(s, \Theta^N_s)- \iota^N(Z^N_s)\pa_z g^N(s, \Theta^N_s)\big]\,d\o_s,$}\\
&\mbox{$ Z^N_t=z+ \int_0^t \big[\pa_x g^N(s, \Theta^N_s)+ \iota^N(Z_s^N)\pa_y g^N(s,\Theta^N_s) \big]\,d\o_s.$}
\end{split}
\end{equation*}
Applying Lemma \ref{lem-RDEx}, but extended to the multidimensional case (using the extended Lemma \ref{lem-RDElinear} as shown in Remark \ref{rem-RDElinear}),
the RDE above
has a unique solution $\Theta^N(\th)=(X^N,Y^N,Z^N)(\th)\in C^{k_0}_{\a,\b}([0,T]; \dbR^3)$ and satisfies \eqref{Charest} with a constant $C_N:= C(N, T, \|g^N\|_{k_0})$.
Now set $N := K_0+1$.  For $(t, \th)\in [0, T]\times Q$,  it follows from \eqref{norminfinity} that  
 \begin{equation*}
 \lvert Y^N_t(\th)\rvert \le K_0 +C_N t^\a,\q \lvert Z^N_t(\th)\rvert \le K_0 +C_N t^\a.
 \end{equation*}
Set $\delta_0:=C_N^{-1/\alpha}\wedge T$. Then, for $t\le \d_0$, 
we have $|Y^N_t(\th)|$ and $|Z^N_t(\th)|\le N$. Thus $g^N(\Th^N_t) = g(\Th^N_t)$. Therefore, $\Th^N$ 
solves the original untruncated equation \eqref{E:Char} on $[0, \d_0]$. 
\qed

Next, we linearize 
system 
\eqref{E:Char}.  To this end, put
\begin{align}
\label{UVW}
U := [\pa_x X, \pa_y X, \pa_z X], \, V := [\pa_x Y, \pa_y Y, \pa_z Y],\, W:= [\pa_x Z, \pa_y Z, \pa_z Z].
\end{align}
Then
\begin{equation}
\label{E:LinChar}
\begin{split}
 U_t &=[1, 0, 0] - \int_0^t \Big[ \pa_{xz} g U_s+\pa_{yz} g V_s + \pa_{zz} g\, W_s\Big](s, \Th_s)\,d\o_s,\\
 V_t &=[0,1,0] +\int_0^t \Big[  \big[\pa_x g - \pa_{xz} g\,Z_s\big] U_s+
\big[\pa_yg-Z_s\pa_{yz}g\big]\,V_s\\
&\qquad\qquad\qquad\qquad\qquad - Z_s\pa_{zz}g\,W_s \Big](s, \Th_s)\,d\omega_s,\\
 W_t&=[0,0,1] +\int_0^t \Big[ \big[ \pa^2_{xx} g+Z_s \pa_{xy} g\big]\,U_s+ [\pa_{xy} g+\pa_{yy} g\,Z_s] V_s\\
& \qq\qq \qq\qq\qq + \big[\pa_{xz} g+Z_s\pa_{yz} g+\pa_y g\big]\,W_s \Bigr](s, \Th_s)\,d\o_s.
\end{split}
\end{equation}
The next result is due to Peter Baxendale. It is a slight generalization of 
\cite[(14) on p.~291]{Kunita} (which corresponds to \eqref{VZU}  below). 

\begin{lem}
\label{L:LinCharLemma}
Let Assumption~\ref{assum-g} hold and let $K_0$, $\d_0$ be as in Proposition \ref{P:CharExistence}.
For every $(t,\th)\in [0,\delta_0]\times Q$ with $\th=(x,y,z)$
 and every $h=(h_1,h_2,h_3)\in\dbR^3$,
\begin{equation}\label{E:LinCharLemma}
\begin{split}
V_t(\th)\cdot h -Z_t U_t(\th)\cdot h=(h_2-z\cdot h_1)
\exp\Big\{
 \int_0^t \pa_y g(s,\Theta_s(\th))\,d\o_s
\Big\}.
\end{split}
\end{equation}
\end{lem}

\proof Fix $\th\in Q$, $h\in \dbR^3$. Put
$\G_t := V_t \cdot h -  Z_t U_t \cdot h$.
By 
Lemma~\ref{lem-chain}, 
\beaa
\pa_\o \G_t  &=&   \Big[  \big[\pa_x g-\pa_{xz}g\,Z_t\big] U_t+\Bigl[\pa_y g-Z_t\pa_{yz} g\Bigr]\,V_t- Z_t\pa_{zz} g\,W_t \Big] \cdot h \\
&&- [\pa_x g + Z_t \pa_y g] U_t \cdot h  + Z_t \Big[ \pa_{xz} g  U_t+\pa_{yz} gV_t + \pa_{zz} g\, W_t\Big]\cdot h\\
&=& \pa_y g V_t \cdot h -Z_t \pa_y g U_t \cdot h  = \pa_y g \G_t.
\eeaa
Clearly, $\pa^\o_t \G_t = 0$ and $\G_0 = h_2 - z h_1$. Then Lemma \ref{lem-RDElinear} yields
 \eqref{E:LinCharLemma}.
\qed

\subsection{RPDEs and PDEs}
\label{sect-PDE}
Our goal is to associate RPDE \eqref{RPDE} with a function $v$ satisfying 
\bea
\label{paov1}
\pa_\o v(t, x) = 0,
\eea
which would imply that $v$ solves a standard PDE.
To  illustrate this idea, let us first derive the PDE for $v$ heuristically. Assume 
that $u$ is a classical solution of RPDE \eqref{RPDE} with sufficient regularity.  Recall \eqref{E:Char}.
We want to find $v$ satisfying \eqref{paov1} and
\bea
\label{E:vYZ} 
&\dis u\big(t,X_t(\theta_t(x))\big)=Y_t(\theta_t(x)),\q \partial_x u\big(t,X_t(\theta_t(x))\big)=Z_t(\theta_t(x)),&\\
&\mbox{where}~ \th_t(x) := (x, v(t,x), \pa_x v(t,x)).&\nonumber
\eea
In fact, recall  \eqref{UVW} and write
\bea
\label{hatTh}
\hat \Phi_t(x) := \Phi(t,\th_t(x)) \q\text{for $\Phi = \Th$, $U$, $V$, $W$.}
\eea
Applying the operator $\pa^\o_t$ on both sides of the first equality of \eqref{E:vYZ}
yields together with Lemma \ref{lem-chain}
\begin{equation*}
\label{patuY}
\begin{split}
0&=\pa^\o_t\Big[u(t,\hat X_t) - \hat Y_t\Big] = \pa^\o_t u (t, \hat X_t) + \pa_x u  (t,\hat X_t) \hat U_t \cd \pa_t \th_t(x)- \hat V_t \cd \pa_t \th_t(x) \\
&= f(t, \hat X_t, u, \pa_x u, \pa^2_{xx} u) - \big[ V_t(\th_t(x)) - Z_t(\th_t(x))  U_t(\th_t(x))\big] \cd \pa_t \th_t.
\end{split}
\end{equation*}
 By Lemma~\ref{L:LinCharLemma} with $h := \pa_t \th_t(x) = [0, \pa_t v(t,x), \pa_{tx} v(t,x)]$ and $z=\partial_x v(t,x)$,
\begin{align*}
& \big[ V_t(\th_t(x)) - Z_t(\th_t(x))  U_t(\th_t(x))\big] \cd \pa_t \th_t\\
& = [h_2 - zh_1] e^{\int_0^t \pa_y g(s, \Th_s(\th_t(x))) \,d\o_s} =\pa_t v(t,x) e^{\int_0^t \pa_y g(s, \Th_s(\th_t(x))) \,d\o_s}.
 \end{align*}
 We emphasize 
  that 
  {\color{black} the variable  $\th_t(x)$ above}
   is fixed 
   {\color{black} when Lemma~\ref{L:LinCharLemma} is applied,}
   while the {\color{black} variable $t$ in} $V_t$ is viewed as the running time. In particular, in the last term above $\Th_s(\th_t(x))$ involves both times $s$ and $t$. Then, by \eqref{patuY}, 
\beaa
\pa_t v(t,x) \exp\big(\int_0^t \pa_y g(s,  \Th_s(\th_t(x)))d\o_s\big) &=& f(t, \hat X_t, u, \pa_x u, \pa^2_{xx} u).
\eeaa
By \eqref{E:vYZ}, $u(t, \hat X_t)$ and $\pa_x u(t, \hat X_t)$ are functions of $(t, \th_t(x))$.  Moreover, by applying the operator $\pa_x$ on both sides of the second equality of \eqref{E:vYZ}, 
\beaa
\pa^2_{xx} u(t, \hat X_t)  \hat U_t \cd \pa_x \th_t (x) = \hat W_t  \cd \pa_x \th_t(x).
\eeaa
Note that $\pa_x \th_t(x) = [1, \pa_x v, \pa^2_{xx} v](t,x)$. Then, provided $\hat U_t \cd \pa_x \th_t (x)\neq 0$,
\beaa
\pa^2_{xx} u(t, \hat X_t) = {\hat W_t  \cd \pa_x \th_t(x) \over \hat U_t \cd \pa_x \th_t (x) } = {\partial_x Z_t(\th_t)+\partial_y Z_t(\th_t)\, \pa_x v +\partial_z Z_t(\th_t)\,\pa^2_{xx} v\over  \partial_x X_t(\th_t)+\partial_y X_t(\th_t)\, \pa_x v +\partial_z X_t(\th_t)\,\pa^2_{xx} v}.
\eeaa
Therefore, formally $v$ should satisfy the PDE
\begin{align}
\label{PDE}
\pa_t v (t,x) = F(t,x, v(t,x), \pa_x v(t,x), \pa^2_{xx} v(t,x)),\q v(0,x) = u_0(x),
\end{align}
where, for $\th=(x,y,z)$,
\begin{equation}
\label{F}
F(t,\th,\g) :=  f\Big(t, \Th_t(\th), {W_t(\th)\cd [1, z, \g] \over  U_t(\th)\cd [1, z, \g]}\Big)
e^{-\int_0^t \pa_y g(s, \Th_s(\th)) d\o_s}.
\end{equation}

Now, we    carry out the  analysis above rigorously.
We  start from PDE \eqref{PDE} and derive the solution for RPDE \eqref{RPDE}. Recall   \eqref{Ckb} and that $k$ is a generic, sufficiently large regularity index that  may vary from line to line. 

\begin{lem}
\label{L:hatS}
Let Assumption \ref{assum-g} hold. Let
 $v\in C^{k,0}_{\a,\b}(\dbR_T)$ for some large $k$. Put $K_0:= \|v\|_\infty \vee \|\pa_x v\|_\infty$.
  Let $\d_0$ be determined by Proposition~\ref{P:CharExistence}.  Then there exists a constant $\d\in (0, \d_0]$ such that the following holds:

(i) For every $(t, x)\in [0, \d] \times \dbR$,  
$\pa_x \hat X_t(x)=U_t(\th_t(x)) \cd \pa_x \th_t(x) \ge {1/2}$.

(ii) For every $t\in [0,\delta]$,  
$\hat{X}_t:\dbR\to\dbR$ is a $C^1$-diffeomorphism and $\hat{S}_t(x):=\hat{X}^{-1}_t(x)$ belongs
to $C^{k, loc}_{\alpha,\b}
([0,\delta]\times\dbR)$ (for a possibly different $k$) and satisfies
\begin{equation}
\label{hatS}
\begin{split}
\hat{S}_t(x)&=x-
\int_0^t \frac{[\hat U_s \cd \pa_t \th_s ](\hat S_s(x))}{(\partial_x\hat{X}_s)\,(\hat{S}_s(x))}\,ds + \int_0^t \frac{\pa_z g(s,\hat\Th_s(\hat{S}_s(x)))}{(\partial_x\hat{X}_s)\,(\hat{S}_s(x))}\,d\o_s.
\end{split}
\end{equation}
\end{lem}

\proof (i) Note that $\th_t(x) \in Q$ for all $(t,x)\in \dbR_T$.  By \eqref{Charest},  it is clear that $U \in C^k_{\a,\b}([0, \d_0]\times Q; \dbR^3)$. 
Recall that, by  Definition~\ref{defn-Cx}~(i),
 the regularity here is uniform in $x$. Thus, together with  the regularity of $v$, we have
\begin{equation}
\label{Ureg}
\begin{split}
&| U_t(\th_t(x))- U_0(\th_0(x)) | \\&\qquad\le |U_t(\th_t(x)) - U_t(\th_0(x))|+|U_t(\th_0(x)) - U_0(\th_0(x))|  \le C t^\alpha.
\end{split}
\end{equation}
Since $\pa_x \th_t(x) = [1, \pa_x v, \pa^2_{xx} v](t,x)$ is bounded, 
\beaa
| U_t(\th_t(x)) \cd \pa_x \th_t(x) - U_0(\th_0(x)) \cd \pa_x  \th_0(x) | \le C t^\alpha.
\eeaa
Note that $U_0(\th_0(x)) \cd \pa_x  \th_0(x) = [1, 0, 0] \cd [1, \pa_x v, \pa^2_{xx} v](t,x) = 1$.  Hence,
 there  exists  a $\delta\le \delta_0$ such that,  for every $(t,x)\in [0,\delta]\times\dbR$,
$\partial_x\hat{X}_t(x) =  U_t(\th_t(x)) \cd \pa_x \th_t(x) \ge 1/2$.

(ii) First, by (i) we see that $\hat X_t$ is one to one for $t\in [0,\d]$. Choose $\iota\in C^\infty_b(\dbR)$ with $\iota(y)=1/y$ for $y\ge 1/4$ and  $\iota(y)=0$ for $y\le 1/5$.
Define functions $\hat{a}$, $\hat{b}:[0,\delta]\times\dbR\to\dbR$ by
\begin{equation}\label{E:hatab}
\begin{split}
\hat{a}(t,x)&:=-\iota\big(\partial_x\hat{X}_t(x)\big)~ \hat U_t(x) \cd \pa_t \th_t(x),\\
\hat{b}(t,x)&:=\iota\big(\partial_x\hat{X}_t(x)\big)~ \pa_z g(t,\hat{\Theta}_t(x)).
\end{split}
\end{equation}
Note that $\hat{a}$, $\hat{b}\in C^{k}_{\a,\b} ([0,\delta]\times\dbR)$.
Then, by  Lemma \ref{lem-RDEx},  the RDE
\begin{align*}
\tilde{S}_t(x)=x+\int_0^t \hat{a}(s,\tilde{S}_s(x))\,ds +\int_0^t \hat{b}(s,\tilde{S}_s(x))\,d\o_s,\quad (t,x)\in [0,\delta]\times\dbR,
\end{align*}
has a unique solution $\tilde{S}\in C^{k,loc}_{\a,\b}([0,\delta]\times\dbR)$. Now,  by 
(i), we see that  $\tilde S$  actually satisfies RDE \eqref{hatS}.

It remains to verify that  $\hat{X}_t\circ\tilde{S}_t=\textrm{id}$, $t\in [0,\delta]$. Indeed, note that 
\beaa
\hat{X}_t\circ\tilde{S}_t = X_t(\th_t( \tilde S_t(x))) \q\mbox{and}\q \pa_\o v=0,\q \pa^\o_t X =0.
\eeaa
 Then, by \eqref{E:Char} and \eqref{hatS}, 
\begin{align*}
&\pa_\o [\hat X_t(\tilde S_t(x))] = \pa_\o X_t(\th_t(\tilde S_t(x)))  + U_t(\th_t(\tilde S_t(x))) \cd \pa_x \th_t( \tilde S_t(x)) \pa_\o \tilde S_t(x)\\
&\q = - \pa_z g(t, \hat \Th_t(\tilde S_t(x))) + \pa_x \hat X_t(\tilde S_t(x))  {\pa_z g(t,\hat\Th_t(\tilde {S}_t(x))) \over (\partial_x\hat{X}_t)\,(\tilde {S}_t(x))} =0;\\
&\pa^\o_t[\hat X_t(\tilde S_t(x))]  =  U_t(\th_t(\tilde S_t(x))) \cd \Big[\pa_t \th_t (\tilde S_t(x)) + \pa_x \th_t( \tilde S_t(x)) \pa^\o_t \tilde S_t(x)\Big]\\
&\q= \hat U_t(\tilde S_t(x)) \cd \pa_t \th_t (\tilde S_t(x))  +\pa_x \hat X_t(\tilde S_t(x)) \Big[- { \hat U_t(\tilde S_t(x)) \cd \pa_t \th_t (\tilde S_t(x))\over (\partial_x\hat{X}_t)\,(\tilde {S}_t(x))}\Big] =0.
\end{align*}
Thus $\hat X_t(\tilde S_t(x))= \hat X_0(\tilde S_0(x)) = x$. This concludes the proof. 
\qed

\begin{thm}
\label{T:EU2RPDE}
 Let Assumption \ref{assum-g} hold and $v$ and $\d$ be as in Lemma \ref{L:hatS}. Assume further that $v$ is a classical solution (resp. subsolution, supersolution) of PDE \eqref{PDE}.  Put $u(t,x) := \hat Y_t \circ \hat X_t^{-1}(x)$. Then $u\in C^k_{\a,\b}([0,\delta]\times\dbR)$ is a classical solution (resp. subsolution, supersolution) of RPDE \eqref{RPDE}.
 \end{thm}

\proof  It is clear that $u\in C^{2, loc}_{\a,\b}([0,\delta]\times\dbR)$. To show the uniform properties in terms of $x$, define  first 
$\check S_t(x) := \hat S_t(x) - x$, $\check Y_t(\th) := Y_t(\th)-y$. Then, by Lemma \ref{lem-RDEx}, $\check S \in C^k_{\a,\b}([0, \d]\times \dbR)$ and $\check Y \in  C^k_{\a,\b}([0, \d]\times \dbR^3)$. Note that
\begin{align*}
u(t,x) &= \hat Y_t(\hat S_t(x)) \\ &= \tilde Y_t\big(\check S_t(x)+x, v_t(\check S_t(x)+x), \pa_x v_t(\check S_t(x)+x)\big) + v_t(\check S_t(x)+x).
\end{align*}
Since $v\in C^{k,0}_{\a,\b}(\dbR_T)$, it is clear that $u\in C^k_{\a,\b}([0, \d]\times \dbR)$.

We  prove only the subsolution case. The other statements can be proved similarly. 
First, note that $u(t,x) = Y_t\big(\th_t(\hat S_t(x))\big)$. Then, denoting $\hat x := \hat S_t(x)$,
\begin{align*}
&\pa_\o u(t,x) = \pa_\o Y_t\big(\th_t(\hat x)\big) + V_t\big( \th_t(\hat x)\big) \cd \pa_x \th_t(\hat x) ~ \pa_\o \hat S_t(x) \\
&= \big[g(t, \hat \Th_t(\hat x) - \hat Z_t(\hat x) \pa_z g(t, \hat \Th_t(\hat x))\big] + V_t\big(\th_t(\hat x)\big) \cd \pa_x \th_t(\hat x) {\pa_z g(t,\hat\Th_t(\hat x)) \over (\partial_x\hat{X}_t)\,(\hat x)}\\
&= g(t, \hat \Th_t(\hat x)) +  {\pa_z g(t,\hat\Th_t(\hat x)) \over (\partial_x\hat{X}_t)\,(\hat x)}\big[\hat V_t(\hat x) -  \hat Z_t(\hat x) \hat U_t(\hat x)\big]  \cd \pa_x \th_t(\hat x). 
\end{align*}
Note that, for $(x, y, z) := \th_t(x) = [x, v, \pa_x v]$ and $h:=  \pa_x \th_t(x) = [1, \pa_x v, \pa^2_{xx} v]$, we have $h_2 - z h_1 = \pa_x v - \pa_x v =0$. Then, by Lemma \ref{L:LinCharLemma}, we have
\bea
\label{VZU}
\Big[\hat V_t(\hat x) -  \hat Z_t(\hat x) \hat U_t(\hat x)\Big]  \cd \pa_x \th_t(\hat x) =0,
\eea
and thus 
\bea
\label{paou}
\pa_\o u(t,x) &=& g(t, \hat \Th_t(\hat x)). 
\eea
Similarly, note that $\pa_t \th_t(x) = [0, \pa_t v, \pa_{tx} v]$, 
\begin{align*}
&\pa^\o_t u(t,x) =  V_t\big(\th_t(\hat x)\big) \cd \Big[\pa_t \th_t(\hat x) + \pa_x \th_t(\hat x) ~ \pa^\o_t \hat S_t(x) \Big]\\
&=\hat V_t\big(\hat x) \cd \pa_t \th_t(\hat x)  +  \hat Z_t(\hat x)(\partial_x\hat{X}_t)\,(\hat x)\Big[- \frac{\hat U_t(\hat x) \cd \pa_t \th_t (\hat x)}{(\partial_x\hat{X}_t)\,(\hat x)}\Big]\\
&= \big[\hat V_t(\hat x) -  \hat Z_t(\hat x) \hat U_t(\hat x)\big] \cd \pa_t \th_t(\hat x) = \pa_t v(t, \hat x) \exp\Big(\int_0^t \pa_y g(s, \hat \Th_s(\hat x)) d\o_s\Big).
\end{align*}
Since $v$ is a classical subsolution of \eqref{PDE}-\eqref{F}, the definition of $F$ yields
\begin{align}
\label{patu}
\pa^\o_t u(t,x) \le    f\Big(t, \hat \Th_t(\hat x), {\hat W_t(\hat x)\cd  [1, \pa_x v(t, \hat x), \pa^2_{xx} v(t, \hat x)] \over  \hat U_t(\hat x)\cd [1, \pa_x v(t, \hat x), \pa^2_{xx} v(t, \hat x)]}\Big).
\end{align}

Now, we  identify the functions inside $g$ and $f$ in \eqref{paou} and \eqref{patu}. First, by definition
\bea
\label{xu}
\hat X_t (\hat S_t(x)) = x\text{ and }\hat Y_t (\hat S_t(x)) = u(t,x).
\eea
Next, differentiating \eqref{xu} with respect to $x$, we have
\beaa
&1 = \pa_x \big[X_t( \th_t(\hat x))\big] = U_t( \th_t(\hat x)) \cd \pa_x \th_t(\hat x) ~ \pa_x \hat S_t(x),&\\
&\pa_x u(t,x) =\pa_x \big[Y_t(\th_t(\hat S_t(x)))\big] = V_t(\th_t(\hat x)) \cd \pa_x \th_t(\hat x) ~ \pa_x \hat S_t(x).&
\eeaa
Thus, by \eqref{VZU},
\begin{equation}
\label{paxu}
\pa_x u(t,x) - \hat Z_t (\hat x) = \Big[\hat V_t(\hat x) -  \hat Z_t (\hat x) \hat U_t(t, \hat x)\Big] \cd \pa_x \th(t, \hat x) =0.
\end{equation}
Moreover,
\begin{equation}
\label{paxxu}
\begin{split}
&\pa^2_{xx} u(t,x) = \pa_x[\pa_x u(t,x)] = \pa_x \Big[ Z_t (\th_t(\hat S_t(x)))\Big]\\ 
&\q=  W_t (\th_t(\hat x)) \cd \pa_x \th_t(\hat x) ~ \pa_x \hat S_t(x)= {\hat W_t(\hat x) \cd [1, \pa_x v(t, \hat x),  \pa^2_{xx} v(t, \hat x)]  \over \hat U_t(\hat x) \cd [1, \pa_x v(t, \hat x),  \pa^2_{xx} v(t,\hat x)] }.
\end{split}
\end{equation}
Plugging \eqref{xu}-\eqref{paxxu} into \eqref{paou}-\eqref{patu}, we see that $u$ satisfies the desired subsolution properties.
\qed

  Now, we proceed in the opposite direction, namely deriving $v$ from $u$. Assume that
  $u\in C^k_{\a,\b}(\dbR_T)$ for some large $k$ and define $K_0 :=\|u\|_\infty \vee \|\pa_x u\|_\infty$.  Let $Q_2$ and $Q$ be as in \eqref{Q} and $\d_0$  as in Proposition \ref{P:CharExistence}.   For any fixed $(t,x) \in [0, \d_0]\times \dbR$, consider the  mapping 
  \begin{equation}
  \label{yzmapping}
  (y, z) 
   \mapsto 
   [  Y - u(t, X), Z - \pa_x u(t, X)]\,
   (t,x,y,z)
  \end{equation}
  from $Q_2$ to $\dbR^2$. The Jacobi matrix of  this mapping 
is given by
  \beaa
  J(t,x,y,z) := \left[ \ba{lll} 
  \pa_y Y - \pa_x u(t, X) \pa_y X & \pa_y Z - \pa^2_{xx} u(t, X) \pa_y X \\ 
   \pa_z Y - \pa_x u(t, X) \pa_z X & \pa_z Z - \pa^2_{xx} u(t, X) \pa_z X
   \ea\right]\,{\color{black}(t,x,y,z)}.
   \eeaa
   Note that $\det(J(0,x,y,z)) = 1$. 
    {\color{black} Thus, noting also that $\pa_x u$ and $\pa^2_{xx} u$ are bounded, one
    can see, similarly to \eqref{Ureg}, that there exists a $\d\le \d_0$ such that  $\det(J(t,x,y,z)) \ge {1/2}$ for all $(t,x,y,z)\in [0, \d]\times Q$. }
    This implies that the mapping \eqref{yzmapping} is one to one and the inverse mapping has sufficient regularity. Denote by $R(t,x)$ the range of the mapping \eqref{yzmapping}. Then 
   \beaa
   R(0,x)  = \{(y-u(0,x), z- \pa_x u(0,x)): (y,z) \in Q_2\} \supset  \dbR \times (-1, 1).
   \eeaa 
   Thus,  by \eqref{Ureg} and the boundedness of  $\pa_x u, \pa^2_{xx} u$ again, 
   and by choosing a smaller $\d$ if necessary, we may assume that
    $(0, 0) \in R(t,x)$ for all $(t,x) \in [0, \d] \times \dbR$. Therefore, for any $(t,x) \in [0, \d] \times \dbR$, there exists a unique $(v(t,x), w(t,x))\in Q_2$ such that, denoting $\tilde \th_t(x):= (x, v(t,x), w(t,x))$,
   \begin{align}
 \label{YZvw}
Y_t(\tilde \th_t(x)) = u(t, X_t(\tilde \th_t(x))),\qq Z_t(\tilde \th_t(x)) = \pa_x u(t, X_t(\tilde \th_t(x)))
   \end{align}
   Differentiating the first equation in \eqref{YZvw} with respect to $x$ and applying the second, we obtain
   \beaa
   0 &=& \pa_x \Big[Y_t(\tilde \th_t(x)) - u(t, X_t(\tilde \th_t(x))\Big] \\
   &=& \Big[V_t(\tilde \th_t(x)) - \pa_x u (t, X_t(\tilde \th_t(x))) U_t(\tilde \th_t(x))\Big] \cd \pa_x \tilde \th_t(x) \\
   &=& \Big[V_t(\tilde \th_t(x)) - Z_t U_t( \tilde \th_t(x))\Big]  \cd [1, \pa_x v(t,x), \pa_x w(t,x)] \\
   &=& \big[\pa_x v(t,x) - w(t,x)\big] \exp\big(\int_0^t \pa_y g(s, \Th_s(\tilde \th_t(x)))d\o_s\big),
  \eeaa
  where the last equality  holds true thanks to Lemma \ref{L:LinCharLemma}. Then $w(t,x) = \pa_x v(t,x)$ and thus  \eqref{E:vYZ}  holds. In particular,  we may use the notation $\th_t(x)$ in \eqref{E:vYZ}   again to replace $\tilde \th_t(x)$.

   We  verify now that $v$  indeed satisfies PDE \eqref{PDE}.
  
    \begin{thm}
  \label{thm-utov} 
  Let Assumption \ref{assum-g} hold, let $u\in C^k_{\a,\b}(\dbR_T)$ for some large $k$, and
  let $\d$ and $v$  be determined as above. Assume further 
 that
 $u$ be a classical solution (resp.~subsolution, supersolution) of RPDE \eqref{RPDE}.  Then, for a possibly smaller $\d>0$, 
  we have 
  $ U_t(\th_t(x)) \cd \pa_x \th_t(x) \ge {1/ 2}$ for all $(t, x)\in [0, \d] \times 
  \dbR$
   and $v \in C^{k,0}_{\a,\b}([0,\d]\times 
   \dbR)$ is a classical solution (resp.~subsolution, supersolution) of PDE \eqref{PDE} on 
   $[0, \d] \times \dbR$.
  \end{thm}
\proof  The regularity of $v$ is straightforward. We  prove only the case that $u$ is a classical subsolution. The other cases can be proved similarly.

 Recall the notations in \eqref{hatTh}. Differentiating the first equality of  \eqref{E:vYZ}  with respect to $\o$ and applying the second equality, we obtain
   \begin{align*}
   0 &= \pa_\o \Big[Y_t(\th_t(x)) - u(t, X_t(\th_t(x)))\Big] 
   = \pa_\o Y_t(\th_t(x)) + \hat V_t  \cd  \pa_\o \th_t(x)\\&\qquad  - \pa_\o u(t, \hat X_t ) - \pa_x u (t, \hat X_t )\big[\pa_\o X(t, \th_t(x)) +  \hat U_t \cd  \pa_\o \th_t(x)\big].
   \end{align*}
   By \eqref{classical} and   \eqref{E:vYZ}, 
 $
   \pa_\o u(t, \hat X_t ) = g(t, \hat X_t, u(t, \hat X_t), \pa_x u(t, \hat X_t)) =  g(t, \hat \Th_t).
 $
   Then, by \eqref{E:Char} and Lemma \ref{L:LinCharLemma},
    \begin{align*}
 0&=   [g(t, \hat \Th_t) - \hat Z_t \pa_z g(t, \hat \Th_t) ] + \hat V_t \cd  \pa_\o \th_t(x) - g(t, \hat \Th_t) - \hat Z_t [-  \pa_z g(t, \hat \Th_t)]\\&\qquad - \hat Z_t \hat U_t \cd  \pa_\o \th_t(x)\\
   &= \big[\hat V_t  -  \hat Z_t \hat U_t \big] \cd  [0,  \pa_\o v(t,x), \pa_\o\pa_x v(t,x) ]  = \pa_\o v(t,x) e^{\int_0^t \pa_y g(s, \Th_s(\th_t(x))) d\o_s}.
  \end{align*}
 Thus $\pa_\o v(t,x) = 0$
 and Lemma \ref{L:hatS} can be applied. In particular, for a possibly smaller  $\d>0$, 
  $ U_t(\th_t(x)) \cd \pa_x \th_t(x) \ge {1/ 2}$ for all $(t, x)\in [0, \d] \times \dbR$.
  
  Finally, following exactly the same arguments as for deriving  \eqref{PDE}, one can complete the proof  that $v$ is a classical subsolution of PDE \eqref{PDE}.
  \qed
 
 {\color{black}

\begin{rem}
\label{rem-siz}
{\rm We shall investigate the case with semilinear $g$ in details in Section~\ref{sect-semilinear} below. 
Here we consider another special case:
\bea
\label{glineary}
g =  \si(z).
\eea
which has received strong attention in the literature.  Let $\si'$ and $\si''$ denote the first and second order derivatives of $\si$, respectively.
 In this case the characteristic equations \eqref{E:Char} becomes
\begin{align*}
X_t=x- \int_0^t \si'(Z_s)\,d\omega_s,\, Y_t=y + \int_0^t \big[\si(Z_s) -Z_s\si'(Z_s)\big]\,d\omega_s,\,
Z_t=z,
\end{align*}
which has  the explicit global solution
\bea
\label{E:Charlineary}
X_t = x - \si'(z) \o_t,\q Y_t = y + [\si(z) - z\si'(z)]\o_t,\q  Z_t = z.
\eea
Moreover, in this case \eqref{F} becomes
\begin{equation*}
 F(t,x,y,z,\g) := f\big(t,  x- \si'(z)\o_t, ~ y+ [\si(z)-z\si'(z)]\o_t, ~z, ~ {\g  \over  1 - \si''(z)\o_t  \g}\big).
 \end{equation*}
}
\end{rem}
}

\subsection{Local wellposedness of  PDE~\eqref{PDE}}
To study the wellposedness of PDE \eqref{PDE} and hence that of RPDE \eqref{RPDE}, we first establish a PDE result. Let $K_0 >0$ and, similar to \eqref{Q}, consider
\bea
\label{Q3}
Q_3:= \{(y,z,\g) \in \dbR^3: {\color{black} \max\{|y|, |z|, |\g|\}}\le K_0+1\}.
\eea

\begin{lem}
\label{lem-PDElocal}
Let $k\ge 2$ and $\d_0>0$.

(i) Suppose that $u_0 \in C^{k+1+\b}_b(\dbR)$ with $|u_0|$, $|\pa_x u_0|$, $|\pa^2_{xx}u_0|\le K_0$.

(ii)  Suppose that $F\in  C^{k+1}_{\a,\b}([0, \d_0]\times \dbR\times Q_3)$  and $\pa_\g F \ge c_0>0$.

\no Then there exists a constant $\d\le \d_0$, depending on $K_0$, $c_0$, and the norm $\|F\|_2$ 
on $[0, \d_0]\times \dbR\times Q_3$, such that PDE \eqref{PDE} has a classical solution $v\in C^{k,0}_{\a,\b}([0, \d]\times \dbR)$ on $[0, \d] \times\dbR$.
\end{lem}

\proof It suffices to prove  $v\in C_b^{2+\b}([0,\delta]\times\dbR)$. The further regularity of $v$ when $k\ge 2$ follows from  standard bootstrap arguments
 (cf.~Gilbarg \& Trudinger \cite[Lemma~17.16]{GT}) together with Remark~\ref{rem-Cx}.
Since the proof is very similar to that of  Lunardi \cite[Theorem~8.5.4]{Lunardi}, which considers a similar boundary-value problem,
we shall present only the main ideas for the more involved existence part of the lemma. 
The first step is to linearize our equation and set up an appropriate fixed point problem. To this end, let $\delta>0$ and define an operator $A: C_b^{2+\b}([0,\delta]\times \dbR)\to C_b^{\b}([0,\delta]\times \dbR)$ by
\begin{equation}
\label{Aw}
\begin{split}
 (A v)\,(t,x)&:= \partial_y F(\hat\th_0(x))\, v(t,x)+\partial_z F(\hat\th_0(x))\, \partial_x v(t,x)\\
 &\qquad +\partial_\gamma F(\hat\th_0(x))\, \partial^2_{xx}v (t,x),
\end{split}
\end{equation}
where  $\hat\th_0(x) := \big(0,x,u_0(x),\partial_x u_0(x),\partial_{xx}^2u_0(x)\big)$.
Next, define 
\begin{align}
\label{B1}
B_1 := 
\{ v\in   C_b^{2+\b}([0,\delta]\times\dbR): v(0,\cd)=u_0,  \|v- u_0\|_{C_b^{2+\b}} \le 1
\}.
\end{align}
Now given $v\in B_1$, consider the solution $w$ of the  linear PDE  
\begin{align}
\label{PDEvw}
\partial_t w =Aw +[F(t,x,v,\partial_x v,\partial^2_{xx} v)-Av]\text{ on $[0,\d]\times\dbR$}
\end{align}
with $w(0,\cdot)=u_0$.
Following the arguments in  Lunardi \cite[Theorem~8.5.4]{Lunardi}, when $\d>0$ is small enough,  
PDE \eqref{PDEvw}  has a unique solution $w \in B_1$. This defines a mapping $\G(v) := w$ for $v\in B_1$. Moreover, when $\d>0$ is small enough, $\G$ is a contraction mapping, and hence there exists a unique fixed point $v\in B_1$. Then $v=w$ and, by \eqref{PDEvw},  $v$ solves \eqref{PDE} on $[0, \d]\times\dbR$.
\qed

We now turn back to PDE \eqref{PDE}-\eqref{F} and RPDE \eqref{RPDE}.
\begin{thm}
\label{thm-localclassical}
Let Assumption \ref{assum-g} hold and let  $k\ge 2$, $\d_0>0$. 

(i) Suppose that $u_0 \in C^{k+1+\b}_b(\dbR)$ with $|u_0|$, $|\pa_x u_0|, |\pa^2_{xx}u_0|\le K_0$.

(ii) Suppose that  $f\in C^{ k+1}_{\a,\b}([0, \d_0]\times \dbR\times Q_3)$  and $\pa_\g f \ge c_0>0$.

\no Then there exists a constant $\d\le \d_0$, depending on $K_0$, $c_0$, the regularity of $f$  on $[0, \d_0]\times Q_3$, and the regularity of $g$ on $[0, \d_0]\times Q$,  such that PDE \eqref{PDE}-\eqref{F} has a classical solution $v\in C^{k,0}_{\a,\b}([0, \d]\times \dbR)$ on $[0, \d]\times\dbR$, 
and consequently, for a possibly smaller $\d>0$, RPDE \eqref{RPDE} has a classical solution $u\in C^k_{\a,\b}([0,\d]\times \dbR)$.
\end{thm}

\proof Recall \eqref{F}. By the uniform regularity of $\Th$ in Proposition \ref{P:CharExistence}, one can verify straightforwardly that, for $\d>0$ small enough, $F$ satisfies the conditions in  Lemma  \ref{lem-PDElocal} (ii). Then, by Lemma  \ref{lem-PDElocal},  PDE \eqref{PDE}-\eqref{F}  has a classical solution $v\in B_1$ for a possibly smaller $\d$. Finally, it follows from Theorem \ref{T:EU2RPDE} that RPDE \eqref{RPDE} has a local classical solution.
\qed

\subsection{The first order case}
\label{sect-1st}
We consider the case  $f$ being of first order, i.e.,
\bea
\label{f1st}
f =  f(t,\th) = f(t,x,y,z).
\eea
This case is completely degenerate in terms of $\g$. It is not covered by
Theorem~\ref{thm-localclassical}.  However, in this case  PDE \eqref{PDE}-\eqref{F} is  also of first order,
i.e.,
\begin{equation}
\label{F1st}
F(t,\th) := f\big(t, \Th_t(\th)\big)\,
{\color{black}\exp\big(-\int_0^t \pa_y g(s, \Th_s(\th)) d\o_s\big)}.
\end{equation}
When $f$ is smooth, so is $F$.  Thus
we can modify the characteristic equations \eqref{E:Char} to  solve PDE \eqref{PDE}-\eqref{F1st} explicitly. 
Put  $\tilde \Th = (\tilde X, \tilde Y, \tilde Z)$ and consider  
\bea
\label{Char1st}
\begin{split}
\tilde X_t&=x- \int_0^t \pa_z F(s,\tilde \Theta_s)\,ds,\\ 
\tilde Y_t&=y + \int_0^t \big[F(s,\tilde \Theta_s)-\tilde Z_s\pa_z F(s,\Theta_s)\big]\,ds,\\
\tilde Z_t&=z +\int_0^t \big[\pa_x F(s, \tilde \Theta_s)+ \tilde Z_s \pa_yF(s, \tilde \Theta_s) \big]\,ds.
\end{split}
\eea
Similar to \eqref{E:vYZ}, let $\tilde v$ be determined (locally in time) by
\begin{align*}
 \tilde v\big(t, \tilde X_t(\tilde \theta_t(x))\big)=\tilde Y_t( \tilde \theta_t(x)),~ \partial_x\tilde v\big(t,\tilde X_t( \tilde \theta_t(x))\big)=\tilde Z_t(\tilde \theta_t(x)), 
\end{align*}
where $\tilde \th_t(x) := (x, \tilde v(t,x), \pa_x \tilde v(t,x))$.
Then one can see that \eqref{paov1} should be replaced with 
$\pa_t \tilde v =0$, and thus $\tilde v(t,x) = u_0(x)$.
By similar (actually 
easier)
 arguments as in previous subsections, one can  prove
the following statement.

\begin{thm}
\label{thm-1st}
Let Assumption \ref{assum-g} hold, $f$ take the form \eqref{f1st} with 
$f\in C^{k+1}_{\a,\b}([0, T]\times Q)$,  and
 $u_0 \in C^{k+1+\b}_b(\dbR)$ for some large $k$ with $|u_0|$, $|\pa_x u_0|$, $|\pa^2_{xx}u_0|\le K_0$. Then there is a constant $\d>0$ such that the following holds:

(i) The system of ODEs \eqref{Char1st} is wellposed on $[0,\d]$ for all 
$\th\in Q$.

(ii) For each $t\in [0, \d]$,  the mapping $x\in \dbR\mapsto \tilde X_t(x, u_0(x), \pa_x u_0(x))\in \dbR$
is invertible 
and thus possesses an inverse function, 
to be denoted by $\tilde S_t$.

(iii) The map $v$ defined by $v(t,x) := \tilde Y_t(\tilde \th_t(\tilde S_t(x)))$ belongs to $C^{k,0}_{\a,\b}([0, \d]\times \dbR)$ and is a classical solution to PDE \eqref{PDE}-\eqref{F1st}. Consequently  RPDE \eqref{RPDE}-\eqref{f1st} has a classical solution $u\in C^{k}_{\a,\b}([0, \d]\times \dbR)$.
\end{thm}

\section{Viscosity solutions of rough PDEs: definitions and basic properties}
\label{sect-viscosity}
\setcounter{equation}{0}
We introduce a notion of viscosity solution for RPDE \eqref{RPDE} and study its basic properties.  For any $(t_0,x_0)\in (0,T]\times\dbR$ and $\d \in (0, t_0)$, define
\beaa
 D_\d(t_0,x_0) := [t_0-\d, t_0]\times O_\d(x_0) := [t_0-\d, t_0]\times \{x\in \dbR:  |x-x_0|\le \d\}.
 \eeaa

\subsection{The definition}

For $u\in C(\dbR_T)$ and $(t_0,x_0)\in (0,T]\times\dbR$,   put  
\begin{equation}
\label{cAg}
\begin{split}
 \cA_g^0 u(t_0,x_0; \d) &:=
\Big\{\varphi\in C^{2}_{\alpha,\b}(D_\d(t_0,x_0)):  \varphi(t_0,x_0)=u(t_0,x_0) \\
& \qq\qq \mbox{and}\q \partial_\omega \varphi=g(\cdot,\varphi,\partial_x \varphi)  ~\mbox{on} ~ D_\d(t_0,x_0)\Big\},\\
 \overline{\mathcal{A}}_g u(t_0,x_0) &:= \bigcup_{0<\d\le t_0}\Big\{ \f \in \cA_g^0 u(t_0,x_0; \d):   \f \le u ~\mbox{on} ~ D_\delta(t_0, x_0)\Big\},\\
 \underline{\mathcal{A}}_g u(t_0,x_0) &:= \bigcup_{0<\d\le t_0}\Big\{ \f \in \cA_g^0 u(t_0,x_0;\d):   \f \ge u ~\mbox{on} ~ D_\delta(t_0, x_0)\Big\}.
\end{split}
\end{equation}

\begin{defn}\label{D:2RPDEviscosity}
Let $u\in C(\dbR_T)$ 
and recall the operator $\cL$ in \eqref{classical}.

(i) We say $u$ is a \emph{viscosity supersolution} (resp.~\emph{subsolution}) of  RPDE \eqref{RPDE} if, for every $(t_0,x_0)\in (0,T]\times\dbR$ and $\varphi\in\overline{\mathcal{A}}_g u(t_0,x_0)$ (resp.~$\underline{\mathcal{A}}_g u(t_0,x_0)$), we have
$\cL\f(t_0,x_0)\ge\text{ (resp.~$\le$) } 0.
$

(ii) We say $u$ is a \emph{viscosity solution} of RPDE \eqref{RPDE} if it is both a viscosity supersolution and a viscosity subsolution of \eqref{RPDE}.
\end{defn}
We remark that it is possible to consider semi-continuous viscosity solutions as in the standard literature. However, for simplicity, in this paper we restrict ourselves to continuous 
solutions
 only.

\begin{prop}[Consistency]
\label{prop-consistency}
Let Assumptions~\ref{assum-g} and \ref{assum-f} hold and let
$u\in C^{2}_{\alpha,\b}(\dbR_T)$. Then $u$ is a classical subsolution (resp.~classical supersolution) of RPDE \eqref{RPDE} if and only if it is a viscosity subsolution (resp.~viscosity supersolution) of \eqref{RPDE}.
\end{prop}

\proof We prove only the subsolution case. The supersolution case can be proved similarly.

First, assume that $u$ is a viscosity subsolution. By choosing $u$ itself as a test function, we can immediately infer that $u$ is a classical subsolution.

Next, assume that $u$ is a classical subsolution. Let $(t,x)\in (0,T]\times\dbR$ and $\varphi\in\underline{\mathcal{A}}_g u(t,x)$ with corresponding $\delta_0\in (0,t]$. Then, at $(t,x)$, 
\begin{equation}\label{E:Cons1}
\begin{split}
&u-\varphi=0,\quad\partial_x [u-\varphi]=0,\qquad  \partial^2_{xx} [u-\varphi]\le 0, \\
&~~\pa_\o [u-\f] =0;\q \pa_{x\o} [u-\f] = c~ \pa^2_{xx} [u-\f], 
\end{split}
\end{equation}
where $c:=\pa_z g(t,x, u, \pa_x u)$.
For any $(\delta,h)\in [0,\delta_0]\times O_{\delta_0}(x)$,
 by Lemma \ref{lem-Taylor},
\begin{equation}\label{E:ConsTaylor}
\begin{split}
 0\ge &[u-\varphi](t-\delta,x+h)\\ &=-\partial_t^\o [u-\f](t,x)~\delta+\frac{1}{2}\partial^2_{xx} [u-\f](t,x)
 |h- c~\omega_{t-\delta,t}|^2 +R^{2,u-\f}_{\delta,h},
\end{split}
\end{equation}
where $R^{2,u-\f}_{\delta,h}=O((\delta^\alpha+h)^{2+\beta})$.
Fix a number $\delta_1\in (0,\delta_0]$ such that, for every $\delta\in (0,\delta_1]$,
we have $|c ~\omega_{t-\delta,t}|<\delta_0$.
 From now on, let $\delta\in (0,\delta_1]$. 
Setting $h:= c~\omega_{t-\delta,t}$ in \eqref{E:ConsTaylor} yields
\begin{align*}
-\partial_t^\o [u-\f](t,x)\,\delta\le -R^{2, u-\f}_{\delta, h}\le C(\delta^\alpha+|c\omega_{t-\delta,t}|)^{2+\beta}\le C\,\delta^{\alpha(2+\beta)}.
\end{align*}
Recall \eqref{ab}. Dividing the inequality above by
$\delta$ and sending $\delta$  to $0$, we have  
$\partial_t^\o u(t,x)\ge \partial_t^\o \varphi(t,x)$.
By Assumption~\ref{assum-f}~(i) and by \eqref{E:Cons1},
\begin{align*}
\Bigl[\partial_t^\o \varphi-f(\cdot,\varphi,\partial_x\varphi,
\partial^2_{xx}\varphi)\Bigr](t,x)&\le
\Bigl[\partial_t^\o u-f(\cdot,u,\partial_x u,
\partial^2_{xx} u)\Bigr](t,x)\le 0,
\end{align*}
i.e., $u$ is a viscosity subsolution at $(t,x)$. 
\qed

\subsection{Equivalent definition through semi-jets} 
As in the standard PDE case (cf.~Crandall, Ishii, \& Lions \cite{CIL}), viscosity solutions can also be defined through semi-jets. To see this, we first note that, for $\f\in \cA^0_gu(t_0, x_0;\d)$, our   second order Taylor expansion (Lemma \ref{lem-Taylor}) yields
\begin{align*}
\f(t,x) &= \f(t_0, x_0) + \pa^\o_t \f(t_0, x_0) (t-t_0) + \pa_\o \f(t_0,x_0) \o_{t_0,t} \\
& + \pa_x \f(t_0,x_0) (x-x_0)+ \pa^2_{\o\o}\f(t_0,x_0) \underline \o_{t_0,t} + {1\over 2} \pa^2_{xx} \f(t_0,x_0) |x-x_0|^2 
\\&+ \pa_{x\o} \f(t_0, x_0) \o_{t_0, t} (x-x_0) + R(t,x),
\end{align*}
where $(t,x) \in D_\d(t_0, x_0)$.  Since $\pa_\o \f (t,x)= g(t,x, \f, \pa_x\f)$, we have
\begin{equation}
\label{paxo}
\begin{split}
\pa_{x\o} \f &= \pa_x g + \pa_y g \pa_x \f + \pa_z g \pa^2_{xx} \f,\\
\pa_{\o\o} \f &= \pa_\o g + \pa_y g \pa_\o \f + \pa_z g \pa_{\o x} \f \\&= \pa_\o g + \pa_y g g + \pa_z g [\pa_x g + \pa_y g \pa_x \f + \pa_z g \pa^2_{xx}\f ].
\end{split}
\end{equation}
Motivated by this, we define semijets as follows. Given  $u\in C(\dbR_T)$, $(t_0,x_0)\in (0,T]\times\dbR$, and $(a,z,\gamma)\in \dbR^3$, put
\begin{align*}
&\psi^{ a, z,\g}_{g,u, t_0, x_0} (t,x) := y + a[t-t_0] +b\, \o_{t_0,t} + z [x-x_0] \\ 
&\qquad\qquad\qquad+c\,\underline \o_{t_0,t} + {1\over 2} \g |x-x_0|^2 + q\, \o_{t_0, t} [x-x_0],\q\mbox{where}\\
&y:=u(t_0,x_0),~ b:=g(t_0,x_0,y,z),~ q:= [\partial_x g+\partial_y g z +\partial_z g \gamma](t_0,x_0,y,z), \\
&c:=[\partial_\omega g+\partial_y g g+ \partial_z g  (\partial_x g+\partial_y g z +\partial_z g \gamma)] (t_0,x_0,y,z).
\end{align*}
We then define the 
$g$-superjet $\overline{\mathcal{J}}_g u(t_0,x_0)$  and the 
$g$-subjet $\underline{\mathcal{J}}_g u(t_0,x_0)$ by
\begin{equation}
\label{jet}
\begin{split}
 \overline{\mathcal{J}}_g u(t_0,x_0) &:=\bigcup_{0<\d\le t} \Big\{(a,z,\g) \in \dbR^3:    \psi^{a, z,\g}_{g, u, t_0, x_0} \le u ~\mbox{on}~ D_\d(t_0,x_0)\Big\},\\
\dis \underline{\mathcal{J}}_g u(t_0,x_0) &:=\bigcup_{0<\d\le t} \Big\{(a,z,\g) \in \dbR^3:    \psi^{a, z,\g}_{g, u, t_0, x_0}  \ge u~\mbox{on}~ D_\d(t_0,x_0)\Big\}.
\end{split}
\end{equation}
Notice that $\pa_\o  \psi^{a, z,\g}_{g, u, t_0, x_0}  = g(\cd, \psi^{a, z,\g}_{g, u, t_0, x_0}, \pa_x \psi^{a, z,\g}_{g, u, t_0, x_0})$ holds true only at $(t_0,x_0)$, but not in $D_\d(t_0,x_0)$, so in general $ \psi^{a, z,\g}_{g, u, t_0, x_0} \notin \cA^0_g u(t_0,x_0; \d)$. Nevertheless, we still have the following equivalence.

\begin{prop}\label{P:jetEq}
Let Assumptions~\ref{assum-g} and \ref{assum-f} be in force and $u\in C(\dbR_T)$.
Then $u$ is a viscosity supersolution (resp.~subsolution) of \eqref{RPDE} at $(t_0,x_0)\in (0,T]\times\dbR$ 
if and only if, for every  $(a,z,\gamma)\in\overline{\mathcal{J}}_g u(t_0,x_0)$ (resp.~$\underline{\mathcal{J}}_g u(t_0,x_0)$),
\begin{align}\label{E:jetCriterion} 
a-f(t_0,x_0,u(t_0,x_0),z,\gamma)\ge 0\text{ (resp.~$\le 0$).}
\end{align}
\end{prop}

\proof We prove only the supersolution case. The subsolution case can be proved similarly.

First, we prove the if part.  Assume that \eqref{E:jetCriterion} holds for every $(a,z,\gamma)\in\overline{\mathcal{J}}_g u(t_0,x_0)$.
Let $\varphi\in\overline{\mathcal{A}}_g u(t_0,x_0)$.
Then there exists a $\delta_0\in (0,t_0\wedge 1]$ such that,
whenever $0\le \d\le \d_0$, $|h|\le \d_0$,
\begin{align*}
u(t_0-\d,x_0+h) - u(t_0,x_0)\ge \varphi(t_0-\d,x_0+h) -\varphi(t_0,x_0),\q  
\end{align*}
By Lemma ~\ref{lem-Taylor}, there exists a $C>0$ such that
\begin{align*}
&\varphi(t_0-\d,x_0+h) -\varphi(t_0,x_0)\ge \Big[ - \partial_t^\o \varphi~ \d +\partial_x \varphi ~h +\frac{1}{2} \partial_{xx}^2\varphi~|h|^2\\
&+\partial_\omega \varphi~\omega_{t_0,t}+\partial^2_{\omega\omega}\varphi ~ \underline{\omega}_{t_0,t} +\partial_{x\omega} \varphi~h~\omega_{t_0,t}\Big](t_0,x_0)- C\big[\d^\alpha+|h|\big]^{(2+\beta)}.
\end{align*}
For any  $\eps>0$, by \eqref{ab},  there exists a $\delta_\eps\in (0,\delta_0)$,
such that, for every $0\le \d\le \d_\e$, $|h|\le \d_\e$,
\begin{align*}
&u(t_0-\d,x_0+h) - u(t_0,x_0)\ge \Big[ - \partial_t^\o \varphi~ \d +\partial_x \varphi ~h +\frac{1}{2} \partial_{xx}^2\varphi~|h|^2\\
&+\partial_\omega \varphi~\omega_{t_0,t}+\partial^2_{\omega\omega}\varphi ~ \underline{\omega}_{t_0,t} +\partial_{x\omega} \varphi~h~\omega_{t_0,t}\Big](t_0,x_0) -\eps~\d -{1\over 2}\eps~|h|^2,
\end{align*}
By \eqref{paxo}, the above  inequality  implies 
$(\partial_t^\o \varphi+\eps,\partial_x \varphi,\partial_{xx}^2 \varphi
-\eps)(t_0,x_0)
\in \overline{\mathcal{J}}_g u(t_0,x_0)$. Thus, by \eqref{E:jetCriterion},  
$[\partial_t^\o \varphi + \e -f(\cdot,\varphi,\partial_x\varphi,
\partial_{xx}^2\varphi-\e) ](t_0,x_0)\ge 0$.
Sending $\e\to 0$ yields 
 $\cL\f(t_0,x_0)\ge 0$, i.e.,
$u$ is a viscosity supersol. at $(t_0, x_0)$.

Next, we prove the only if part. Assume that $u$ is a viscosity supersolution at  $(t_0,x_0)\in (0,T]\times\dbR$. Let $(a,z,\gamma)\in\overline{\mathcal{J}}_g u(t_0,x_0)$ and consider the 
RPDE
\begin{align*}
\varphi(t,x) &= u(t_0,x_0)+[a+\e]\,[t-t_0]
+ z\,[x-x_0]+\frac{1}{2}[\gamma-\e]\,\lvert x- x_0\rvert^2\\
&\qquad - \int_t^{t_0} g(\cdot,\varphi,\partial_x\varphi)(s,x)\,d\o_s.
\end{align*}
By Theorem \ref{thm-1st}, the  RPDE above has a classical solution
 $ \f\in C^{2}_{\a,\b}(D_\d(t_0,x_0))$
for some   $\delta\in (0,\delta_0]$. It is clear that $\f \in \cA^0_g u(t_0,x_0;\d)$. Moreover, by using
our Taylor expansion (Lemma \ref{lem-Taylor}), one may easily verify that
\begin{align*}
\f = \psi^{a+\e, z, \g-\e}_{g, u, t_0,x_0}  + R~\mbox{on}~ D_\d(t_0,x_0),
\end{align*}
 where $|R(t,x)|\le C[|t-t_0|^\a + |x-x_0|]^{2+\b}$.
Then, by choosing $\d>0$ small enough, we have
$\f \le \psi^{a, z, \g}_{g, u, t_0,x_0} \le u~\mbox{on} ~ D_\d(t_0,x_0)),
$
where the second inequality is due to the assumption  $(a,z,\gamma)\in\overline{\mathcal{J}}_g u(t_0,x_0)$. This implies $\f \in \overline\cA_g u(t_0,x_0)$. Thus
$0  \le \cL\f(t_0, x_0) =  a+\e - f(t_0, x_0, u(t_0,x_0), z, \g-\e)$.
Sending $\e\to 0$  yields \eqref{E:jetCriterion}.
\qed

\begin{rem}\label{R:jet}
{\rm By Proposition~\ref{P:jetEq} and its proof, we can see that,
depending on the regularity order $k_0$ of $g$ as specified in Assumption \ref{assum-g}, it is equivalent to use test functions of 
 class $C^{k}_{\a,\b}(D_\d(t_0,x_0))$ for any $k$ between $2$ and $k_0$. This is crucial  for  Theorem \ref{T:Stability} below.
 }
\end{rem}

\subsection{Change of variables formula}

Let $\lambda \in C([0, T])$  and $n\ge 2$ be an even integer. 
For any $u:\dbR_T\to\dbR$, define
\begin{align}\label{E:tildeu}
\tilde u(t,x):={e^{\eta_t}\over 1+x^n} u(t,x),\mbox{ where } \eta_t := \int_0^t \lambda_s\,ds.
\end{align}
If $u\in C^{2}_{\a,\b}(\dbR_T)$, then 
\begin{align}\label{E:changeVariables1}
\left.\ba{c}
u = e^{-\eta_t}(1+x^n) \tilde u,\q
\pa_x u = e^{-\eta_t} \big[(1+x^n) \pa_x \tilde u + n  x^{n-1}  \tilde u\big],\\
 \pa^2_{xx} u = e^{-\eta_t}\big[(1+x^n) \pa^2_{xx} \tilde u + 2n x^{n-1} \pa_x \tilde u + n(n-1) x^{n-2} \tilde u\big],\\
  \pa_\o u = e^{-\eta_t}(1+x^n) \pa_\o \tilde u,\q \pa^\o_t u = e^{-\eta_t} (1+x^n) [\pa^\o_t \tilde u - \l \tilde u].
\ea\right.
\end{align}
Define $\tilde{f}:\dbR_T^4\to\dbR$ and $\tilde{g}:\dbR_T^3\to\dbR$ by
\begin{equation}
\label{E:tildefg}
\begin{split}
&\tilde{f}(t,x,y,z,\gamma):=\lambda_t\,y+{e^{\eta_t}\over 1+x^n} \,
f\Big(t,x,{(1+x^n)y\over e^{\eta_t}} ,\\ 
&\qquad {(1+x^n) z + n  x^{n-1}  y\over e^{\eta_t}},
 {(1+x^n) \g + 2n x^{n-1} z + n(n-1) x^{n-2} y\over e^{\eta_t}}\Big),\\
&\tilde{g}(t,x,y,z):= {e^{\eta_t}\over 1+x^n}\, g\Big(t,x,{(1+x^n)y\over e^{\eta_t}}, {(1+x^n) z + n  x^{n-1}  y\over e^{\eta_t}}\Big).
\end{split}
\end{equation}
Clearly, $\tilde f$ and $\tilde g$ inherit the 
 regularity of $f$ and $g$. Whenever they are smooth,
\begin{equation}
\label{tildefgderivative}
\begin{split}
&\pa_y \tilde f = \l + \pa_y f + \pa_z f {nx^{n-1}\over 1+x^n} + \pa_\g f {n(n-1)x^{n-2}\over 1+x^n},\\
&\pa_z \tilde f =\pa_z f + \pa_\g f {2 nx^{n-1}\over 1+x^n};\q  \pa_\g \tilde f = \pa_\g f;\q \pa^2_{\g\g} \tilde f = e^{-\eta_t} (1+x^n) \pa^2_{\g\g} f.
\end{split}
\end{equation}
Then it is straightforward to verify that $\tilde f$ and $\tilde g$ inherit most  desired properties of $f$ and $g$ that we  utilize later.

\begin{lem}
\label{lem-poly}
(i) If $g$ is of the form of \eqref{gsemilinear} or of \eqref{glinear}, 
then so is $\tilde g$; and if $f$ is of the form of  \eqref{flinear}, then so is $\tilde f$.

(ii) If $f$ is convex in $\g$, then so is $\tilde f$. 

(iii) If $f$ is uniformly parabolic, then so is $\tilde f$.

(iv) If $f$ is uniformly Lipschitz continuous in $y$, $z$, $\g$, then so is $\tilde f$.

(v) If $\|f(\cd, y,z,\g)\|_{C^\a(\dbR_T)} \le C[1+|y|+|z|+|\g|]$, then so is $\tilde f$.
\end{lem}
 In particular, if $f$ and $g$ satisfy 
 Assumptions~\ref{assum-g} and \ref{assum-f}, then so do $\tilde f$ and $\tilde g$. 
However, we remark that $\tilde g$ does not inherit the same form when $g$ is in the form of \eqref{glineary}.
Now consider the RPDE for $\tilde u$:
\begin{equation}\label{E:tilde2RPDE}
\begin{split}
\tilde{u}(t,x)&=u_0(x)+\int_0^t
\tilde{f}(s,x,\tilde{u}(s,x),\partial_x\tilde{u}(s,x),\partial^2_{xx}\tilde{u}(s,x))\,ds\\
&\qquad+
\int_0^t \tilde{g}(s,x,\tilde{u}(s,x),\partial_x\tilde{u}(s,x))\,d\o_s,\quad
(t,x)\in \dbR_T.
\end{split}
\end{equation}

\begin{prop}\label{P:ChangeVariables}
Let  Assumptions~\ref{assum-g} and \ref{assum-f} be in force, $\lambda\in C([0,T])$, $n\ge 2$ even, and $u\in C(\dbR_T)$.
Then $u$ is a viscosity subsolution  (resp.~classical subsolution) of RPDE~\eqref{RPDE} if and only if
$\tilde{u}$ is a viscosity subsolution (resp.~classical subsolution) of RPDE~\eqref{E:tilde2RPDE}.
\end{prop}

\proof The equivalence of the classical solution properties is straightforward. Regarding
the viscosity solution properties, we  prove the if part; the only if part can be proved similarly.

Assume that $\tilde u$ is a viscosity subsolution of RPDE \eqref{E:tilde2RPDE}. For any $(t_0, x_0)\in (0, T]\times \dbR$ and $\f \in \underline \cA_g u(t_0, x_0)$, put $\tilde \f(t,x) := {e^{\eta_t} \over 1+x^n} \f(t,x)$. It is straightforward to check that $\tilde \f \in \underline \cA_{\tilde g} \tilde u(t_0, x_0)$. Then, by the viscosity subsolution property of $\tilde u$ at $(t_0, x_0)$, 
\beaa
0 \ge \pa^\o_t\tilde \f - \tilde f(t_0, x_0,  \tilde \f, \pa_x \tilde\f, \pa^2_{xx} \tilde \f)
= {e^{\eta_{t_0}}\over 1+x_0^n} \Big[\pa^\o_t \f - f(t_0, x_0, \f, \pa_x \f, \pa^2_{xx}\f)\Big]. 
\eeaa
This implies that $u$ is a viscosity subsolution of RPDE \eqref{RPDE}.
\qed

\begin{rem}
\label{rem-poly}
{\rm Let $(f, g)$ satisfy  Assumptions~\ref{assum-g} and \ref{assum-f} and let $u$ be a viscosity semi-solution of RPDE~\eqref{RPDE}.

(i)  If $u$ has polynomial growth, by choosing $n$  large enough, we have
\bea
\label{vanish}
\lim_{|x|\to \infty} \sup_{0\le t\le T}|\tilde u(t,x)| = 0.
\eea

(ii)  If $f$ is uniformly Lipschitz continuous in $y$,  by choosing $\lambda$ 
sufficiently large (resp.~small), we have
\bea
\label{fmonotone}
\mbox{$\tilde{f}$ is strictly increasing (resp.~decreasing) in $y$}.
\eea
In particular, $\tilde f$ will be proper in the sense of Crandall, Ishii,  \& Lions \cite{CIL}.
}
\end{rem}

\subsection{Stability}

The following technical lemma will be crucial for the stability result. 
Given $(t_0,x_0)\in [0,T)\times\dbR$ and $\e>0$ small,  put
\begin{align*}
&\qquad\qquad\qquad\qquad   D^+_{\e}(t_0, x_0) := [t_0, t_0+\e^3) \times O_{\e}(x_0),\\
 &\pa D^+_{\e }(t_0, x_0) := \big\{(t,x): t\in [t_0, t_0+\e^3], |x-x_0| = \e~\mbox{or}~ t= t_0+\e^3, |x|\le \e\big\}.
\end{align*}

\begin{lem}\label{L:Stability}
Let  Assumption \ref{assum-g} hold, let $(t_0,x_0)\in [0,T)\times\dbR$, 
and let $\delta_0\in (0,T-t_0]$. Assume that $\f\in C^k_{\a,\b}(D^+_{ \d_0^{1/3}}(t_0,x_0))$ for some large $k$ and $\partial_\omega\varphi=g(\cdot,\varphi,\partial_x\varphi)$ in $D^+_{ \d_0^{1/ 3}}(t_0,x_0)$. Define  
\bea
\label{gphi}
g^\varphi(t,x,y,z):=[g(\cdot,\varphi+y,\partial_x\varphi+z)- g(\cdot,\varphi,\partial_x\varphi)](t,x). 
\eea
Then there exists an $\eps_0\in (0,1]$ such that, for every  $\eps\in(0,\eps_0]$, there exists  a function $\psi^\eps\in C^{4}_{\a,\b} (D^+_{\e}(t_0,x_0))$  that satisfies the following properties:
\bea
\label{psieest1}
&\dis \partial_\omega\psi^\eps =g^\varphi(\cdot,\psi^\eps,\partial_x\psi^\eps),\q \pa^\o_t \psi^\e = \e;&\\
\label{psieest2}
&\dis \lvert\psi^\eps\rvert+ \lvert\partial_x\psi^\eps\rvert+\lvert\partial^2_{xx}\psi^\eps\rvert \le C\eps^2~ \mbox{in}~ D^+_{\e}(t_0,x_0);&\\
\label{psieest3}
&\dis\psi^\eps(t_0,x_0)<0<\inf_{ (t,x)\in\partial D^+_{\e} (t_0,x_0)} \psi^\eps(t,x).&
\eea
\end{lem}

\proof Without loss of generality, we let $(t_0,x_0)=(0,0)$. Since our results are local, without loss of generality, we can  assume  that $\f\in C^k_{\a,\b}(\dbR_T)$.
Let $\iota\in C^\infty(\dbR)$ be  such that $\iota(x) = x^4$ for $|x|\le 1$ and $\iota(x)=0$ for $|x|\ge 2$. For any $\e>0$ small, consider the 
RPDE
\bea
\label{psie}
 \psi^\e(t,x) = \iota(x) - \e^5 + \e t + \int_0^t g^\f(s,x, \psi^\e, \pa_x \psi^\e) d\o_s.
 \eea
 By Theorem \ref{thm-1st}, there exists  a $\d_1\le \d_0$ such that $\psi^\e \in   C^{4}_{\a,\b}([0,\d_1]\times\dbR)$ for all $\e \le 1$ and 
 \bea
 \label{psienorm}
 \sup_{\e\le 1} \|\psi^\e\|_4 <\infty.
 \eea
The equalities in \eqref{psieest1} are obvious.  Now, we  verify that $\psi^\e$ satisfies \eqref{psieest2} and \eqref{psieest3}.
 Recall the fourth-order Taylor expansion in Lemma \ref{lem-Taylor}:
  \bea
  \label{psieTaylor}
  \psi^\e(t,x) = \sum_{ \ell+|\nu|\le 4} {1\over \ell!}  (\pa^\ell_x\cD_\nu\psi^\eps)(0,0) x^\ell\cI^\nu_{0,t} + R_4(t,x).
  \eea
  We claim that
  \begin{equation}
  \label{papsiClaim}
  \begin{split}
&  \pa^4_x\psi^\e(0,0) = 24,\q  \pa^{\ell_1}_x \pa^{\ell_2}_\o    \pa_\o \psi^\e(0,0)  = O(1)\q\mbox{for}~ \ell_1 + \ell_2 = 3,\\
 &  \pa^\o_t \psi^\e(0,0) = \e,\q  \pa^\o_t  \pa_\o \psi^\e(0,0) =  O(\e),\q    \pa^\o_t  \pa^2_\o \psi^\e(0,0)  = O(\e),\\
 &  \pa^\ell_x\cD_\nu\psi^\eps(0,0)  = O(\e^5) ~ \mbox{for all other terms such that}~ \ell + |\nu|\le 4.
   \end{split}
   \end{equation}
   Note further that, by \eqref{psienorm}, for every $(t,x) \in \overline D^+_\e(0,0)$,
   \begin{align}
   \label{R4}
  | R_4(t,x) | \le C[t^\a + |x|]^{4+\b} \le C[ \e^{3\a} + \e]^{4+\b} \le C\e^{4+\b}.
  \end{align}
  Then, for $(t,x) \in \overline D^+_\e(0,0)$, plugging \eqref{papsiClaim} and \eqref{R4} into \eqref{psieTaylor}, we obtain
  \beaa
 | \psi^\e(t,x)| &\le&  C \Big[\sum_{\ell_1+\ell_2=4}  |x|^{\ell_1}  t^{\a\ell_2}+ \e [t + t^{1+\a} + t^{1+2\a}] + \e^5+\e^{4+\b}\Big]\\
 &\le&   C\Big[ \sum_{\ell_1+\ell_2=4}  \e^{\ell_1+ 3\a \ell_2}  + \e^4\Big] \le C\e^4.
  \eeaa
  Similarly, applying the third-order Taylor expansion of $\pa_x \psi^\e$ and the second-order Taylor expansion of $\pa^2_{xx} \psi^\e$, we obtain 
  \beaa
   |\pa_x \psi^\e(t,x)| \le C\e^3,\q  |\pa^2_{xx} \psi^\e(t,x)|\le C\e^2,\q (t,x) \in \overline D^+_\e(0,0).
   \eeaa
   Thus we have proved \eqref{psieest2}.
  
  To prove \eqref{psieest3}, let $(t,x) \in \pa D^+_{\e}(0,0)$.  By 
  \eqref{psieTaylor}, \eqref{papsiClaim}, and \eqref{R4},
\begin{align*}
\psi^\e(t,x) \ge \e t + {24\over 4!} x^4 - C\big[\sum_{\ell_1 + \ell_2 = 3}  |x|^{\ell_1} t^{\a (\ell_2+1)} +\e [ t^{1+\a} + t^{1+2\a}]+\e^5 + \e^{4+\b}\big].
\end{align*}
 Note that in both, the case   $t\in [0, \e^3]$ and $|x| = \e$ and the case $t= \e^3$ and $|x|\le \e$, we have 
$  \e t + {24\over 4!} x^4 \ge \e^4.
 $
  Hence 
 \beaa
\psi^\e(t,x)  &\ge& \e^4 - C\sum_{\ell_1 + \ell_2 = 3}  \e^{\ell_1 + 3\a (\ell_2+1)} - C\e^{1+3(1+\a)} - C\e^5 - C \e^{4+\b}\\
&\ge& \e^4 - C\big[\e^{3 + 3\a} + \e^{4+\b}] \ge {1\over 2} \e^4,
\eeaa
 when $\e$ is small, thanks to the assumption that $3\a>1$.  Moreover, it is clear that $\psi^\e(0,0) = -\e^5 <0$. Thus we have proved  \eqref{psieest3}.

 It remains to prove \eqref{papsiClaim}. First, by \eqref{psie} it is clear that 
 \begin{align}
 \label{papsi0}
  \psi^\e(0,x) = \iota(x) -\e^5,\,  \pa^\o_t \psi^\e(t,x) = \e,\, \pa_\o \psi^\e(t,x) = g^\f(t,x, \psi^\e, \pa_x \psi^\e).
 \end{align}
  By the first two equalities of \eqref{papsi0}, one can easily see that
   \begin{equation}
 \label{papsi1}
 \begin{split}
& \psi^\e(0,0) = \e^5 <0, \, \pa^4_x\psi^\e(0,0) = 24,\, \pa^\ell_x \psi^\e(0,0)=0,~\mbox{for}~  1\le \ell\le 3,\\
 &  \pa^\o_t \psi^\e(0,0) = \e, \,  \pa^\ell_x \cD_\nu \pa^\o_t \psi^\e(0,0)=0~\mbox{for}~ 1\le \ell + |\nu|\le 2,\\
 &  \pa_\o \psi^\e(0,0) = g\big(0,0, \f(0,0) - \e^5, \pa_x\f(0,0)\big)\\
   &\qquad\qquad\quad - g\big(0,0, \f(0,0), \pa_x\f(0,0)\big)  = O(\e^5).
 \end{split}
 \end{equation}
 All  terms above satisfy \eqref{papsiClaim}.
 The remaining derivatives
  involved 
 in \eqref{psieTaylor} take the form
    \begin{equation}
 \label{papsi2}
 \pa^\ell_x \cD_\nu \pa_\o \psi^\e(0,0) = \pa^\ell_x \cD_\nu\big[g^\f(\cd, \psi^\e, \pa_x \psi^\e)\big]\big|_{(t,x) = (0,0)},\,1\le \ell + |\nu|\le 3.
 \end{equation}
 Note that
 \begin{align}
 \label{papsi3}
 [\pa^\ell_x \cD_\nu g^\f](\cd, \psi^\e, \pa_x \psi^\e)\big]\big|_{(t,x) = (0,0)} = [\pa^\ell_x \cD_\nu g^\f](0, 0, -\e^5, 0) = O(\e^5).
 \end{align}
 Then \eqref{papsi2} becomes
  \begin{align}
 \label{papsi4}
\pa^\ell_x \cD_\nu \pa_\o \psi^\e(0,0) = \pa^\ell_x \cD_\nu\big[g^\f(0, 0, \psi^\e(0,x), \pa_x \psi^\e(0,x))\big]\big|_{x = 0} + O(\e^5) \\
   =  \pa^\ell_x \cD_\nu\big[g\big(0, 0, \f(0,0)+ \psi^\e(0,x),  \pa_x\f(0,0)+ \pa_x \psi^\e(0,x)\big)\big]\big|_{x = 0}+ O(\e^5).\nonumber
 \end{align}
 Thus the derivatives are combinations of terms involving the derivatives of $g$ with respect to $(y,z)$, which are all bounded by our assumption, and the derivatives of $(\psi^\e, \pa_x \psi^\e)$. By a tedious but quite straightforward computation of the derivatives, we obtain from \eqref{papsi1} and \eqref{papsi4} 
 with the abbreviation
 $\eta(0):= \eta(0,0, \f(0,0), \pa_x \f(0,0))$
  for any function $\eta$,
 \begin{align*}
&  \pa^\o_t  \pa_\o \psi^\e(0,0) =  \pa_y g(0) \e + O(\e^5),\q    \pa^\o_t  \pa^2_\o \psi^\e(0,0)  = |\pa_y g(0)|^2 \e + O(\e^5),\\
&\qquad \pa^{\ell_1}_x \pa^{\ell_2}_\o    \pa_\o \psi^\e(0,0)  = 24 [\pa_z g(0)]^{\ell_2 + 1}  + O(\e^5)\q\mbox{for}~ \ell_1 + \ell_2 = 3,
 \end{align*}
 and all other terms  either contain $\pa^\ell_x \psi^\e(0,0)=0$ for some $1\le \ell\le 3$, 
 or $\pa^\ell_x \cD_\nu \pa^\o_t \psi^\e(0,0)=0$ for some $1\le \ell + |\nu|\le 2$, or 
 \begin{align*}
 [\pa^\ell_x \cD_\nu g^\f](\cd, \psi^\e, \pa_x \psi^\e)\big]\big|_{(t,x) = (0,0)} = O(\e^5)
 \end{align*}
  for some $\ell + |\nu|\le 3$. We thus prove \eqref{papsiClaim} for all the cases, and hence complete the proof of the lemma.
 \qed
 
\begin{thm}[Stability]\label{T:Stability}
Let Assumption \ref{assum-g} hold and $(f_n)_{n\ge 1}$ be a sequence of functions satisfying Assumption \ref{assum-f}.  For each $n\ge 1$, let $u_n$ be a viscosity subsolution of RPDE \eqref{RPDE} with generator $(f_n, g)$.   Assume further that, for some functions $f$ and $u$, 
\begin{align}\label{E:Stability}
\lim_{n\to\infty} [f_n-f](t,x,y,z,\gamma)=0\quad\text{and}\quad
\lim_{n\to\infty} [u_n-u](t,x)=0
\end{align}
locally uniformly in $(t,x,y,z,\gamma)\in \dbR_T^4$.
Then $u$ is a viscosity subsol.~of \eqref{RPDE}.
\end{thm}

\proof By the locally uniform convergence, 
$f$ and $u$ are continuous.
Let $(t_0,x_0)\in (0,T]\times\dbR$ and $\varphi\in\underline{\mathcal{A}}_g u(t_0,x_0)$. We 
 apply Lemma \ref{L:Stability} at 
$(t_0,x_0)$, but in the left neighborhood
\begin{align}
\label{De-}
 D^-_{\e}(t_0, x_0) := (t_0-\e^3, t_0] \times O_{\e}(x_0),\hskip3cm\\
\pa D^-_{\e }(t_0, x_0) := \big\{(t,x):  (t_0-\e^3, t_0], |x-x_0| = \e~\mbox{or}~ t= t_0-\e^3, |x|\le \e\big\}.\nonumber
\end{align}
We emphasize that, while for notational simplicity we established Lemma \ref{L:Stability} in the right neighborhood $D^+_{\e}(t_0, x_0)$, we may easily reformulate it to the left neighborhood  by using tbackward rough paths  introduced in \eqref{back}.  
By Remark \ref{R:jet}, we may assume without loss of generality that $\f\in C^k_{\a,\b}(\dbR_T)$ for some large $k$. Then, for any $\e>0$ small, by Lemma \ref{L:Stability}, there exists  $\psi^\eps\in C^{4}_{\a,\b}
(D^-_{\e}(t_0,x_0))$
such that the following holds:
\bea
\label{stability-psie}
\left.\ba{c}
\dis  \partial_\omega\psi^\eps =g^\varphi(\cdot,\psi^\eps,\partial_x\psi^\eps),\q \pa_t^\o\psi^\eps = \e;\\
\dis \lvert\psi^\eps\rvert+ \lvert\partial_x\psi^\eps\rvert+ \lvert\partial^2_{xx}\psi^\eps\rvert \le C\eps^2~ \mbox{in}~ D^-_{\e}(t_0,x_0);\\
\dis\psi^\eps(t_0,x_0)<0<\inf_{ (t,x)\in\partial D^-_{\e} (t_0,x_0)} \psi^\eps(t,x).
\ea\right.
\eea
This together with setting $\f^\e := \f + \psi^\e$ yields
\begin{align}
\label{phie}
\sup_{(t,x)\in\partial D^-_{\e} (t_0,x_0)} \big[[u-\varphi^\e](t,x)\big]<0<[u-\varphi^\e](t_0,x_0). 
\end{align}
Since $u^n$ converges to $u$ locally uniformly, we have, for $n=n(\eps)$ large enough,
\begin{align*}
\sup_{(t,x)\in\partial D^-_{\e} (t_0,x_0)} \big[[u_n-\varphi^\e](t,x)\big]<0<[u_n-\varphi^\e](t_0,x_0).
\end{align*}
Then there exists  $(t_\e, x_\e) = (t^n_\eps, x^n_\eps)\in D^-_{\e} (t_0,x_0)$ such that
\begin{align*}
[u_n-\varphi^\e](t_\eps,x_\eps) =0= \max_{[t_0-\e^3,t_\eps]\times\bar{O}_\eps(x_0)} \big[[u_n-\varphi^\e](t,x)\big].
\end{align*}
Note that 
\begin{align*}
\pa_\o \f^\e &= \pa_\o \f + \pa_\o \psi^\e = g(\cd, \f, \pa_x\f) \\
&+ \big[g(\cd, \f + \psi^\e, \pa_x \f + \pa_x \psi^\e) - g(\cd, \f, \pa_x\f)\big] = g(\cd, \f^\e, \pa_x\f^\e).
\end{align*}
Then $\f^\e\in \underline{\mathcal{A}}_g u_n(t_\eps,x_\eps)$.
 By the viscosity subsolution property of $u_n$, 
\begin{align*}
\partial_t^\o\varphi^\eps(t_\eps,x_\eps)-
f_n(\cdot,\varphi^\eps,\partial_x{\varphi^\eps},
\partial^2_{xx}\varphi^\eps)
(t_\eps,x_\eps)\le 0.
\end{align*}
Fix $n$ and send  $\eps\to 0$. Then,  by the convergence of $\psi^\e$ and its derivatives, 
\beaa
\partial_t^\o\varphi(t_0,x_0)-f_n(\cdot,\varphi,\partial_x\varphi,\partial^2_{xx}\varphi)
(t_0,x_0)\le 0,
\eeaa
Now  by sending $n\to \infty$ we get
$
\partial_t^\o\varphi(t_0,x_0)-
f(\cdot,\varphi,\partial_x\varphi,\partial^2_{xx}\varphi)
(t_0,x_0)\le 0,
$
i.e.,  $u$ is a viscosity subsolution of \eqref{RPDE}. \qed

\section{Viscosity solutions of rough PDEs: comparison principle}
\label{sect-comparison}
\setcounter{equation}{0}
Let
\begin{equation}
\label{boundary}
\begin{split}
&\text{$u_1$ (resp. $u_2$)  be  a viscosity subsol. (resp. supersol.) of RPDE \eqref{RPDE},}\\
&u_1(0,\cd) \le u_2(0,\cd),\q\mbox{and $u_1, u_2$ have polynomial growth in $x$}.
\end{split}
\end{equation}
Our goal is to show that $u_1 \le u_2$ on $\dbR_T$.  

When both $u_1$ and $u_2$ are smooth, 
$u:= u_1-u_2$ solve a linear RPDE. Then the function $F$ corresponding to this linear RPDE becomes linear, see \eqref{linearF} below.  Thus, by using the representation formula \eqref{linearABCF0}-\eqref{linear-rep2}  below, one can easily show that $u_1 \le u_2$.

\subsection{Partial comparison principle}
Here, we assume that one of $u_1$ and $u_2$ is smooth. We need the following result
(cf.~Lemma \ref{L:Stability}).

\begin{lem}
\label{lem-PartCompGeneral}
Let Assumption \ref{assum-g} be in force. Let $\f\in C^k_{\a,\b}(\dbR_T)$ for some large $k$. Let 
$\partial_\omega\varphi=g(\cdot,\varphi,\partial_x\varphi)$ on $\dbR_T$. For any $0\le t_0 < T$, $0<\d\le T-t_0$,  and $\e>0$, recall \eqref{gphi}, and consider the  RPDE 
\begin{align}
\label{psie2}
\psi^\e(t,x) = \e + t - t_0+ \int_{t_0}^t g^\f(s, x, \psi^\e, \pa_x \psi^\e) d\o_s,\qquad [t_0, t_0+\d] \times \dbR.
\end{align}
Then $\psi^\e\in C^{2}_{\a,\b}([t_0, t_0+\d] \times \dbR)$ with $\|\psi^\e\|_{C^{2}_{\a,\b}([t_0, t_0+\e] \times \dbR)} \le C$, where $C$ depends only on $g$ and $\f$, but not on $t_0$, $\e$, and $\d$. 
Moreover, $\psi^\e$ satisfies
\bea
\label{psieest}
\left.\ba{c}
\dis \partial_\omega\psi^\eps =g^\varphi(\cdot,\psi^\eps,\partial_x\psi^\eps),\q \pa^\o_t \psi^\e = 1,\\
\dis \lvert\psi^\eps\rvert+ \lvert\partial_x\psi^\eps\rvert+\lvert\partial^2_{xx}\psi^\eps\rvert \le C[\eps+\d^{\a\b}]~ \mbox{in}~ [t_0, t_0+\d]\times \dbR,\\
\dis \inf_{x\in \dbR} \psi^\eps(t,x) >0 \q\mbox{for all}~ t\in (t_0, t_0+\d].
\ea\right.
\eea
\end{lem}
\proof The uniform regularity of $\psi^\e$ and the first line of \eqref{psieest} are clear.  Note that $\psi^\e(t_0,x) = \e$, $\pa_x\psi^\e(t_0,x) = 0$, $\pa^2_{xx} \psi^\e(t_0,x) = 0$. The second line of \eqref{psieest} follows from the H\"{o}lder continuity of the functions in terms of $t$. Moreover, since $g^\f(t,x, 0, 0)=0$, we may write it as $g^\varphi(t,x,\psi^\eps,\partial_x\psi^\eps) = \si(t, x) \psi^\e + b(t,x) \pa_x \psi^\e$, where $\si$ and $b$ depend on $\psi^\e$. Then we may view \eqref{psie2} as a linear RPDE with coefficients $\si$ and $b$.
Thus, by \eqref{linearABCF0}-\eqref{linear-rep2}, we have a representation formula for $\psi^\e$.    The uniform regularity of $\psi^\e$ implies the uniform regularity of $\si$ and $b$, which 
leads to the third line of \eqref{psieest}. 
\qed

\begin{thm}\label{thm-partial}
Let Assumptions~\ref{assum-g} and \ref{assum-f} and \eqref{boundary} be in force. If one of $u_1$ and $u_2$ is in $C^k_{\a,\b}(\dbR_T)$ for some large $k$, then $u_1 \le u_2$.
\end{thm}

\proof For the sake of a contradiction, assume that 
$[u_1-u_2](t_0, x_0) >0$ for some $(t_0, x_0) \in (0, T]\times \dbR$.  By Remark \ref{rem-poly} (i), without loss of generality we may assume both $u_1$ and $u_2$ satisfy \eqref{vanish}. Put
\begin{align*}
c_\d := \sup_{(t,x) \in [0, \d] \times \dbR} [u_1-u_2](t,x),\q \d_0 := \inf\{\d \ge 0: c_\d >0\}.
\end{align*}
Then $c_\d$ is nondecreasing in $\d$,  $c_0  \le 0 < c_{t_0}$, and thus $\d_0 < t_0$. For any $0<\d \le t_0-\d_0$, $c_{\d_0 + \d} >0$. By \eqref{vanish} and since $u(0,\cd) \le 0$, there exists $(t_\d, x_\d) \in (\d_0, \d_0 + \d]\times \dbR$ such that $(u_1-u_2)(t_\d, x_\d) = c_{\d_0+\d}$. Set $\e := c_{\d_0+\d}\wedge \d^{\a\b}$. Applying Lemma \ref{lem-PartCompGeneral} with $\f := u_2$ on $[\d_0,  t_\d]$,  but again backwardly in time, we have $\psi^\e$ satisfying
\beaa
\left.\ba{c}
\dis\psi^\e(t_\d,x) = \e,\q   \pa^\o_t \psi^\e = 1,\q \inf_{x\in \dbR} \psi^\eps(\d_0,x) >0,\\
  \lvert\psi^\eps\rvert+ \lvert\partial_x\psi^\eps\rvert+\lvert\partial^2_{xx}\psi^\eps\rvert \le C\d^{\a\b}~ \mbox{in}~ [t_0, t_\d]\times \dbR.
\ea\right.
\eeaa
Define $\f^\e := u_2 + \psi^\e$. Note that
$
[u_1- \f^\e](t_\d, x_\d) \ge 0 > \sup_{x\in \dbR} [u_1 - \f^\e](\d_0, x).
$
Then there exists $(t_\d^*, x_\d^*) \in (\d_0, t_\d]\times \dbR$ such that 
\beaa
[u_1- \f^\e](t^*_\d, x^*_\d) = 0 = \sup_{(t,x) \in [\d_0, t^*_\d]\times \dbR} [u_1- \f^\e](t^*_\d, x^*_\d).  
\eeaa
By the definition of $g^\f$, it is clear that $\pa_\o \f^\e = g(\cd, \f^\e, \pa_x\f^\e)$. Then $\f^\e \in \underline{\mathcal{A}}_g u_1(t^\ast_\d, x^\ast_\d)$. Thus, by using the classical supersolution property of $u_2$ and the viscosity subsolution property of $u_1$ we have
\begin{align*}
\Bigl[
&\partial_t^\o u_2-f(\cdot,u_2,\partial_x u_2, \partial^2_{xx} u_2) \Bigr](t^\ast_\d,x^\ast_\d)\\&\qquad\ge 0
\ge  \Bigl[ \partial_t^\o \f^\e - f(\cdot, \f^\e, \partial_x \f^\e, \partial^2_{xx} \f^\e)\Bigr](t^\ast_\d, x^\ast_\d).
\end{align*}
Now, at $(t^*_\d, x^*_\d)$, we have
\beaa
1 &=& \partial_t^\o \f^\e  - \partial_t^\o u_2 \le f(\cdot, \f^\e, \partial_x \f^\e, \partial^2_{xx} \f^\e) - f(\cdot,u_2,\partial_x u_2, \partial^2_{xx} u_2)\\
&\le& C\Big[ \lvert\psi^\eps\rvert+ \lvert\partial_x\psi^\eps\rvert+\lvert\partial^2_{xx}\psi^\eps\rvert\Big] \le C\d^{\a\b},
\eeaa
which is an obvious contradiction when $\d$ is small.
 \qed

\begin{rem}
\label{rem-PartComp}
{\rm When $g$ is independent of $y$, we can prove Proposition \ref{thm-partial}  much easier without invoking Lemma \ref{lem-PartCompGeneral}.  In fact, in this case, assuming  to the contrary that  $(u_1-u_2)(t_0,x_0)>0$ for some $(t_0,x_0)\in \dbR_T$. Then
\begin{align*}
c:=\sup_{(t,x)\in [0,t_0]\times\dbR} [u_1-u_2](t,x)\ge [u_1-u_2](t_0,x_0)>0.
\end{align*}
By \eqref{vanish} and  $[u_1-u_2](0,\cd) \le 0$, there exists $(t^*, x^*) \in (0, t_0] \times \dbR$ such that 
\beaa
[u_1-u_2](t^*, x^*) = c := \sup_{(t,x)\in [0,t_0]\times\dbR} [u_1-u_2](t,x)\ge [u_1-u_2)]t_0,x_0)>0.
\eeaa
Define $\f = u_2 + c$. Since $g$ is independent of $y$, we have 
\beaa
\pa_\o \f = \pa_\o u_2  = g(t,x, \pa_x u_2) = g(t,x,\pa_x \f).
\eeaa
Then one can easily verify that $\varphi\in\underline{\mathcal{A}}_g u_1(t^\ast,x^\ast)$. Moreover,  by Remark \ref{rem-poly} (ii),  we can assume without loss of generality that $f$  is strictly decreasing in $y$.
Now it follows from the classical supersolution property of $u_2$ and the viscosity subsolution property of $u_1$ that, taking values at $(t^\ast, x^\ast)$,
\begin{align*}
&\partial_t^\o u_2-f(\cdot,u_2,\partial_x u_2, \partial^2_{xx} u_2) ~\ge~ 0~\ge~ \partial_t^\o \f - f(\cdot, \f, \partial_x \f, \partial^2_{xx} \f)\\
&\qquad\qquad\qquad
= \partial_t^\o u_2-f(\cdot,u_2+c,\partial_x u_2, \partial^2_{xx} u_2),
\end{align*}
which is the desired contradiction since   $f$ is strictly decreasing  in $y$. 
}
\end{rem}

The following comparison result follows immediately from Theorem~\ref{thm-partial}.

\begin{cor}\label{cor-partial}
Let Assumptions~\ref{assum-g} and \ref{assum-f} and \eqref{boundary} be in force. If  RPDE \eqref{RPDE}  has a classical solution $u\in C^k_{\a,\b}(\dbR_T)$ for some large $k$ and $u_1(0,\cd) \le u(0,\cd) \le u_2(0,\cd)$, then $u_1 \le u\le  u_2$. In particular, $u$ is the unique solution in the viscosity sense.
\end{cor}

\subsection{Full comparison}
We shall follow the approach in Ekren, Keller, Touzi, \& Zhang \cite{EKTZ}.  For this purpose, we strengthen Assumption \ref{assum-g} slightly by imposing some uniform property of $g$ in terms of $y$.

\begin{assum}
\label{assum-g1}
The diffusion coefficient $g$ belongs $C^{k_0, loc}_{\a,\b}(\dbR_T^3)$ for some $k_0$ large enough, and

(i) $\pa_z g \in C^{k_0-1}_{\a,\b}(\dbR_T^3)$,

(ii) for $i=0,\cds, k_0$ and $z \in \dbR$, $\pa^i_x g(\cd,z), \pa^i_y g(\cd,z) \in C^{k_0-i}_{\a,\b}(\dbR_T^2)$ with 
$\|\pa^i_x g(\cd,z) \|_{k_0-i} + \|\pa^i_y g(\cd,z) \|_{k_0-i} \le C [1+|z|]$.
\end{assum}
We remark that, under Assumption \ref{assum-g}, all the results in this subsection hold true if we assume instead that $T$ is small enough.

Given an initial condition $u_0$, motived by the partial comparison, we fix a large $k$ and define
\begin{equation}
\label{ubar}
\begin{split}
\overline u(t,x) &:= \inf\big\{ \f(t,x):  \f  \in \overline \cU\big\},\q \underline u(t,x) := \sup\big\{ \f(t,x):  \f  \in \ul \cU\big\},
\end{split}
\end{equation} 
where
\begin{equation}
\label{Ubar}
\begin{split}
\cU &:= \big\{\f \in \dbL^0(\dbR_T):   \text{\cag in $t$, with polynomial growth in $x$,}\\
&\qq\text{ $\f(0,\cd) = u_0$,}
\text{ and $\exists~ 0=t_0<\cds<t_n=T$ such that}\\
&\qq\qq\text{$\f \in C^k_{\a,\b}((t_{i-1}, t_{i}]\times \dbR)$ for $i=1$, $\ldots$, $n$} \big\},\\
\overline\cU \!&:=\! \big\{\f \in\cU:  \D \f_{t_i} \ge 0, ~\text{$\f$ is a classical supersolution of }\\
&\qq\qq\text{RPDE \eqref{RPDE} on each $(t_{i-1}, t_i]$}   \big\},\\
\underline\cU &:=\big\{\f \in\cU:  \D \f_{t_i} \le 0, ~\text{$\f$ is a classical subsolution of }\\
&\qq\qq\text{RPDE \eqref{RPDE} on each $(t_{i-1}, t_i]$}   \big\}.
\end{split}
\end{equation}

\begin{lem}
\label{lem-cU}
Let Assumptions~\ref{assum-g1}, \ref{assum-f}, and \ref{assum-u0}  hold. Then  $\overline\cU$, $\underline \cU\neq\emptyset$.
\end{lem}
\proof  We  prove $\ol\cU \neq \emptyset$  in several steps. The proof for $\ul\cU$ is similar.

{\it Step 1.} Put $Q_1 := \dbR^2 \times \{z\in \dbR: |z|\le 1\}$. Then $g \in C^{k_0}_{\a,\b}([0, T]\times Q_1)$ and let $N_0$ denote its $k_0$-norm. Under our strengthened conditions in Assumption~\ref{assum-g1}, it follows from the arguments in Proposition \ref{P:CharExistence} that there exist $\d_0, C_0$, depending only on $N_0$, such that
\begin{align}
\label{Charest1}
 \Th^0 \in C^{k_0}_{\a,\b}([0,\delta_0]\times Q_1; \dbR^3)~  \mbox{with}~ \| \Th^0\|_{k_0}\le C_0,  
\end{align}
where  $\Th^0_t(\th) := \Th_t(\th) - \th$.
Moreover, for a possibly smaller $\d_0>0$, again depending only on $N_0$, we have 
\bea
\label{paXest1}
\pa_x X_t(x,y,0) \ge {1/ 2} \q\mbox{for all}~ (t,x,y) \in [0, \d_0]\times \dbR^2.
\eea

{\it Step 2.} Recall \eqref{F} and put
\begin{align*}
&\ol F(t,y) := \sup_{x\in \dbR} F(t,x,y,0,0)\\ &= \sup_{x\in \dbR} f\big(t, \Th_t(x,y,0), {\pa_x Z_t(x,y,0) \over \pa_x X_t(x,y,0)}\big)\,
 \exp\big(-\int_0^t \pa_y g(s, \Th_s(x,y,0))d\o_s\big).
\end{align*}
By  Assumption \ref{assum-f},  \eqref{Charest1} and \eqref{paXest1} yield, for all $(t,x,y)\in [0, \d_0] \times \dbR^2$,
\begin{align*}
|F(t,x,y,0, 0)| &\le C_1\Big[1+ |Y_t(x,y,0)|+ |Z_t(x,y,0)| + \big|{\pa_x Z_t(x,y,0) \over \pa_x X_t(x,y,0)}\big|\Big]\\
 &\le C_1[1+|y|],
\end{align*}
where $C_1$ depends on $N_0$ and the $K_0$ and $L_0$  in Assumption \ref{assum-f}. Then
$|\ol F(t,y)|\le C_1[1+|y|]$.
Moreover, it is clear that $F(t,x,y,0,0)$ is differentiable in $y$. However, due to the exponential term outside of $f$, in general $\pa_y F(t,x,y,0,0)$ may not be bounded, and we can only claim that $|\pa_y F(t,x,y,0,0)| \le C_1[1+|y|]$.  By the regularity of $f$, it is clear that $\ol F$ is continuous in $t$. Let $\hat F$ be a smooth mollifier of $\ol F$  such that  
\begin{align}
\label{olFest1}
\ol F \le \hat F \le \ol F+1, \q |\hat F(t,y)| \le C_1[1+|y|],\q |\pa_y  \hat F(t,y)| \le C_1[1+|y|].
\end{align}
Set $K_1 := \|u_0\|_\infty e^{C_1 T}$ for the above $C_1$. We see that $\hat F(t,y)$ is uniformly Lipschitz 
in $y$ on $[0, \d_0] \times [-K_1-1, K_1+1]$. Let $\iota$ be a smooth truncation function such that $\iota(x) = x$ for $|x|\le K_1$, $\iota(x) = \mbox{sign}(x) [K_1+1]$ for $|x|\ge K_1+1$, and $|\iota(x)|\le |x|$ for all $x$. Now  consider the 
ODE
\begin{align}
\label{bpsi}
\ol \psi_t =  \|u_0\|_\infty + \int_0^t \hat F(s, \iota(\ol\psi_s))ds,\q  0\le t\le \d_0.
\end{align}
Clearly, \eqref{bpsi} has a solution $\ol \psi \in C^\infty([0, \d_0])$. Since $|\hat F(s, \iota(y))|\le C_1[1+|\iota(y)|] \le C_1[1+|y|]$, 
 $|\ol \psi_t| \le \|u_0\|_\infty e^{C_1 t} \le K_1$ for $t\le \d_0$. Thus $\iota(\ol\psi_t) = \psi_t$ and  
\bea
\label{olpsi}
\ol \psi_t =  \|u_0\|_\infty + \int_0^t \hat F(s, \ol\psi_s)ds,\q  0\le t\le \d_0.
\eea
Abusing the notation by letting $\ol \psi(t,x) := \ol \psi(t)$, we have
\begin{align*}
&\qquad\qquad \pa_\o \ol \psi_t =0,\q \pa_x \ol\psi_t = 0,\q   \pa^2_{xx} \ol\psi_t = 0, \\
& \pa_t \ol\psi (t,x) = \hat F(t, \ol\psi_t)\ge  \ol F(t, \ol\psi_t) \ge F(t,x, \ol \psi, \pa_x \ol \psi, \pa^2_{xx} \ol\psi),
\end{align*}
i.e., $\ol \psi$ is a classical supersolution of PDE \eqref{PDE}-\eqref{F}.  Thus \eqref{hatTh} becomes
\beaa
\hat \Th_t(x) := \Th_t(x, \ol\psi(t,x), \pa_x\ol\psi(t,x)) = \Th_t(x, \ol\psi_t, 0).
\eeaa
in this case. By \eqref{paXest1},  for any $t\in [0, \d_0]$,  $x\mapsto \hat X_t(x)$ is invertible. Put 
\beaa
\ol \f(0,x) := u_0(x),\q \ol \f(t,x):=  \hat Y_t(\hat X_t^{-1}(x)), ~ t\in (0, \d_0].
\eeaa
Then $\D \ol \f_0(x) \ge 0$, $\ol \f \in C^k_{\a,\b}((0, \d_0]\times \dbR)$, and by  Theorem \ref{T:EU2RPDE}, $\ol \f$ is a classical supersolution of PPDE \eqref{PPDE} on $(0, \d_0] \times \dbR$.

{\it Step 3.}  Let $\d_0$ be as in Step 1. We emphasize that $\d_0$ depends only on $N_0$, in particular 
not on $\|u_0\|_\infty$. Let $0=t_0<\cds<t_n=T$ be a partition such that $t_i-t_{i-1} \le \d_0$, $i=1$, $\ldots$, $n$. By Step 2, we have a desired function $\ol\f$ on $[t_0, t_1]\times \dbR$. In particular,
$\ol\f(t_1,\cd)$ is bounded. Now consider RPDE \eqref{RPDE} on $[t_1, t_2]$ with initial condition $\ol\f(t_1, \cd)$. Following the same arguments, we may extend $\ol \f$ to $[t_0,t_2]$ such that $\D \ol\f(t_1,\cd) \ge 0$, $\ol\f \in C^k_{\a,\b}((t_1, t_2]\times \dbR)$, and $\ol \f$ is a classical supersolution of RPDE \eqref{RPDE} on $(t_1, t_2] \times \dbR$. Repeating the arguments yields the desired $\ol\f$ on $\dbR_T$, i.e., $\ol\f\in \ol \cU$.
 \qed

Now, by Theorem \ref{thm-partial} and Proposition \ref{prop-consistency}, it is clear that 
\bea
\label{ulu<olu}
\ul u \le \ol u.
\eea
We next establish the viscosity solution property of $\ol u$ and $\ul u$. We shall follow the arguments in Theorem \ref{T:Stability}, which relies on the crucial Lemma \ref{L:Stability}.
\begin{lem}
\label{lem-oluvis}
Let Assumptions~\ref{assum-g1}, \ref{assum-f}, and \ref{assum-u0}  hold. 

(i) $\ol u$ (resp.  $\ul u$) is  bounded and upper (resp. lower) semi-continuous. 

(ii) Moreover, if $\ol u$ (resp. $\ul u$) is continuous,
 then $\ol u$ (resp. $\ul u$)  is a viscosity supersolution (resp.  viscosity subsolution) of RPDE \eqref{RPDE}.
\end{lem}
We remark that  it is possible to extend our definition of viscosity supersolutions to lower semicontinuous functions. However,  here (i) shows that $\ol u$ is upper semicontinuous. So it seems that the continuity of $\ol u$ in (ii)  is intrinsically required in this approach.

\ms
\proof By the proof of Lemma \ref{lem-cU}, $\ol u$ is bounded from above. Similarly $\ul u$ is bounded from below. Then it follows from \eqref{ulu<olu} that $\ol u, \ul u$ are bounded.

We  establish next the  upper semicontinuity  for $\ol u$. The regularity for $\ul u$ can be proved similarly. Fix $(\ol t, \ol x)\in \dbR_T$. For any $\e>0$, there exists $\f_\e\in \ol \cU$ such that $\f_\e(\ol t, \ol x) <\ol u(\ol t, \ol x) +\e$.   By the structure of $\ol \cU$, it is clear that $\f_\e \ge \ol u$ on $\dbR_T$. Assume that $\f_\e\in \cU$ corresponds to the partition $0=t_0<\cds<t_n=T$ as in \eqref{Ubar}. 
We distinguish between two cases.
 
 {\it Case 1.} Assume $\ol t \in (t_{i-1}, t_i)$ for some $i=1$, $\ldots$, $n$. Since $\f_\e$ is continuous in $(t_{i-1}, t_i) \times \dbR$, there exists $\d>0$ such that $|\f_\e(t,x) - \f_\e(\ol t, \ol x)|\le \e$ whenever $|t-\ol t|+|x-\ol x|\le \d$. Then, for such $(t,x)$, 
\beaa
\ol u (t,x) \le \f_\e(t,x) \le \f_\e(\ol t, \ol x)+\e \le \ol u(\ol t, \ol x) + 2\e.
\eeaa
This implies that $\ol u$ is  upper semi-continuous  at $(\ol t, \ol x)$. 

 {\it Case 2.}  Assume $\ol t = t_i$ for some $i=0$, $\ldots$, $n$. By the same arguments as in Case 1, for any $\e>0$, there exists $\d>0$ such that $\ol u(t, x) \le  \ol u(\ol t, \ol x) + 2\e$ for all $(t,x) \in (t_i-\d, t_i]\times O_\d(\ol x)$. To see the regularity in the right neighborhood,  assume 
 for notational simplicity that $\ol t=0$. Let $\hat F$ be as in the proof of Lemma \ref{lem-cU}. Consider the following ODE (with parameter $x$):
 \begin{align}
 \label{ODEvHatF}
 v(t,x) = \f_\e(0, x) + \int_0^t \hat F(s,v(s, x)) ds.
 \end{align}
By the arguments in Subsection \ref{sect-1st}, there exists a $\d_0>0$ such that 
\eqref{ODEvHatF} has a classical solution $v \in C^{k,0}_{\a,\b}([0, \d_0]\times \dbR)$, which clearly leads to a classical supersolution $u\in C^k_{\a,\b}([0, \d_0]\times \dbR)$ of the original RPDE \eqref{RPDE} with initial condition $\f_\e(0, x)$.  Now consider  RPDE \eqref{RPDE} on $[\d_0, T]$ with initial condition $u(\d_0, \cd)$. By the arguments in Lemma \ref{lem-cU}, there exists  a $\tilde \f_\e \in \ol \cU$ such that $\tilde \f_\e(t,x) = u(t,x)$ for $(t,x) \in [0, \d_0]\times \dbR$.  Then  $\ol u \le u$ on  $[0, \d_0]\times \dbR$. Now by the continuity of $u$, there exists a $\d\le \d_0$ such that, whenever $|t|+|x-\ol x|\le \d$,
\beaa
\ol u(t,x) \le u(t,x) \le u(0, \ol x) + \e = \f_\e(0, \ol x) +\e\le \ol u(0, \ol x) + 2\e,
\eeaa
implying the regularity in the right neighborhood, and thus $\ol u$ is  upper semicontinuous.

We finally show that $\ul u$ is a viscosity subsolution, provided its continuity, 
The viscosity supersolution property of $\ol u$ follows similar arguments. 
Fix $(t_0, x_0)\in (0, T]\times \dbR$. 
Let $\f\in \ul\cA_g \ul u(t_0, x_0)$. For any $\e>0$, let  $(\pa D_\e^-(t_0,x_0),\psi^\e)$ be as in \eqref{De-}-\eqref{stability-psie}. 
By  definition, 
there exists a $u_\e\in \ul\cU$ with $\ul u(t_0, x_0) -u_\e(t_0, x_0) \le - \psi^\e(t_0, x_0)$. Let $\f^\e := \f + \psi^\e$, 
$\pa D_\e^- := \pa D_\e^-(t_0,x_0)$. Then 
\begin{align*}
 \f^\e(t_0, x_0) &= \ul u(t_0, x_0) + \psi^\e(t_0, x_0) \le u_\e(t_0, x_0\text{ and} \\
 \inf_{(t,x)\in \pa D_\e^-} [\f^\e - u_\e] (t, x) &\ge  \inf_{(t,x)\in\pa D_\e^-} [\ul u + \psi^\e - u_\e] (t, x)\ge  \inf_{(t,x) \in\pa D_\e^-} \psi^\e(t, x) >0.
\end{align*}
follow. Thus there exists a $(t_\e, x_\e) \in (t_0-\e^3, t_0] \times O_\e(x_0)$ such that
\beaa
[\f^\e - u_\e] (t_\e, x_\e) = 0 =  \inf_{(t,x)\in [t_0-\e^3, t_\e]\times O_\e(x_0)} [\f^\e - u_\e] (t, x).
\eeaa
Then $\f^\e \in  \ul\cA_g u_\e (t_\e, x_\e)$. Since $u_\e$ is a classical subsolution, hence a viscosity subsolution,  $\cL \f^\e(t_\e,x_\e) \le 0$. Sending $\e\to 0$ yields $\cL\f(t_0,x_0) \le 0$, i.e.,
 $\ul u$ is a viscosity subsolution.
\qed

\begin{thm}
\label{thm-comparison}
Let Assumptions~\ref{assum-g1}, \ref{assum-f}, and \ref{assum-u0}  hold. 
Let \eqref{boundary} be in force with  $u_1(0,\cd) \le u_0 \le u_2(0,\cd)$. Assume further that
\bea
\label{olu=ulu}
\ol u = \ul u.
\eea
Then $u_1 \le \ol u = \ul u \le u_2$ and $\ol u$  is the  unique   viscosity solution  of RPDE \eqref{RPDE}.
\end{thm}
\proof By Lemma \ref{lem-oluvis} and \eqref{olu=ulu}, it is clear that $\ol u = \ul u$  is continuous
 and  is a viscosity solution  of RPDE \eqref{RPDE}.   By Theorem \ref{thm-partial}
 (partial comparison),  $u_1 \le \ol u$ and $\ul u \le u_2$. Thus \eqref{olu=ulu} leads to the comparison principle immediately.
\qed 

 \begin{rem}
  \label{rem-perron}
{\rm   The introduction of $\ol u$ and $\ul u$ is motivated from  Perron's approach in PDE viscosity theory. However, there are several differences.
   
   (i) In Perron's approach, the functions in $\ol\cU$ are viscosity supersolutions, rather than classical supersolutions. So our $\ol u$ is in principle larger than the counterpart in PDE theory. Similarly our $\ul u$ is smaller than the  the counterpart in PDE theory. Consequently,  it is more challenging to verify the condition  \eqref{olu=ulu}. 
  
  (ii) The standard Perron's approach is mainly used for the existence of viscosity solution
  in the case the PDE satisfies the comparison principle. Here we use $\ol u$ and $\ul u$ to prove both the comparison principle and the existence. 
  
  (iii) In  the standard Perron's approach, one shows directly that $\ol u$  is a viscosity solution, while in Lemma \ref{lem-oluvis} we are only able to show $\ol u$ is a viscosity supersolution.
 }
  \end{rem}  

The condition \eqref{olu=ulu} is  in general  quite challenging. In the next section, we establish the complete result  when the diffusion coefficient $g$ is semilinear.

\section{Rough PDEs with semilinear diffusion}
\label{sect-semilinear}
We study RPDE \eqref{RPDE} and PPDE \eqref{PPDE}  in the case that  $g$ is semilinear, i.e.,
\bea
\label{gsemilinear}
g(t,x,y,z)=\sigma(t,x)\,z+g_0(t,x,y).
\eea

We  employ the following assumption.
\begin{assum}
\label{assum-g2}
$\si \in C^{k_0}_{\a,\b}(\dbR_T)$ and $g_0\in C^{k_0}_{\a,\b}(\dbR_T^2)$ for some large $k_0$. 
\end{assum}

 Clearly, Assumption \ref{assum-g2} implies Assumption \ref{assum-g1}. Note that in this section,
  we  obtain global result. Thus we require $g_0$ and its derivatives are uniformly bounded in $y$ as well.

\subsection{Global equivalence with the PDE}

Here, \eqref{E:Char} becomes
\begin{equation}
\label{characteristic1}
\begin{split}
 X_t(x) &= x-\int_0^t \sigma(s,X_s(x))\,d\o_s,\\Y_t(x,y)&=y+\int_0^t g_0\big(s,X_s(x),Y_s(x,y)\big)\,d\o_s,
 \end{split}
 \end{equation}
where $X$ (resp.~$Y$) depends only on $x$ (resp.~$(x,y)$),
and 
 \begin{align*}
Z_t(\th) &= z+\int_0^t \Big[Z_s(\th) \big[\pa_x\sigma(s,X_s(x))+\pa_y g_0\big(s,X_s(x),Y_s(x,y)\big)\big]
\\&\qquad\qquad\qquad +\pa_x g_0\big(s,X_s(x),Y_s(x,y)\big) \Big]\,d\o_s.
\end{align*}
where $\th=(x,y,z)$.
By Lemma \ref{lem-RDElinear}, we have, omitting the variable $\th$,
\begin{align}
\label{Z1}
Z_t = \G_t z + \int_0^t   {\G_t\over \G_s}\pa_x g_0(s,X_s,Y_s)\, d\o_s,\,
\end{align}
where $\G_t:= 
 \exp\big(\int_0^t  [\pa_x\sigma(s,X_s)+\pa_y g_0(s,X_s,Y_s)]\,d\o_s\big)$.

\begin{lem}
\label{lem-characteristic1}
Let Assumption \ref{assum-g2} hold.

(i)  RDE \eqref{characteristic1} has a classical solution $(X, Y)$ satisfying 
\bea
\label{XYest}
X - x \in  C^k_{\a,\b}(\dbR_T),\q Y-y \in C^k_{\a,\b}(\dbR_T^2)
\eea

(ii) There exists a $c>0$ such that 
\bea
\label{paXY}
\left.\ba{c}
\dis \pa_x X_t(x) = \exp\Big(-\int_0^t \pa_x \si(s, X_s(x)) \,d\o_s\Big) \ge c,\\
\dis \pa_y Y_t(x,y) = \exp\Big(\int_0^t \pa_y g_0(s, X_s(x), Y_s(x, y)) \,d\o_s\Big)\ge c.
\ea\right. 
\eea

(iii) For each $t$, the mapping $x \mapsto X_t(x)$ has inverse function $X^{-1}_t(\cd)$; and for each $(t,x)$, the mapping $y \mapsto Y_t(x, y)$ has inverse function $Y^{-1}_t(x, \cd)$.
\end{lem}
We remark that the proof below uses \eqref{paXY}. One can also use the backward rough path in \eqref{back} to construct the inverse functions directly. This argument works in multidimensional settings as well (see \cite{KZ}).
\ms

\proof (i)  follows directly from Lemma \ref{lem-RDEx}, which also implies
\beaa
&\dis\pa_x X_t(x) = 1 - \int_0^t \pa_x \si(s, X_s(x)) \pa_x X_s(x) \,d\o_s;&\\
&\dis\pa_y Y_t(x, y) = 1 + \int_0^t \pa_y g_0\big(s, X_s(x), Y_s(x,y)) \pa_y Y_s(x,y\big) \,d\o_s.&
\eeaa
Then the representations in \eqref{paXY} follow from Lemma \ref{lem-RDElinear}.   Moreover, set $\check X := X-x\in C^k_{\a,\b}(\dbR_T)$ and  $\tilde \si(t, x):= \si(t, X_t(x)) = \si(t, x + \check X_t(x))$. Then by the uniform regularity of $\si$, 
$\sup_{x\in \dbR} \|\tilde\si(\cd, x)\|_k \le C$. This implies that $ \int_0^t \pa_x \si_s(X_s(x)) \pa_x X_s(x) \,d\o_s$ is uniformly bounded, uniformly in $(t,x)$. Therefore, we obtain the first estimate for $\pa_x X$ in \eqref{paXY}. The second estimate for $\pa_y Y$ in \eqref{paXY} follows from the similar arguments. 

Finally,  for each $t$, the fact $\pa_x X_t(x) \ge c$  implies that  $x\mapsto X_t(x)$ is one to one and the range  is the whole real line $\dbR$. Thus $X^{-1}_t: \dbR\to \dbR$ exists. Similarly one can show that $Y^{-1}_t(x, \cd)$ exists. 
\qed

One can easily check,  omitting  $(x, y,z)$ in $X_t(x), Y_t(x, y)$, $Z_t(x,y,z)$,
\beaa
&\pa_y X_t = 0;\q  \pa_z X_t = 0;\q  Z_t = {\pa_y Y_t \over \pa_x X_t} z + {\pa_x Y_t \over \pa_x X_t};  \q \pa_z Z_t = {\pa_y Y_t \over \pa_x X_t};&\\
&\pa_y Z_t =  {\pa_{yy} Y_t \over \pa_x X_t} z + {\pa_{x y}Y_t \over \pa_x X_t},~  \pa_x Z_t =  {\pa_{xy}Y_t\pa_x X_t - \pa_y Y_t \pa^2_{xx} X_t  \over (\pa_x X_t)^2} z + {\pa_{x x}Y_t \pa_x X_t - \pa_x Y_t \pa^2_{xx} X_t\over (\pa_x X_t)^2},&
\eeaa
and then \eqref{F} becomes 
\begin{align}
\label{F1} 
F(t,x,y,z,\g) := {1\over \pa_y Y_t}~ f\Big(t, X_t, Y_t, {\pa_y Y_t \over \pa_x X_t} z + {\pa_x Y_t \over \pa_x X_t},  
  {\pa_y Y_t \over (\pa_x X_t)^2}\g \qquad \\ 
 + {\pa_{yy} Y_t \over (\pa_x X_t)^2} z^2+ {2\pa_{xy} Y_t \pa_x X_t - \pa_y Y_t \pa^2_{xx} X_t\over (\pa_xX_t)^3} z+ {\pa^2_{xx} Y_t \pa_x X_t- \pa_x Y_t  \pa^2_{xx} X_t\over (\pa_xX_t)^3}\Big).\nonumber
\end{align}

Under our conditions, $F$  has typically quadratic growth in $z$ and is not uniformly Lipschitz 
 in $y$.  Moreover, the first equality of \eqref{E:vYZ}  becomes
\begin{equation}
\label{uvrep}
\begin{split}
v(t,x) = Y^{-1}_t\big(x, u(t, X_t(x))\big)\text{ or}~
u(t,x) = Y_t\big(X^{-1}_t(x), v(t, X^{-1}_t(x))\big).
\end{split}
\end{equation}

By using similar arguments as in Subsection \ref{sect-PDE}, we obtain the following result which is global in this semilinear case.
\begin{thm}
\label{thm-semilinear}
Let Assumptions \ref{assum-g2} and \ref{assum-f} hold. Assume
 $u\in C^k_{\a,\b}(\dbR_T)$ and $v\in C^{k,0}_{\a,\b}(\dbR_T)$ satisfy \eqref{uvrep}. Then $u$ is a classical solution (resp. subsolution,  supersolution)  of RPDE \eqref{RPDE}-\eqref{gsemilinear}   if and only if  $v$ is a classical solution (resp.  subsolution, supersolution) of  PDE \eqref{PDE}-\eqref{F1}.
\end{thm}

The next result establishes  equivalence in the viscosity sense. 
\begin{thm}
\label{thm-equiv1}
Let Assumptions \ref{assum-g2} and \ref{assum-f} hold. Assume $u$, $v\in C(\dbR_T)$ satisfy \eqref{uvrep}.  Then $u$ is a viscosity solution (resp.~subsol.,  supersol.)  of RPDE \eqref{RPDE}-\eqref{gsemilinear}  at $(t_0,x_0)\in (0,T]\times\dbR$  if and only if $v$  is a viscosity solution  (resp.~subsolution,  
supersolution) of PDE \eqref{PDE}-\eqref{F1} at $(t_0, X^{-1}_{t_0}(x_0))$.
\end{thm}
\proof We  prove  the statement only for supersolutions. First, we prove the if part.  Let $\tilde x_0 := X^{-1}_{t_0}(x_0)$ and $v$ be a  viscosity supersolution  of PDE \eqref{PDE}-\eqref{F1} at $(t_0, \tilde x_0)$.  Let $\f \in \overline{\mathcal{A}}_g u(t_0,x_0)$ with corresponding $\d_0$.  Define
\begin{align}
\label{psi1}
\psi (t, x) :=  Y^{-1}_t\big(x, \f(t, X_t(x))\big),\mbox{ i.e., } Y_t(x, \psi(t,x)) = \f(t, X_t(x)).
\end{align}
It is clear that $\psi(t_0, \tilde x_0) = v(t_0, \tilde x_0)$. By the continuity of $X$, there exists a 
$\d>0$ such that $X_t(x) \in O_{\d_0}(x_0)$ for all $(t,x) \in D_\d(t_0, \tilde x_0)$.  By the same arguments as for \eqref{paov1},  $\pa_\o \psi = 0$. Moreover, for $(t,x) \in D_\d(t_0, \tilde x_0)$, since $\f \in \overline{\mathcal{A}}_g u(t_0,x_0)$, we have $\f(t, X_t(x)) \le u(t, X_t(x))$. By Lemma \ref{lem-characteristic1}, the mapping $y\to Y_t(x,y)$ is increasing. Thus $Y^{-1}$ is also increasing and  $\psi(t,x) \le v(t,x)$,  i.e., $\psi$ is a test function for $v$ at $(t_0, x_0)$ and 
\bea
\label{psiviscosity1}
\pa_t \psi(t_0, \tilde x_0) \ge  F\big(t_0, \tilde x_0, \psi(t_0,\tilde x_0), \pa_x \psi(t_0, \tilde x_0), \pa^2_{xx}\psi(t_0, \tilde x_0)\big).
\eea
By the derivation of $F$, this implies
\bea
\label{phiviscosity1}
\pa^\o_t \f(t_0, x_0) \ge  f\big(t_0, x_0, \f(t_0,x_0), \pa_x \f(t_0, x_0), \pa^2_{xx}\f(t_0, x_0)\big),
\eea
i.e., $u$ is a viscosity supersolution at $(t_0, x_0)$.

For the opposite direction, assume that $u$ is a  viscosity supersolution  of RPDE \eqref{RPDE} at $(t_0, x_0)$. For $\psi \in \overline{\mathcal{A}}_0 v(t_0,\tilde x_0)$ corresponding to $g=0$, define
$\f (t, x) :=  Y_t\big(X^{-1}_t(x), \psi(t, X^{-1}_t(x))\big)$,
which still implies $Y_t(x, \psi(t,x)) = \f(t, X_t(x))$. By similar arguments,  \eqref{psiviscosity1} 
follows from \eqref{phiviscosity1}. 
\qed

\begin{rem}
\label{rem-visPDE}
{\rm In the general case,  there are two major differences:

(i) The transformation determined by \eqref{E:vYZ} involves both $u$ and $\pa_x u$, i.e.,
 to extend Theorem \ref{thm-equiv1}, one has to assume that the candidate viscosity solution $u$ is differentiable in $x$.

(ii) The transformation is local, in particular, the $\d$ in Theorem \ref{thm-utov}  depends on $\|\pa^2_{xx} u\|_\infty$, i.e.,  unless $\pa^2_{xx} u$  is bounded and  the solution is essentially classical, we have difficulty to extend Theorem \ref{thm-equiv1} to the general  case, even just in local sense.
}
 \end{rem}

\subsection{Some a priori estimates}
Here, we establish  uniform a priori estimates for $v$ that will be crucial for the comparison principle of viscosity solutions in the next subsection. First, we estimate  the $\dbL^\infty$-norm  of $v$.
\begin{prop}
\label{prop-uvbound}
Let Assumptions \ref{assum-g2}, \ref{assum-f}, and \ref{assum-u0} hold and $f$ be smooth. Assume  further that $v\in C^{k,0}_{\a,\b}(\dbR_T)$ is a classical solution of   PDE \eqref{PDE}-\eqref{F1}. Then there exists a constant $C$, which depends only on the constants $K_0$,  $L_0$ in Assumption \ref{assum-f},  and the regularity of $\si, g_0$ in Assumption \ref{assum-g2}, but does not depend on $u_0$ or the further regularity  of $f$, such that
\bea
\label{vbound}
 |v(t,x)| \le e^{Ct} \big[\|u_0\|_\infty + Ct\big].
\eea
\end{prop}
\proof First, we write \eqref{PDE}-\eqref{F1} as
\begin{align*}
&\pa_t v =  a(t,x) \pa^2_{xx} v + b(t,x) \pa_x v +F(t,x,v, 0, 0), \q\mbox{where} \\
&a(t,x) := \tilde a(t,x, v(t,x), \pa_x v(t,x), \pa^2_{xx} v(t,x)),\\
&b(t,x) := \tilde b(t,x, v(t,x), \pa_x v(t,x), \pa^2_{xx} v(t,x)), \\
& \tilde a(t,x,y,z,\g) := \int_0^1 \pa_\g F(t,x,y,\l z, \l \g) d\l,\\
& \tilde b(t,x,y,z,\g) := \int_0^1 \pa_z F(t,x,y,\l z, \l \g) d\l.
\end{align*}
Since $v$ is a classical solution, $a$ and $b$ are smooth functions. Reversing the time by 
setting $\hat \f(t,x) := \f(T-t, x)$ for $\f = v$, $a$, $b$, $F$,  we have
\beaa
\pa_t \hat v +  \hat a(t,x) \pa^2_{xx}\hat v + \hat b(t,x) \pa_x\hat v +\hat F(t,x,\hat v, 0, 0) =0,\q \hat v(T,x) = u_0(x).
\eeaa
Let $B$ be a standard Brownian motion.  Consider  the  SDE
\bea
\hat X_t = x + \int_0^t \hat b(s, \hat X_s) ds + \int_0^t \sqrt{2\hat a}(s, \hat X_s) dB_s.
\eea
Then $\hat Y_t := \hat v(t, \hat X_t)$ solves the BSDE
\beaa
\hat Y_t = u_0(\hat X_T) + \int_t^T \hat F(s, \hat X_s, \hat Y_s,0,0) ds - \int_t^T \hat Z_s dB_s.
\eeaa
Since 
$F(t,x,y,0,0) = {1\over \pa_y Y_t}f\big(t, X_t, Y_t, {\pa_x Y_t \over \pa_x X_t},  ~  {\pa^2_{xx} Y_t \pa_x X_t- \pa_x Y_t  \pa^2_{xx} X_t\over (\pa_xX_t)^3}\big)$, we have
\bea
\label{F0est}
|F(t,x,y, 0, 0)|\le C[1+ |y|],
\eea
following from Lemma \ref{lem-characteristic1}. Then, by  standard BSDE estimates,  
\beaa
|v(T,x)| = |\hat v(0, x)| = |\hat Y_0| \le e^{CT} \big[\|u_0\|_\infty +  CT \big].
\eeaa
This implies \eqref{vbound} for $t=T$. Similarly we may prove \eqref{vbound} for all $t>0$.  
\qed

\begin{rem}
\label{rem-blowup}
{\rm (i) We are not able to establish similar a priori estimates for $\pa_x v$.  Besides the possible insufficient regularity of $u_0$, we emphasize that the main difficulty here is not that $F$ has quadratic growth in $z$, but that $F$ is not uniformly Lipschitz continuous in $y$. Nevertheless,  we  obtain some local estimate for  $\pa_x v$ in Proposition \ref{prop-semilinear-local}, which will be crucial for the comparison principle of viscosity solutions later.

(ii)  To overcome the  difficulty above and apply standard techniques, 
\cite{LS3} imposed  technical conditions on $f$ in the case $f= f(z, \g)$ (cf.~\cite[(1.12)]{LS3}):
\begin{align}
\label{LScondition}
 \g \pa_\g f + z \pa_z f - f ~\mbox{is either bounded from above or from below},
 \end{align} 
This is  essentially satisfied when $f$ is convex or concave in $(z, \g)$.  Our $f$ in \eqref{fIsaacs} below does not satisfy \eqref{LScondition}, in particular,  we do not require $f$ to be convex or concave in $z$. See also Remark \ref{rem-Isaacs}. 
 }
\end{rem}

The next result relies on representation of $v$ and BMO estimates for BSDEs  with quadratic growth. For this purpose, we restrict $f$ to Bellman-Isaacs type with the  Hamiltonian
\begin{equation}
\label{fIsaacs}
\begin{split}
 f(t,x,y,z,\g) &= \sup_{e_1 \in E_1} \inf_{e_2\in E_2}  \Big[{1\over 2} \si^2_f(t, x, e) \g + b_f(t, x, e) z \\  
 &\qquad\qquad+f_0\big(t,x,y,   \si_f(t, x,e) z, e\big)\Big],
 \end{split}
 \end{equation}
where $E:= E_1 \times E_2 \subset \dbR^2$ is the control set and $e=(e_1, e_2)$.

\begin{assum}
\label{assum-f1}
(i) $\si_f$, $b_f \in C^0(\dbR_T\times E)$ are bounded by $K_0$, uniformly Lipschitz continuous in $x$ with Lipschitz constant $L_0$, and $\si_f\ge 0$;

(ii) $f_0\in C^0(\dbR_T^3\times E)$ is uniformly Lipschitz continuous in $(x, y,z)$  with Lipschitz constant $L_0$,  and $f_0(t,x,0,0,e)$ is bounded by $K_0$.
\end{assum}

Assumption \ref{assum-f1} obviously implies Assumption \ref{assum-f}.

\begin{prop}
\label{prop-semilinear-local}
Let Assumptions \ref{assum-g2},    \ref{assum-f1}, and \ref{assum-u0} hold, and $(g, f)$ take the form \eqref{gsemilinear}-\eqref{fIsaacs}.   Assume  $v\in C^{k,0}_{\a,\b}(\dbR_T)$ is a classical solution of   PDE \eqref{PDE}-\eqref{F1}. Then there exist  constants $\d_0>0$ and $C_0$,  which depend only on $K_0$, $L_0$ in Assumption \ref{assum-f1},  the regularity of $\si$, $g_0$ in Assumption  \ref{assum-g2},  and  $\|u_0\|_\infty$, but  not  on the further regularity  of $f$ and $u_0$, such that
\bea
\label{palest}
|\pa_x v(t,x)|\le C_0[1+ \|\pa_x u_0\|_\infty]\q\mbox{for all}~(t,x) \in [0, \d_0]\times \dbR.
\eea
\end{prop}
\proof Under \eqref{gsemilinear} and \eqref{fIsaacs}, \eqref{F} and the equivalent \eqref{F1}  becomes
\begin{equation}
\label{FIsaacs}
\begin{split}
F(t,x,y,z,\g) &= \sup_{e_1\in E_1} \inf_{e_2\in E_2} \Big[{1\over 2} \hat \si^2_f(t, x, e) \g + \hat b_f(t,x,e) z\\
&\qquad\qquad + F_0\big(t,x,y, \hat \si_f(t,x,e) z,e\big)\Big], 
\end{split}
\end{equation}
where, omitting $(x, y)$ in$X_t(x)$ and $Y_t(x, y)$, 
\begin{align*}
&\hat \si_f(t,x,e) :={ \si_f(t, X_t, e)\over \pa_x X_t},\q \hat b_f(t,x,e) :=  { b_f(t, X_t, e)\over \pa_x X_t},\\
&F_0(t,x,y,z,e) :=   {1\over 2}\hat \si^2_f(t, x, e) {\pa^2_{yy} Y_t \over \pa_y Y_t } z^2 +\hat \si^2_f(t, x, e) [{\pa^2_{xy} Y_t \over \pa_y Y_t} - { \pa^2_{xx} X_t\over  2\pa_x X_t}]  z\\
& \qq +  {1\over \pa_y Y_t} f_0\big(t,X_t, Y_t,   \pa_y Y_t  z + \hat\si_f (t, x, e) \pa_x Y_t, e\big)\\
& \qq  +\hat \si^2_f(t, x, e)  {\pa^2_{xx} Y_t \pa_x X_t - \pa_x Y_t \pa^2_{xx} X_t \over  2 \pa_y Y_t \pa_xX_t}+\hat b_f(t, x, e) {\pa_x Y_t \over \pa_y Y_t}.
\end{align*}
By \eqref{paXY}, we have,  again omitting $(x, y)$ in $X_t(x), Y_t(x,y)$, 
\begin{align*}
&\pa^2_{xx} X_t = -\pa_x X_t\int_0^t \pa^2_{xx} \si(s, X_s)  \pa_x X_s d\o_s,\\
&\pa^2_{yy} Y_t =\pa_y Y_t \int_0^t \pa^2_{yy} g_0(s, X_s, Y_s) \pa_y Y_s d\o_s,\\
&\pa^2_{xy} Y_t = \pa_y Y_t \int_0^t [\pa^2_{xy} g_0(s, X_s, Y_s)\pa_x X_s + \pa^2_{yy} g_0(s, X_s, Y_s) \pa_x Y_s] d\o_s.&
\end{align*}
Then, by \eqref{XYest} we can easily verify that
\begin{equation}
\label{paF0est}
\begin{split}
\left.\ba{c}
\dis \mbox{$\hat\si_f$, $\hat b_f$, and $F_0(\cd, 0, 0,\cd)$ are bounded,}\q |\pa_x \hat\si_f |\le C,\q |\pa_x \hat b_f|\le C,\\
\dis |\pa_z F_0(t,x,y,z)| \le C[1+\rho(t) |z|],\\
 \dis |\pa_x F_0(t,x,y,z)| + |\pa_y F_0(t,x,y,z)|\le C[1+|y|+|z|+ \rho(t)|z|^2].
\ea\right.
\end{split}
\end{equation}
where $\rho\ge 0$ is a continuous function with $\rho(0)=0$. Here for notational simplicity we are assuming the relevant functions are differentiable, but actually we only need their uniform Lipschitz continuity.

Now, let $\bar B$ be a standard Brownian motion and $\cE = \cE_1\times \cE_2$ be the set of $\dbF^{\bar B}$-progressively measurable $E$-valued processes. Fix  $\d>0$ and define  $\bar \f(t,x,y,z,e) := \f(\d-t, x, y,z, e)$ for $\f = \hat\si_f, \hat b_f, F_0$. For any $e\in \cE$, introduce the following decoupled FBSDE on $[0, \d]$:
\bea
\label{cXY}
\left.\ba{lll}
\dis \cX^e_t = x +   \int_0^t \bar b_f(s, \cX^e_s, e_s) ds +\int_0^t \bar \si_f(s, \cX^e_s, e_s) d\bar B_s;\\
\dis \cY^e_t = u_0(\cX^e_\d) + \int_t^\d \bar F_0(s, \cX^e_s, \cY^e_s, \cZ^e_s, e_s) ds - \int_t^\d \cZ^e_s d\bar B_s.
\ea\right.
\eea
By  Zhang \cite[Theorems 7.2.1, 7.2.3]{Zhang}, there exist constants $c_0$, $C_0$, depending on $\|u_0\|_\infty$ and $\|f(\cd, 0,0,\cd)\|_\infty$ (the bound of $|f(t,x,0,0,e)|$), such that 
\bea
\label{cZest}
\|\cY^e\|_\infty \le C_0,\q \dbE\Big[\exp\big(c_0\int_0^\d |\cZ^e_s|^2 ds\big)\Big] \le C_0 <\infty.
\eea
Differentiating \eqref{cXY} with respect to $x$ yields
\begin{align*}
 \dis \td \cX^e_t &= 1 +   \int_0^t \pa_x \bar b_f \td \cX^e_s ds +\int_0^t \pa_x \bar \si_f\td \cX^e_s d\bar B_s;\\
\dis \td \cY^e_t &= \pa_x u_0(\cX^e_\d)\td \cX^e_\d + \int_t^\d \Big[\pa_x \bar F_0 \td \cX^e_s + \pa_y \bar F_0 \td \cY^e_s + \pa_z \bar F_0 \cZ^e_s\Big]  ds \\ &\qquad - \int_t^\d \td \cZ^e_s d\bar B_s.
\end{align*}
This implies
\begin{align*}
 \td \cY^e_0 = \dbE\Big[ \G^e_\d  \pa_x u_0(\cX_\d)\td \cX^e_\d + \int_0^\d \G^e_t \pa_x \bar F_0 \td \cX^e_t dt \Big]
\end{align*}
where $\G^e_t:=  \exp\big(\int_0^t  \pa_z \bar F_0 d\bar B_s + \int_0^t [\pa_y \bar F_0 - {1\over 2} |\pa_z \bar F_0|^2] ds \big)$.
By \eqref{paF0est},
\begin{align*}
&\dbE[ |\G^e_t|^4] =  \dbE\Big[\exp\big(4\int_0^t  \pa_z \bar F_0 d\bar B_s + \int_0^t [4\pa_y \bar F_0 - 2 |\pa_z \bar F_0|^2] ds \big)\Big] \\
&= \dbE\Big[\exp\big(4\int_0^t  \pa_z \bar F_0 d\bar B_s  - 16 \int_0^t |\pa_z \bar F_0|^2 ds + \int_0^t [4\pa_y \bar F_0 +14 |\pa_z \bar F_0|^2] ds \big)\Big]\\
&\le \Big(\dbE\Big[\exp\big(8\int_0^t  \pa_z \bar F_0 d\bar B_s  - 32\int_0^t |\pa_z \bar F_0|^2 ds  \big)\Big]\\
&\qquad\qquad\times\dbE\Big[\exp\big( \int_0^t [8\pa_y \bar F_0 +28 |\pa_z \bar F_0|^2]  ds  \big)\Big]\Big)^{1\over 2}\\
&=\Big(\dbE\Big[\exp\big( \int_0^t [8\pa_y \bar F_0 +28 |\pa_z \bar F_0|^2]  ds  \big)\Big]\Big)^{1\over 2}\\
&= \Big(\dbE\Big[\exp\big( C\int_0^\d [ 1+ |\cY^e_s|+ |\cZ^e_s| + \rho(t) |\cZ^e_s|^2 + \rho(t)^2 |\cZ^e_s|^2 ]  ds  \big)\Big]\Big)^{1\over 2}
\end{align*}
Set $\d_0>0$ small enough so that $C[\rho(\d_0) + \rho(\d_0)^2] \le {c_0\over 2}$. Then, for $\d\le \d_0$, by \eqref{cZest}, we obtain  $\dbE\Big[ |\G^e_t|^4\Big] \le C_0$, and by the second line of    \eqref{paF0est} it is clear that $\dbE\Big[\sup_{0\le t\le \d} |\td \cX^e_t|^4\Big] \le C_0$. Thus
\begin{equation}
\label{tdYe}
\begin{split}
|\td \cY^\e_0|&\le C_0 \dbE\Big[\|\pa_x u_0 \|_\infty |\G^e_\d| |\td \cX^e_\d|\\
& +  \int_0^\d |\G^e_t||\td \cX^e_t|[1+|\cY^e_t| +|\cZ^e_t|^2]dt \Big]
\le C_0[1+\|\pa_x u_0\|_\infty].
\end{split}
\end{equation}

Finally, we remark that,  since we know a priori that  $v\in C^{k,0}_{\a,\b}(\dbR_T)$, by the standard truncation arguments, we may assume without loss of generality that $F_0$ is uniformly Lipschitz 
in $(x,y,z)$ (with the Lipschitz constant possibly depending on the regularity of $v$). Then, by Buckdahn \& Li  \cite{BL},
\bea
\label{game}
v(\d,x) = \inf_{\cS}\sup_{e_1 \in \cE_1} \cY^{(e_1, \cS(e_1))}_0,
\eea
where $\cS: \cE_1 \to \cE_2$ is the so called nonanticipating strategies.  This implies 
\beaa
|\pa_x v(\d, x)| \le \sup_{\cS}\sup_{e_1 \in \cE_1}  |\td \cY^{(e_1, \cS(e_1))}_0| \le C_0[1+\|\pa_x u_0\|_\infty].
\eeaa
Since $\d\le \d_0$ is arbitrary, the proof is complete.
\qed

\begin{rem} 
\label{rem-doubly}
{\rm  (i) We reverse the time in \eqref{cXY}. Hence, in  spirit of the backward rough path in \eqref{cXY},  $\ol B$ and the rough path $\o$ (or the original $B$ in \eqref{SPDE}) have opposite directions of time evolvement. Thus \eqref{cXY} is in the line of the backward doubly SDEs of Pardoux  \& Peng 
\cite{PP94BDSDE}. When $E_2$ is a singleton,  
\cite{MPS} provided a representation for the corresponding SPDE \eqref{SPDE}  in the context of second order backward doubly SDEs. We shall remark though, while the wellposedness of backward doubly SDEs holds true for random coefficients, its representation for solutions of SPDEs requires Markovian structure, i.e., the $f$ and $g$ in \eqref{SPDE} depend only on $B_t$ (instead of the path $B_\cd$). The stochastic characteristic approach used in this paper does not have this constraint. Note again that our $f$ and $g$ in RPDE \eqref{RPDE} and PPDE \eqref{PPDE} are allowed to depend on the (fixed) rough path $\o$.

(ii) For \eqref{game}, from a game theoretical point of view, it is more natural to use the 
so-called weak formulation (see Pham \& Zhang \cite{PZ}). However, as we are here mainly concerned about the regularity,  the strong formulation in \cite{BL} is more convenient.
} 
\end{rem}

\subsection{The global comparison principle and existence of viscosity solution}

We need   the following PDE result from Safonov \cite{Safonov88parabolic} (see also  
\cite{MP94} 
 for a corresponding statement
for bounded domains and 
\cite{Safonov89elliptic} for the elliptic case).

\begin{thm}
\label{thm-Safonov}
Consider PDE \eqref{PDE}. Assume that, for some $\b>0$, 

(i) $F$ is convex in $\g$ and uniformly parabolic, i.e., $\pa^2_{\g\g} F \ge 0$
and  $\pa_\g  F \ge c_0 >0$,

(ii) $F$ is uniformly Lipschitz continuous in $(y,z, \g)$,

(iii) $\|F(\cd, y,z,\g)\|_{C_b^{\b}(\dbR_T)} \le C[1+|y|+z|+|\g|]$,

(iv) $u_0 \in C^{2+\b}_b(\dbR)$.

\no Then there exists  $\beta_0\in (0,1)$ depending only on $c_0$ such that, whenever $\beta\in (0, \beta_0]$, PDE   \eqref{PDE} has a classical solution $ v\in C^{2+\b}_{b}(\dbR_T)$.
\end{thm}

\begin{thm}
\label{thm-semilinear-classical}
Let Assumptions \ref{assum-g2},    \ref{assum-f1}, and \ref{assum-u0} hold, and $(g,f)$ take the form \eqref{gsemilinear}-\eqref{fIsaacs}.   Assume further that

(i) $f$ is either convex or concave in $\g$, namely either $E_1$ or $E_2$ in \eqref{fIsaacs}  is a singleton,

(ii) $\si_f \ge c_0>0$,

(iii) $u_0\in  C_b^{k+1+\b}(\dbR)$ and $f\in C^{k+1,loc}_{\a,\b}(\dbR_T^4)$.

\no   Then there exists $\d_0>0$, depending on $K_0$, $L_0$ in Assumption \ref{assum-f1},  the regularity of $\si$, $g_0$ in Assumption  \ref{assum-g2},  and  $\|u_0\|_\infty$, but independent of the further regularity of $u_0$ and $f$, such that PDE \eqref{PDE}-\eqref{F1} has a classical solution  $v\in C^{k,0}_{\a,\b}([0, \d_0]\times \dbR)$.
\end{thm}
\proof  We  prove only the  convex case, i.e., $E_2$ is a singleton. When $f$ is concave,  
one can use following the standard transformation: $\tilde f(t,x,y,z,\g) := -f(t, x, -y, -z, -\g)$ is convex and $\tilde v(t,x) := -v(t,x)$ corresponds to $\tilde f$.   Let $\d_0$ be determined by Proposition \ref{prop-semilinear-local}.

First it is clear that  the $F$ in \eqref{F1} (or the equivalent \eqref{FIsaacs})  satisfies the requirements in Theorem \ref{thm-Safonov} (i). Recall \eqref{vbound},  \eqref{palest}, and the $K_0$ in Assumption \ref{assum-f1} (ii).  Put
\beaa
C_1:= e^{CT}[\|u_0\|_\infty + CT K_0] + C_0[1+\|\pa_x u_0\|_\infty].
\eeaa
Introduce a truncation function $\iota\in C^\infty(\dbR)$ such that $\iota(x) = x$ for $|x|\le C_1$, and $\iota(x) = 0$ for $|x|\ge C_1+1$. Define
\beaa
\tilde F(t,x,y,z,\g) := F(t,x, \iota(y), \iota(z), \g).
\eeaa
Then one may verify straightforwardly that $\tilde F$ satisfies all the conditions in Theorem \ref{thm-Safonov}. Thus the following PDE has a classical solution 
$ \tilde{v}\in C^{2+(\beta\wedge\beta_0)}_b(\dbR_T)$:
\bea
\label{PDEtruncate}
\pa_t \tilde v=\tilde F(t,x,\tilde v, \pa_x \tilde v, \pa^2_{xx} \tilde v),\q \tilde v(0,\cd) = u_0.
\eea
Applying Propositions \ref{prop-uvbound} and  \ref{prop-semilinear-local} on the above PDE
yields  $|\tilde v|\le C_1$, $|\pa_x \tilde v|\le C_1$ on $[0, \d_0] \times \dbR$, i.e., 
$v:= \tilde v$ solves  PDE \eqref{PDE}-\eqref{F1} on  $[0, \d_0] \times \dbR$. 

Finally, the further regularity of $v$ follows from  standard bootstrap arguments
(cf.~\cite[Lemma~17.16]{GT})  together with Remark~\ref{rem-Cx}. 
\qed

\begin{rem}
\label{rem-Isaacs}
{\rm The requirement that $f$ is convex or concave is mainly to ensure the existence of classical solutions for PDE \eqref{PDEtruncate}. 
Theorem \ref{thm-Safonov} holds true for multidimensional case as well. When the dimension of $x$ is $1$ or $2$,  Bellman-Isaacs equations may have classical solution as well, see 
\cite[Theorem 14.24]{Lieberman} for $d=1$ and 
\cite[Lemma 6.5]{PZ} for $d=2$ for bounded domains, and also 
\cite[Theorem 17.12]{GT} for elliptic equations in bounded domains when  $d=2$. We believe such results can be 
extended to the whole space and thus the theorem above as well as Theorem \ref{thm-comparisonSemilinear}  will hold true when $f$ is indeed Bellman-Isaacs type. However, 
 when the dimension is high, the Bellman-Isaacs equation does, in general, not have a classical solution
 (see 
 \cite{NV} for a counterexample).
}
\end{rem}

\begin{thm}
\label{thm-comparisonSemilinear}
Let $(g, f)$ take the form \eqref{gsemilinear}-\eqref{fIsaacs}. Let Assumptions \ref{assum-g2},    \ref{assum-f1}, and \ref{assum-u0} hold. Assume 
that, for any $\e>0$, there exist $\ol f^\e$, $\ul f^\e$ such that

(i)  $\ol f^\e$, $\ul f^\e$ satisfy Assumption  \ref{assum-f1} uniformly, i.e., with the same $K_0, L_0$ for all $\e>0$,

(ii) for each $\e>0$, $\ol f^\e, \ul f^\e$ satisfy all the requirements in Theorem \ref{thm-semilinear-classical}, 

(iii) for each $\e>0$, $f-\e \le \ul f^\e \le f\le \ol f^\e \le f+\e$. 

\no Then RPDE \eqref{RPDE} satisfies the comparison principle and has a unique viscosity solution. 
\end{thm}
 \proof   By Lemma \ref{lem-oluvis}, $\ol u$ and $\ul u$ are bounded by some $C_0$.  
  
 {\it Step 1.} We  prove first \eqref{olu=ulu} locally.  Let $\d_0 >0$ be determined by Proposition \ref{prop-semilinear-local}, corresponding to $K_0$, $L_0$, but with $\|u_0\|_\infty$ replaced with the global bound $C_0$ of $\ol u$ and $\ul u$. 
 For any $\e>0$, let  $\ol f^\e$, $\ul f^\e$ be as in the assumption of the theorem, and $\ol F^\e$, $\ul F^\e$  correspond to $\ol f^\e$, $\ul f^\e$ as in \eqref{FIsaacs}.   In the spirit of Remark \ref{rem-poly} (i),  we may assume without loss of generality that $u_0$ is uniformly continuous. Then $u_0$ has  standard smooth mollifiers $\ol u^\e_0$, $\ul u^\e_0$  such that $u_0 - \e \le \ul u_0^\e \le u_0 \le \ol u^\e_0 \le u_0+\e$. By Theorem \ref{thm-semilinear-classical}, let $\ol v^\e$ (resp.~$\ul v^\e$) be the classical solution to PDE \eqref{PDE} -\eqref{FIsaacs} with coefficients $(\ol F^\e, g)$ and initial condition $\ol u^\e_0$   
 (resp.~coefficients $(\ul F^\e, g)$ and initial condition $\ul u^\e_0$) on $[0, \d_0]$.  Then, by \eqref{ubar},
 $\ul v^\e \le \ul v \le \ol v \le \ol v^\e$, where $\ul v := Y^{-1}_t\big(x, \ul u(t, X_t(x))\big)$ as in \eqref{uvrep}, and similarly for  $\ol v$. By \eqref{F} it is clear that $0 \le \ol F^\e - \ul F^\e \le C\e$. Define $\D v^\e := \ol v^\e - \ul v^\e$, $\D u^\e_0 := \ol u^\e_0 - \ul u^\e_0$,  $\D F^\e := \ol F^\e - \ul F^\e$. Then
 \begin{align*}
 &\D v^\e = \pa_t \D v^\e + F^\e(t,x,  \D v^\e, \pa_x\D v^\e, \pa^2_{xx} \D v^\e), \q \D v^\e(0,\cd) = \D u^\e_0,\\
&\text{where }
 F^\e(t,x,y,z,\g) := \ol F^\e(t,x, \ul v^\e + y, \pa_x \ul v^\e + z, \pa^2_{xx} \ul v^\e + \g)\\
 &\qquad\qquad\qquad\qquad\qquad -  \ul F^\e(t,x, \ul v^\e, \pa_x \ul v^\e, \pa^2_{xx} \ul v^\e).
 \end{align*}
 Now following the arguments of Proposition \ref{prop-uvbound}, we see that there exists a constant $C$, independent of $\e$, such that, for every $(t,x) \in [0, \d_0]\times \dbR$,
 \beaa
\ol v(t,x) - \ul v(t,x)\!\!\! &\le&\!\!\! \D v^\e(t,x) \le C^{Ct}\Big[\|\D v^\e_0\|_\infty + Ct \|F^\e(\cd, 0,0,0)\|_\infty\Big] \\
 \!\!\! &\le&\!\!\!  C\Big[\|\D u^\e_0\|_\infty + \|\D F^\e(\cd, \ul v^\e, \pa_x \ul v^\e, \pa^2_{xx} \ul v^\e)\|_\infty\Big] \le C\e.
 \eeaa
 This implies that  $\ol v(t,x) = \ul v(t,x)$. Thus 
 \eqref{olu=ulu} holds on $[0, \d_0]$. Therefore, by  Theorem  \ref{thm-comparison}, 
  $u_1 \le \ol u\le u_2$   and $u:= \ol u$ is the unique viscosity solution of RPDE \eqref{RPDE}-\eqref{fIsaacs} on $[0, \d_0]$.

{\it Step 2.} We  prove now the global result.  Let $0=t_0<\cds<t_n=T$ be such that $t_i-t_{i-1}\le \d_0$ for each $i=1$, $\ldots$, $n$. By Step~1,  $u_1(t_1,\cd)\le u(t_1, \cd) \le u_2(t_1,\cd)$. Now consider RPDE \eqref{RPDE}-\eqref{fIsaacs}   on $[t_1, t_2]$ with initial condition $u(t_1, \cd)$. Note that $\|u(t_1,\cd)\|_\infty \le C_0$ for the same global bound $C_0$. Since $\d_0$ corresponds to this $C_0$, following the same arguments, we see that the  comparison principle holds on $[t_1, t_2]$. Repeating the arguments 
establishes the result on the whole interval $[0, T]$.
 \qed
 
 When $f$ is semilinear, i.e., linear in $\g$, clearly under natural conditions $f$ satisfies the requirements in Theorem \ref{thm-comparisonSemilinear}.  We provide next a simple fully nonlinear example.
 
 \begin{eg}
 \label{eg-HJB}
 Let $\ol a > \ul a >0$ be two constants. Then 
 \bea
 \label{G}
 f(\g) := {1\over 2} \sup_{\ul a \le a \le \ol a} [a \g] = {1\over 2}[\ol a \g^+ - \ul a \g^-]
 \eea
  satisfies the requirements in Theorem \ref{thm-comparisonSemilinear}.
 \end{eg}
 
 
 \proof Let $\eta$ be a smooth symmetric density function with support  $(-1, 1)$. For any $\e>0$, introduce a smooth mollifier of $f$:
 \beaa
 f_\e(\g) := \int_{-1}^1 f(\g - \e x) \eta(x) dx = {1\over 2} \ul a  \g + {\ol a-\ul a\over 2}\int_{-1}^1 (\g-\e x)^+ \eta(x) dx.
  \eeaa
 It is clear that 
 \beaa
 |f_\e - f|\le \big[{\ol a\over 2}\int_{-1}^1 |x|\eta(x)dx\big]~ \e =: c\e.
 \eeaa
  We next consider the Legendre conjugate of $f_\e$:
 \beaa
  h_\e(a) := \sup_{\g \in \dbR} [{1\over 2} a\g - f_\e(\g)], ~a \in [\ul a, \ol a].
\eeaa
By straightforward calculation, we have $h_\e(a) = \infty$ when $a \notin [\ul a, \ol a]$, and 
\beaa
h_\e(a) = {\e\over 2}[\ol a - \ul a] \int_{-1}^{\Phi^{-1}({a-\ul a\over \ol a-\ul a})} x \eta(x) dx,\q a\in [\ul a, \ol a],
\eeaa
where $\Phi(x) := \int_{-1}^x \eta(y) dy$, $x\in [-1,1]$. Note that $\dis f_\e(\g) = \sup_{\ul a \le a \le \ol a} [{1\over 2} a\g - h_\e(a)]$. Then $\ol f_\e := f_{\e\over 2c} + {\e\over 2}$ and $\ol f_\e := f_{\e\over 2c} -{\e\over 2}$ are the desired approximations.
\qed
 
 \begin{rem}
\label{rem-comarisonGeneral}
{\rm (i) As pointed out in Remark \ref{rem-visPDE}, for general $g=g(t,x,y,z)$, the transformation is local and the  $\d$ in Theorem \ref{thm-utov}  depends on $\|\pa^2_{xx} u\|_\infty$.  Then the connection between RPDE \eqref{RPDE} and PDE \eqref{PDE} exists only for local classical solutions, but is not clear for viscosity solutions. Since our current approach relies heavily on the PDE, we have difficulty to extend Theorem \ref{thm-equiv1} to the general case, even just in local sense. We  investigate this challenging problem by exploring more approaches in our future research. 

(ii) When $f$ is of first order, namely $\si_f = 0$ in \eqref{fIsaacs}, then \eqref{FIsaacs}  becomes:
\begin{align}
\label{FIsaacs1st}
& F(t,x,y,z,\g) = \sup_{e_1\in E_1} \inf_{e_2\in E_2} \Big[ \hat b_f(t,x,e) z + F_0\big(t,x,y,e\big)\Big], \\
&\mbox{where}\quad \hat b_f(t,x,e) :=  { b_f(t, X_t, e)\over \pa_x X_t}, \nonumber \\ 
&F_0(t,x,y,e) :=     {1\over \pa_y Y_t} f_0\big(t,X_t, Y_t,   0, e\big)+\hat b_f(t, x, e) {\pa_x Y_t\over \pa_y Y_t}.\nonumber
\end{align}
Under Assumption \ref{assum-f1}, $F_0$ is uniformly Lipschitz continuous in $y$, and thus the main difficulty mentioned in Remark \ref{rem-blowup} (i) does not exist here. Then, following similar arguments as in this subsection, we can easily show that the results of Theorems \ref{thm-semilinear-classical} and \ref{thm-comparisonSemilinear} still hold true if we replace the uniform nondegeneracy condition $\si_f \ge c_0>0$ with $\si_f = 0$.
}
 \end{rem}

{\color{black}

 \subsection{The case that $g$ is linear}
\label{sect-linear}
In this subsection we study the special case that $g$ is linear in $(y,z)$ (by abusing the notation $g_0$):
\bea
\label{glinear}
g(t,x,y,z) = \si(t,x) z + h(t,x) y + g_0(t,x).
\eea
We remark that rigorously speaking this case does not satisfy Assumption \ref{assum-g2}, because $g_0(t,x,y) := h(t,x)y + g_0(t,x)$ is not bounded in $y$. However, similar to the situation in Lemma \ref{lem-RDElinear}, the linear structure allows us to extend all the results in Section \ref{sect-semilinear} to this case.

First, for the $X$ given by \eqref{characteristic1}, we have
\begin{align}
\label{linearY}
Y_t(x,y) =  e^{H_t(x)} \Big[y+ \int_0^t   e^{-H_s(x)} g_0(s, X_s(x)) d\o_s\Big], 
\end{align}
where $H_t (x):=  \int_0^t h(s, X_s(x))d\o_s$.
By straightforward calculation, we see that \eqref{F1} becomes, with omitting $(x, y)$ in $(X, Y, H)$,
\bea
\label{Flinearg}
&&F(t,x,y,z,\g) := e^{-H_t} f\Big(t, X_t, Y_t, {e^{H_t}\over \pa_x X_t} z + {\pa_x Y_t\over \pa_x X_t}, \\
&&\qq\qq {e^{H_t} \over (\pa_x X_t)^2}\big[\g+  [\pa_x H_t -  {\pa^2_{xx} X_t\over \pa_xX_t}] z+ e^{-H_t} \pa^2_{xx} Y_t - { \pa^2_{xx} X_t\over \pa_xX_t}\big]\Big) \nonumber.
\eea
 We now provide some sufficient conditions for the existence of classical solutions to  PDE \eqref{PDE}-\eqref{Flinearg}.  
 \begin{thm}
 \label{thm-linearg-classical}
 Let all the conditions in Theorem \ref{thm-semilinear-classical} hold
  and let $g$ take the form  \eqref{glinear}.  Then  PDE \eqref{PDE}-\eqref{Flinearg} has a classical solution $v\in C^{k,0}_{\a,\b}(\dbR_T)$. Consequently, RPDE \eqref{RPDE}-\eqref{glinear} has a classical solution $u\in C^k_{\a,\b}(\dbR_T)$.
 \end{thm}
 \proof  As in Theorem \ref{thm-semilinear-classical}, we shall only prove the convex case. By the regularity of $f$, $g$ it is clear that $F$ is smooth. Note that, by omitting the variables,  
 \begin{align*}
 \pa_\g F &= {1 \over (\pa_x X)^2} \pa_\g f,\q  \pa^2_{\g\g} F = { e^{H_t}\over (\pa_x X)^4}\pa^2_{\g\g} f, \\
\pa_z F &= {1\over \pa_x X} \pa_z f + [{\pa_x H \over (\pa_x X)^2} - {\pa^2_{xx} X \over (\pa_x X)^3}] \pa_\g f,\\
 \pa_y F 
 &=\pa_y f  + \pa_z f {\pa_{x} H\over \pa_x X} +\pa_\g f {(\pa_x H)^2 +\pa^2_{xx} H\over (\pa_x X)^2} .  
 \end{align*}
 Then one can easily verify that $F$ satisfies all the conditions in Theorem \ref{thm-Safonov}, thus we obtain $v\in C^{k,0}_{\a,\b}(\dbR_T)$. Finally the existence of 
 the corresponding function $u\in C^k_{\a,\b}(\dbR_T)$ follows from Theorem \ref{T:EU2RPDE}.
 \qed

We now assume further  that $f$ is also linear, i.e.,
\bea
\label{flinear}
 f(t,x,y,z,\g) = a(t,x) \g + b(t,x) z + c(t,x) y + f_0(t,x).
\eea
This case is well understood in the literature. By straightforward calculation,
\begin{align}
\label{linearF}
F(t,x,y,z,\g) =  A(t,x) \g + B(t,x) z + C(t,x) y + F_0(t,x),
\end{align}
where, for the $H$ defined by \eqref{linearY} and again omitting the variable $x$,
\begin{align}
\label{linearABCF0}
 A(t,x) :={a(t,X_t) \over (\pa_x X_t)^2}, \hskip 8cm\quad\nonumber\\
 B(t,x) := a(t,X_t) \big[  {2\pa_x H_t  \over  (\pa_xX_t)^2} - {\pa^2_{xx} X_t\over  (\pa_xX_t)^3} \big] +{b(t,X_t) \over \pa_x X_t},\hskip3cm ~ \nonumber\\
 C(t,x) := a(t,X_t) {(\pa_x H_t)^2 + \pa^2_{xx} H_t\over (\pa_x X_t)^2}  + \big[{b(t,X_t)\over \pa_x X_t} - {a(t,X_t) \pa^2_{xx} X_t \over (\pa_x X_t)^3}\big]\pa_x H_t\nonumber\\
\qq\qq\q+ c(t,X_t),\hskip7cm\\
 F_0(t,x) :=\Big[ {a(t,X_t)\over (\pa_x X_t)^2} [(\pa_x H_t)^2 + \pa^2_{xx} H_t\big] \hskip5cm~\nonumber \\
\q+\big[{b(t,X_t)\over \pa_x X_t} - {a(t,X_t) \pa^2_{xx} X_t \over (\pa_x X_t)^3}\big] \pa_x H_t + c(t,X_t)\Big] \int_0^t e^{-H_s} g_0(s, X_s) d\o_s ~~\nonumber \\
\q+\Big[2{a(t,X_t)\over (\pa_x X_t)^2}\pa_x H_t +{b(t,X_t)\over \pa_x X_t} - {a(t,X_t) \pa^2_{xx} X_t \over (\pa_x X_t)^3}\Big]\int_0^t \pa_x ({g_0(s, X_s)\over e^{H_s}})\, d\o_s\ms\nonumber \\
\q+ {a(t,X_t)\over (\pa_x X_t)^2}   \int_0^t \pa^2_{xx} (e^{-H_s} g_0(s, X_s)) \,d\o_s  + f_0(t, X_t) e^{-H_t}.\hskip2cm~~\nonumber
\end{align}
Thus PDE \eqref{PDE} is  linear  and  we have the representation formula
\begin{equation}
\label{linear-rep2}
\begin{split}
v(t,x) =  \dbE\Big[e^{\int_0^t C(t-r, \cX^{t,x}_r) dr} u_0( \cX^{t,x}_t) + \int_0^t e^{\int_0^s C(t-r, \cX^{t,x}_r) dr} F_0(t-s, \cX^{t,x}_s) ds\Big],
\end{split}
\end{equation}
where, for fixed $(t,x) \in \dbR_T$ and for a Brownian motion $\bar B$, 
\begin{equation*}
\cX^{t,x}_s = x + \int_0^s \sqrt{2A(t-r, \cX^{t,x}_r)} dr +\int_0^s B(t-r, \cX^{t,x}_r) d\bar B_s,~ 0\le s\le t.
\end{equation*}

}

\section{Concluding Remarks}
\label{sect-summary}
\setcounter{equation}{0}

In this paper, we established the viscosity theory for general second order fully nonlinear parabolic  SPDEs and path-dependent PDEs through a unified framework based on rough path analysis. 
We allow the diffusion coefficient $g$ to be a general nonlinear first order differential operator, i.e., $g=g(t,\o, x,u, u_x)$. Assuming that  $g$ is smooth enough, we obtained 
the following general results:   1) {\bf Consistency} (i.e., if a candidate solution is 
smooth, 
viscosity solutions are equivalent to classical solutions); 2)  {\bf Equivalence}  (between the definitions using  test functions and by semi-jets, resp.);
3) {\bf Stability}; and 4) {\bf  Partial comparison principle}.

Since the generality of the diffusion coefficient $g$ corresponds to
 a highly convoluted system of characteristics, the complete wellposedness result, especially  global existence and the comparison principle for viscosity solutions, is rather challenging. We thus investigated several important cases,  summarized as follows. 

\ms
\no{\bf A.}~ Both $f$ and $g$ are linear. In this case, we have a representation formula, which is of course already well understood in the literature.

\ms

\no{\bf B.}~   $g$ is linear. In this case, we  proved 
\begin{itemize}
\item {\it global existence of classical solutions} when $f$ and $u_0$ are smooth,

\item {\it global existence and comparison principle} of viscosity solutions when $f$ and $u_0$ are less smooth.
\end{itemize}

\no{\bf C.}~  $g$ is semi-linear. In this case,  we  established
\begin{itemize}
\item {\it global equivalence} of RPDE and associated PDE,  in both the classical and the  viscosity sense,

\item {\it local existence of classical solutions} when $f$ and $u_0$ are smooth, 

\item {\it global existence and comparison principle} of viscosity solutions when $f$ and $u_0$ are less smooth but $f$ is convex in $\pa^2_{xx} u$.
\end{itemize}

\no{\bf D.} For  the general fully nonlinear case,  we proved
\begin{itemize}
\item {\it local equivalence} of RPDE and associated PDE in classical sense (but not in viscosity sense),

\item {\it local existence of classical solutions} when $f$ and $u_0$ are smooth, where the time interval depends on  the regularity of $u_0$ and $f$.
\end{itemize}

We should note that although we  finally found a framework on which the viscosity theory can be carried out in a fairly general 
way, some challenging issues remain, especially in the fully nonlinear case, which calls for further investigations.
E.g., when $g$ is semilinear, the current approach relies heavily on the (local) existence of classical solutions of certain mollified PDEs, for which we have to assume that $f$ is {\it uniformly non-degenerate}. It is
 desirable to remove both the convexity and 
the nondegeneracy constraints on $f$.   Also, as we pointed out in Remark \ref{rem-visPDE}, in the general fully nonlinear case the equivalence between RPDE and PDE in the 
viscosity sense is by no means clear. We often have to investigate the RPDE directly. Consequently, a direct approach for the comparison principle for  RPDE \eqref{RPDE}, which is currently lacking, would help greatly. In particular, it would be interesting to explore the possibility of extending  the  doubling variable arguments to this situation (cf.~Lions \& Souganidis \cite{LS4} for this approach in the case $g= g(\pa_x u)$).

 Finally, it would also be interesting to investigate the alternative approach by using rough path approximations  as  in Caruana,  Friz, \& Oberhauser \cite{CFO} and many other papers mentioned in the introduction. I.e., instead of transforming the RPDE into a PDE, one considers the  approximating PDE
\begin{equation}
\label{PDEn}
\begin{split}
d u^n(t,x) &= f (t, x, u^n, \pa_x u^n, \pa^2_{xx} u^n) dt + g(t,x, u^n, \pa_x u^n) d \o^n_t,\\
u^n(0,x) &= u_0(x),
\end{split}
\end{equation}
where $\o^n$ is a smooth function of $t$ which converges to $\o$ under a rough path norm. The key issue is then to study the convergence of $u^n$. We hope that we will be able to investigate 
 some of these issues in our future publications.

\section{Appendix: Proofs for some results in  Section \ref{sect-rough}}
\label{Appendix}
\setcounter{equation}{0}

For notational simplicity, 
we may often write  $u_t(x) := u(t,x)$.

\no{\bf Proof of Lemma \ref{lem-commute}.}  First, we prove \eqref{commute1}. For $u \in C^{2, loc}_{\a,\b}(\dbR_T^d)$, 
\beaa
u_{s, t}(x)  = \pa_\o u_s(x) \o_{s,t} + R^{1,u}_{s,t}(x).
\eeaa
Differentiating both sides with respect to $x$ yields
\beaa
\pa_x u_{s, t}( x) =\pa_x \pa_\o u_s(x) \o_{s,t} +\pa_x R^{1,u}_{s,t}(x).
\eeaa
By Definition \ref{defn-Cxloc} (iii),  $\pa_x \pa_\o u(\cd, x) \in C^{0, loc}_{\a,\b}([0, T])$ and $[\pa_x R^{1,u}(x)]_{\a(1+\b)} <\infty$. Thus \eqref{commute1} follows directly from \eqref{C1}. 
Next, we  prove \eqref{Fubini}.   Assume $d=1$ for simplicity.  Put $v_t(x) := \int_0^x u_t(y) dy$.  Then
\beaa
v_{s,t}(x)  = \int_0^x u_{s, t}(y) dy = \int_0^x \pa_\o u_s(y)dy ~\o_{s, t} + \int_0^x R^{1,u}_{s,t}(y) dy.
\eeaa
By  continuity of $x\mapsto \pa_\o u(\cd, x)\in C^0_{\a,\b}([0, T])$ and $x\mapsto  u(\cd, x)\in C^1_{\a,\b}([0, T])$, it is clear that $ \int_0^x \pa_\o u(\cd, y)dy \in C^0_{\a,\b}([0, T])$ and $[\int_0^x R^{1,u}(y)dy]_{\a(1+\b)} <\infty$. Then $\pa_\o v_t(x) =  \int_0^x \pa_\o u_t(y)dy$. For the partition 
$t_i = i 2^{-n} t$, $i=0$, $\ldots$, $2^n$, 
\begin{align*}
&\int_0^t \int_0^x u_s(y) \,dy \,d\o_s = \int_0^t v_s(x) \,d\o_s \\ &= \lim_{n\to \infty} \sum_{i=0}^{2^n-1} \big[v_{t_i}( x) \o_{t_i, t_{i+1}} + \pa_\o v_{t_i}(x) \ul\o_{t_i, t_{i+1}}\big]\\
&= \lim_{n\to \infty} \int_0^x \sum_{i=0}^{2^n-1} \big[u_{t_i}(y) \o_{t_i, t_{i+1}} + \pa_\o u_{t_i}(y) \ul\o_{t_i, t_{i+1}}\big] \,dy.
\end{align*}
By  standard estimates, see, e.g., Keller \& Zhang \cite[Lemma 2.5]{KZ}, 
\beaa
\Big|\int_{t_i}^{t_{i+1}} u_s(y) d\o_s - u_{t_i}(y) \o_{t_i, t_{i+1}} - \pa_\o u_{t_i}(y) \ul\o_{t_i, t_{i+1}}\Big| \le {C\|u(\cd, y)\|_1 \over 2^{\a(2+\b)n}}.
\eeaa
Since $O$ is bounded, by  continuity of $u$ again,  $\sup_{u\in O} \|u(\cd, y)\|_1<\infty$. Then, by  \eqref{ab} and the dominated convergence theorem, we obtain immediately 
\begin{align*}
\int_0^t \int_0^x u_s(y) \,dy \,d\o_s &=\lim_{n\to \infty} \int_0^x 
 \sum_{i=0}^{2^n-1}
\int_{t_i}^{t_{i+1}} u_s(y)\,d\o_s  \,dy \\
&= \int_0^x \int_0^t u_s(y) \,d\o_s \,dy.
\end{align*}
This yields \eqref{Fubini} for general $s$, $t$ and $O$ immediately.

Finally, we  prove \eqref{commute2}. Let $u \in C^{3, loc}_{\a,\b}( \dbR_T^d)$.  By Lemma \ref{lem-Ito}, 
\beaa
u_{s,t}(x) = \int_s^t \pa^\o_t u_r(x) \,dr +  \int_s^t \Big[\pa_\o u_r(0) +  \int_0^x \pa_x\pa_\o u_r(y) \,dy\Big]\,d\o_r. 
\eeaa
Differentiating both sides with respect to $x$ together with applying \eqref{Fubini} on the last term above, 
we obtain
\beaa
\pa_x u_{s,t}(x) = \int_s^t \pa_x \pa^\o_t u_r(x) \,dr +  \int_s^t \pa_x\pa_\o u_r(x) \,d\o_r. 
\eeaa
Then  \eqref{commute2} follows from Lemma \ref{lem-Ito}.
\qed

\ms

\no{\bf Proof of Lemma \ref{lem-Taylor}.} We proceed in two steps.

{\it Step 1.} Let $h=0$. For notational simplicity, we omit the variable $x$ in this step. We shall prove by induction that 
\begin{align}
\label{Taylor1}
u_t = \sum_{\|\nu\|\le k} \cD_\nu u_s I^\nu_{s, t} + R^{k,u}_{s,t} ~\mbox{and}~ |R^{k,u}_{s, t}|\le C\|u(\cd,x)\|_k 
\,\vert t-s\vert^{\a(k+\b)}.
\end{align}

When $k=0$, $1$, $2$, \eqref{Taylor1} follows directly from the definitions of the derivatives.  
 Let $m\ge 2$ and assume  that \eqref{Taylor1} holds for  all $k\le m$. 
Let now  $u \in C^{m+1, loc}_{\a,\b}([0, T]\times \dbR)$. By \eqref{Ito},  
\beaa
u_{s,t} = \int_s^t \pa^\o_t u_r dr + \int_s^t \pa_\o u_r d\o_r.
\eeaa
Note that $\pa^\o_t u \in  C^{m-1, loc}_{\a,\b}([0, T]\times \dbR)$. By the induction hypothesis
\begin{align*}
&\pa^\o_t u_r  =  \sum_{\|\nu\|\le m-1} \!\!\! \cD_\nu \pa^\o_t u_s \cI^\nu_{ s,r}  + R^{m-1, \pa^\o_t u}_{s,r}, \\
&\mbox{with}~|R^{m-1, \pa^\o_t u}_{s,r}|\le C\|\pa^\o_t u(\cd, x)\|_{m-1} \vert r-s\vert^{\a(m-1+\b)}.
\end{align*}
Then, since $2\a < 1$,
\begin{align*}
&\Big| \int_s^t \pa^\o_t u_r\,dr - \sum_{\|\nu\|\le m-1}  \cD_\nu \pa^\o_t u_s \int_s^t\cI^\nu_{s,r}\,dr \Big| = \Big| \int_s^t R^{m-1, \pa^\o_t u}_{s,r} \,dr \Big|\\
&\le  C\|\pa^\o_t u(\cd, x)\|_{m-1}\,
 \Big| \int_s^t (r- s)^{\a(m-1+\b)} \,dr \Big|\\ 
 &\le C\|u(\cd, x)\|_{m+1} \,
 \vert t-s\vert^{\a(m+1+\b)}.
 \end{align*}
 Thus it  remains to show that
 \begin{equation}
\label{Taylor2}
\begin{split}
&\Big|\int_s^t \pa_\o u_r \,d\o_r -  \sum_{\|\nu\|\le m} \cD_\nu\pa_\o u_s \int_s^t \cI^\nu_{s,r} \,d\o_r \Big|  \\
&\le C\|u(\cd, x)\|_{m+1}
 \vert t-s\vert^{\a(m+1+\b)}.
\end{split}
\end{equation}
 Note that $v := \pa_\o u \in C^{m, loc}_{\a,\b}([0, T]\times \dbR)$.  Put 
 $U_{s,t} := \sum_{\|\nu\|\le m} \cD_\nu v_s \int_s^t \cI^\nu_{s,r} \,d\o_r$.  For $s < r<t$, 
\begin{align*}
& U_{s,r} +U_{r,t} - U_{s,t}  =   \sum_{\|\nu\|\le m} \cD_\nu v_r \int_r^t \cI^\nu_{r,l}\, d\o_l  ~ - ~
 \sum_{\|\nu\|\le m} \cD_\nu v_s \int_r^t \cI^\nu_{s,l} \,d\o_l \\
&= \sum_{\|\nu\|\le m}  \Big[\sum_{\|\nu'\|\le m-\|\nu\|} \cD_{\nu'} \cD_\nu v_s  \cI^{\nu'}_{s, r} + R^{m-\|\nu\|, \cD_{\nu} v}_{s, r}\Big] \int_r^t \cI^\nu_{r,l}\, d\o_l  \\ 
&\qq\qq -  \sum_{\|\nu\|\le m} \cD_\nu v_s \int_r^t \cI^\nu_{s,l} d\o_l \\
&=  \sum_{\|\nu\|\le m}  \cD_\nu v_s \int_r^t  \Big[\sum_{(\nu', \tilde \nu)= \nu} \cI^{\nu'}_{s, r}\,  
\cI^{\tilde \nu}_{r, l} - \cI^\nu_{s,l}\Big] d \o_l +  \sum_{\|\nu\|\le m} R^{m-\|\nu\|, \cD_{\nu} v}_{s, r} \int_r^t 
\cI^\nu_{r,l}\, d\o_l.
\end{align*}
By induction, one can verify straightforwardly that
\beaa
\sum_{(\nu', \tilde \nu)= \nu} \cI^{\nu'}_{s, r}\,  \cI^{\tilde \nu}_{r, l} = \cI^\nu_{s,l}\q\mbox{and}\q 
\Big|\int_r^t \cI^\nu_{r,l} \,d\o_l \Big| \le C(t-r)^{(\|\nu\|+1)\a}.
\eeaa
We remark that the former identity here is actually Chen's relation. Then, since $\cD_\nu v \in C^{m-\|\nu\|, loc}_{\a,\b}([0, T]\times \dbR)$, by our induction assumption, 
\begin{align*}
|R^{m-\|\nu\|, \cD_{\nu} v}_{s, r}| &\le C\|\cD_\nu v(\cd, x)\|_{m-\|\nu\|}(r-s)^{\a(m-\|\nu\| + \b)} \\
&\le C\|u(\cd, x)\|_{m+1} (r-s)^{\a(m-\|\nu\| + \b)}.
\end{align*}
Therefore,
\beaa
\big|U_{s,r} +U_{r,t} - U_{s,t}\big| &\le& C\|u(\cd, x)\|_{m+1} (r-s)^{\a(m-\|\nu\| + \b)}(t-r)^{(\|\nu\|+1)\a} \\
&\le& C\|u(\cd, x)\|_{m+1} (t-s)^{\a(m+1+ \b)}.
 \eeaa
 Applying the Sewing Lemma, see, e.g., Friz \& Hairer \cite[Lemma 4.2]{FH},  yields \eqref{Taylor2}.
  Hence \eqref{Taylor1}  follows.

{\it Step 2.}  We now prove the general case. By  standard Taylor expansion in $x$ and by Step~1, 
\begin{align*}
&u(t+\d, x+h)  = [u(t+\d, x+h) - u(t, x+h)] + u(t,x+h) \\
&= \sum_{1\le \|\nu\|\le k} \cD_\nu u(t, x+h)\, \cI^\nu_{t, t+\d} +  R^{k, u(\cd, x+h)}_{t, t+\d}  + \sum_{m=0}^k \frac{1}{m!}\pa^m_x u(t,x) h^m \\ &\qq\qq+ O(|h|^{k+\b}) \\
&= \sum_{1\le \|\nu\|\le k} \Big[\sum_{m=0}^{k-\|\nu\|} 
 \frac{1}{m!}
\pa^m_x \cD_\nu u(t, x) h^m 
+ O(|h|^{k-\|\nu\|+\b}) \Big]\,
 \cI^\nu_{t, t+\d} +  R^{k, u(\cd, x+h)}_{t, t+\d}\\
&\qq  + \sum_{m=0}^k 
 \frac{1}{m!}
\pa^m_x u(t,x) h^m + O(|h|^{k+\b}) \\
&=\!\! \sum_{m+\|\nu\|\le k} 
 \frac{1}{m!}
\cD_\nu \pa^m_x u(t, x) h^m\, \cI^\nu_{t, t+\d} 
 +\!\! \sum_{1\le \|\nu\|\le k}\!\! O(|h|^{k-\|\nu\|+\b})\,\cI^\nu_{t, t+\d}  \\
 &\qq\qq+  R^{k, u(\cd, x+h)}_{t, t+\d}+ O(|h|^{k+\b}),
\end{align*}
where the last equality follows from Lemma \ref{lem-commute}. Note that 
$|\cI^\nu_{t, t+\d}| \le C|\d|^{\a \|\nu\|}$ and, by Step~1,  $|R^{k, u(\cd, x+h)}_{t, t+\d}|\le C\|u(\cd, x+h)\|_k\,
 |\d|^{\a(k+\b)}$. Since $\|u(\cd, x)\|_k$ is continuous in $x$ and thus locally bounded, we obtain \eqref{Taylor} immediately. 
\qed

\ms
\no{\bf Proof of Lemma \ref{lem-RDE}.} The wellposedness for a classical solution $u\in C^1_{\alpha,\beta}([0,T])$ with the corresponding estimate is standard. 
 Automatically, by Lemma~\ref{lem-Ito}, $u\in C^2_{\alpha,\beta}([0,T])$ again with the corresponding estimate.
Note that $\pa_\o u_t = g(t, u_t), \pa^\o_t u_t = f(t, u_t)$. Applying Lemma \ref{lem-chain} repeatedly,  one can easily prove by induction that $u\in C^{k+2}_{\a,\b}([0, T])$ and that
\bea
\label{RDEest0}
\|u\|_{k+2} \le C(T, \|f\|_k, \|g\|_{k+1}, u_0).
\eea
 Moreover, with $\tilde u_t := u_t - u_0$, $\tilde f(t, x) := f(t, x+u_0)$, $\tilde g(t,x) := g(t, x+u_0)$, 
\beaa
 \tilde u_t =  \int_0^t \tilde f(s, \tilde u_s) \,ds + \int_0^t\tilde g(s, \tilde u_s) \,d\o_s.
\eeaa
 Thus, by  \eqref{RDEest0},  
$\|\tilde u\|_{k+2} \le C(T,  \|\tilde f\|_k, \|\tilde g\|_{k+1}, 0) = C(T,  \|f\|_k, \|g\|_{k+1})$.
 \qed

\ms

\no{\bf Proof of Lemma \ref{lem-RDEx}.}  First, we  assume that $u_0 \in C^{k+\b}_b(\dbR)$, namely $u_0$ and all its related derivatives are bounded.  Applying Lemma \ref{lem-RDE} we see that $u(\cd, x) \in C^{k+1}_{\a,\b}([0, T])$ and
\bea
\label{RDExest1}
\sup_{x\in \dbR} \|u(\cd, x)\|_{k+1} \le C(\|u_0\|_\infty, \|f\|_k, \|g\|_k)<\infty.
\eea
Recalling Definition \ref{defn-Cx} and in particular Remark \ref{rem-Cx}, we note that our space $C^k_{\a,\b}(\dbR_T^d)$ requires stronger conditions than those in Keller \& Zhang \cite{KZ}.  Thus one may follow the arguments of  \cite[Theorem 6.1]{KZ} to obtain the equation \eqref{paxuv}.   Hence
\begin{align*}
&\pa_x R^{1,u}_{s,t}(x) = \pa_x \big[u_{s,t}(x) - g(s,x, u_s(x)) \o_{s,t} \big]\nonumber\\
&= \pa_x u_{s,t}(x) - \big[\pa_x  g(s, x, u_s(x)) + \pa_y g(s, x, u_s(x))\pa_x u_s(x)\big]\o_{s,t}\nonumber\\
&= \int_s^t [\pa_x f + \pa_y f \pa_x u_r(x)](r, x, u_r(x)) dr \\
&\qq +  \int_s^t [\pa_x g + \pa_y g \pa_x u_r(x)](r, x, u_r(x)) d\o_r\nonumber\\
&\qq - \big[\pa_x  g(s, x, u_s(x)) + \pa_y g(s, x, u_s(x))\pa_x u_s(x)\big]\o_{s,t}.
\end{align*}
Applying Lemma  \ref{lem-RDElinear} on \eqref{paxuv} yields a representation formula for $\pa_x u$ and 
\bea
\label{RDExest2}
\sup_{x\in \dbR} \|\pa_x u(\cd, x)\|_{k} \le C(\|\pa_x u_0\|_\infty, \|f\|_k, \|g\|_k)<\infty,
\eea
which, together with \eqref{RDExest1}, implies 
further 
 that
\bea
\label{RDExest3}
\sup_{x\in\dbR}[\pa_x R^{1,u}(x)]_{\a(1+\b)} \le C(\|u_0\|_\infty, \|\pa_x u_0\|_\infty, \|f\|_k, \|g\|_k).
\eea
 Together with the representation for $\pa_x u$ and the arguments in Lemma~\ref{lem-chain}, 
\begin{equation}
\label{uxx}
\begin{split}
& \pa^2_{xx} u_t(x) = \pa^2_{xx} u_0(x) \\ 
&\qquad+  \int_0^t \big[\pa^2_{xx} f + 2\pa^2_{xy} f \pa_x u_s(x) + \pa_y f \pa^2_{xx} u_s(x)\big](s, x, u_s(x)) ds \\
& \qquad +  \int_0^t \big[\pa^2_{xx} f + 2\pa^2_{xy} f \pa_x u_s(x) + \pa_y f \pa^2_{xx} u_s(x)\big](s, x, u_s(x))  d\o_s,
\end{split}
\end{equation}
follows.
Then one can easily see that 
 $\| u\|_2<\infty$ and thus
 $u\in C^1_{\a, \b}(\dbR_T)$ in the sense of Definition \ref{defn-Cx}
 (cf.~Remark \ref{rem-Cx}). Repeating the arguments
(up to order $k+1$ to ensure H\"older continuity with respect to $x$
of $\pa^k_x u$, etc.), one can  
 show that 
 $\| u\|_{k+1}<\infty$, which implies that $u\in C^k_{\a, \b}(\dbR_T)$.

Next, if $u_0$ is not bounded but all its related derivatives are bounded. Put $\tilde u := u - u_0$, $\tilde \f(t,x,y) := \f(t,x, y+u_0(x))$ for $\f = f$, $g$. One can easily see that 
$\tilde f\in C^{k+1}_{\a,\b}( \dbR_T^2)$ with $\|\tilde f\|_k$ being dominated by $\|f\|_{k+1}$ and the derivatives of $u_0$, similarly for $\tilde g$.  Note that $\tilde u$ satisfies RDE \eqref{RDEx} with coefficients $(\tilde f, \tilde g)$ and initial condition $0$. Thus 
$\tilde u\in C^k_{\a, \b}(\dbR_T)$. 

Finally, if $u_0\in C^{k+\b}(\dbR)$, let $\iota\in C^\infty(\dbR)$ such that $\iota(x) = 1$ when $|x|\le 1$ and $\iota(x)=0$ when $|x|\ge 2$. Let $u^\e\in C^k_{\a, \b}(\dbR_T)$ be the solution to RDE \eqref{RDEx} with coefficients $(f, g)$ and initial condition
 $u^\e_0(x) = u_0( \iota(\e x)\,x)$. Note that $u^\e_0(x) = u_0(x)$ and  hence $u^\e_t(x) = u_t(x)$ whenever $|x|\le {1/ \e}$. Since $\e>0$ is arbitrary, we see that $u\in C^{k,loc}_{\a, \b}(\dbR_T)$.
\qed

\end{document}